\documentclass[english]{article}
\usepackage[T1]{fontenc}
\usepackage[latin9]{inputenc}
\usepackage{float}
\usepackage{amsmath}
\usepackage{amsthm}
\usepackage{amssymb}
\usepackage{graphicx}

\makeatletter
 \theoremstyle{definition}
  \newtheorem{example}{\protect\examplename}
\theoremstyle{plain}
\newtheorem{thm}{\protect\theoremname}
  \theoremstyle{remark}
  \newtheorem{rem}{\protect\remarkname}
  \theoremstyle{plain}
  \newtheorem{prop}{\protect\propositionname}

\@ifundefined{showcaptionsetup}{}{%
 \PassOptionsToPackage{caption=false}{subfig}}
\usepackage{subfig}
\makeatother

\usepackage{babel}
  \providecommand{\examplename}{Example}
  \providecommand{\propositionname}{Proposition}
  \providecommand{\remarkname}{Remark}
\providecommand{\theoremname}{Theorem}

\begin{document}

\title{Exact Model Reduction by a Slow-Fast Decomposition\\
of Nonlinear Mechanical Systems}

\author{George Haller\thanks{Corresponding author. Email: georgehaller@ethz.ch}
and Sten Ponsioen}
\maketitle
\begin{center}
Institute for Mechanical Systems, ETH Zürich \\
Leonhardstrasse 21, 8092 Zürich, Switzerland\\
\par\end{center}
\begin{abstract}
We derive conditions under which a general nonlinear mechanical system
can be exactly reduced to a lower-dimensional model that involves
only the most flexible degrees of freedom. This Slow-Fast Decomposition
(SFD) enslaves exponentially fast the stiff degrees of freedom to
the flexible ones as all oscillations converge to the reduced model
defined on a slow manifold. We obtain an expression for the domain
boundary beyond which the reduced model ceases to be relevant due
to a generic loss of stability of the slow manifold. We also find
that near equilibria, the SFD gives a mathematical justification for
two modal-reduction methods used in structural dynamics: static condensation
and modal derivatives. These formal reduction procedures, however,
are also found to return incorrect results when the SFD conditions
do not hold. We illustrate all these results on mechanical examples.
\end{abstract}

\section{Introduction}

While often hoped otherwise, a typical multi-degree-of-freedom mechanical
system cannot necessarily be reduced to a lower-dimensional model.
There is often a good reason why the original model involves several
degrees of freedom, all of which are essential to reproduce the dynamics
at the required level of accuracy.

\begin{figure}[h]
\begin{centering}
\includegraphics[width=0.4\textwidth]{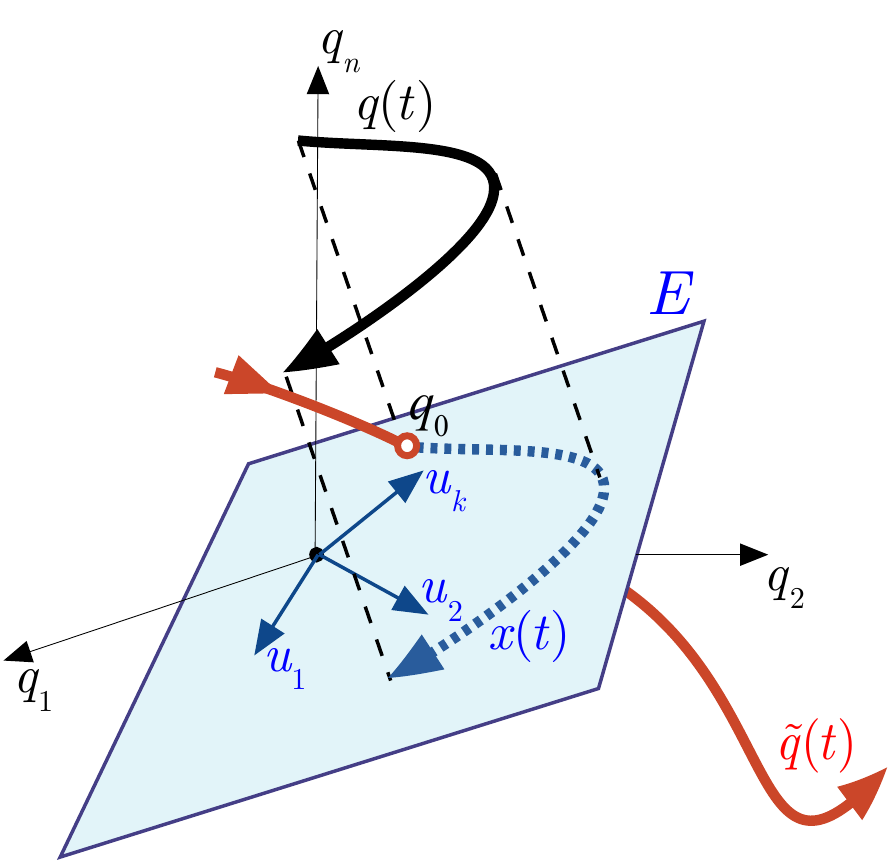}
\par\end{centering}
\caption{Model reduction by projection of a full trajectory $q(t)$ onto a
$k$-dimensional subspace $E$, typically spanned by a few eigenvectors
$u_{1},\ldots,u_{k}$ of the linearized system at the origin. The
model trajectory $x(t)$ starting from a point $q_{0}\in E$ is constrained
to lie in $E$, but the full trajectory $\tilde{q}(t)$ starting from
$q_{0}$ will generally leave the plane $E$. }
\label{fig:projection}

\end{figure}

For any multi-degree-of-freedom system, projections to various linear
subspaces are nevertheless routinely employed for model reduction
purposes (see Besselink et al. \cite{besselink13} for a review of
techniques in structural vibrations, Benner et al. \cite{benner15}
for a more general survey). Most often, however, the accuracy or even
the fundamental validity of these procedures is a priori unknown.
The main reason is that distinguished subspaces identified from linearization
or other considerations are generally not invariant under the nonlinear
dynamics. As a consequence, trajectories of the full system do not
follow those of a projection-based model, as shown in Fig. \ref{fig:projection}.

A model reduction principle can be justified in a strict mathematical
sense if the reduced model is defined on an invariant set of the full
nonlinear system, and hence model trajectories are actual trajectories
of the full system. In addition, the invariant set carrying the model
dynamics should be robust and attracting for the reduced model to
be of relevance for typical trajectories. While numerical or perturbative
approximations to such a set will at best be approximately invariant,
the attractivity of the actual invariant manifold is expected to keep
the impact of non-invariance small, driving trajectories toward the
actual invariant set.

\begin{figure}[h]
\centering{}\includegraphics[width=1\textwidth]{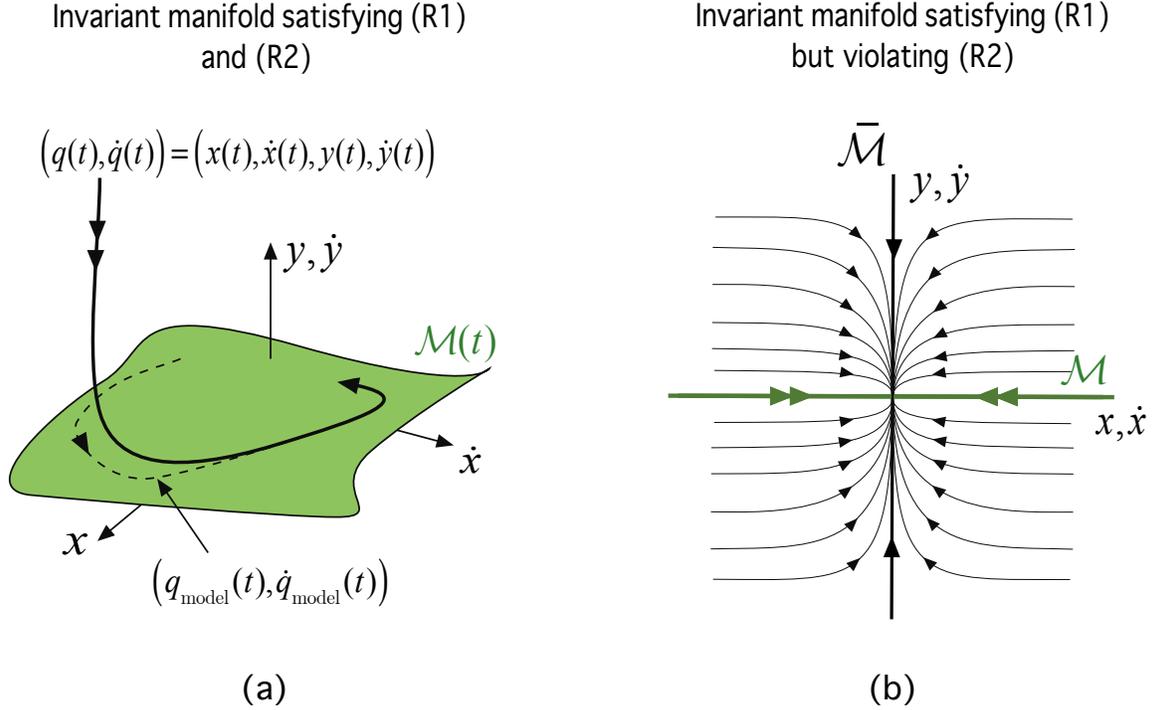}\caption{(a) Illustration of the geometry of requirements (R1) and (R2) for
model reduction in a mechanical system with generalized coordinates
$q$ and associated velocities $\dot{q}$. The reduced model depends
only on a smaller group of degrees of freedom, described by the position
vector $x$ and the corresponding velocity vector $\dot{x}$. The
remaining degrees of freedom are characterized by the positions $y$
and velocities $\dot{y}$. (b) An attracting and persistent invariant
manifold $\mathcal{M}$ that does not provide a faithful reduced-order
model for the full system dynamics. \label{fig:(R1)-(R2)}}
\end{figure}

Motivated by these considerations, we propose here two requirements
for mathematically justifiable and robust model reduction in a nonlinear,
non-autonomous mechanical system:
\begin{description}
\item [{(R1)}] There exists an attracting and persistent lower-dimensional
forward-invariant manifold $\mathcal{M}(t)$. Along the manifold $\mathcal{M}(t)$,
the modeled degrees of freedom (with generalized coordinates $y$
and velocities $\dot{y}$) are smooth functions of the modeling degrees
of freedom (with generalized coordinates $x$ and velocities $\dot{x}$)
. 
\item [{(R2)}] General trajectories approaching $\mathcal{M}(t)$ synchronize
with model trajectories at rates that are faster than typical rates
within $\mathcal{M}(t)$.
\end{description}
By the requirement (R1), the construction of a smooth, lower-dimensional
dynamical model should be equivalent to a reduction to a lower-dimensional
invariant manifold, as illustrated in Fig. \ref{fig:(R1)-(R2)}a.
The dynamics on this manifold, however, is only relevant for the full
system dynamics if nearby motions $q(t)$ approach model trajectories
on $\mathcal{M}(t)$, i.e., the manifold has a domain of attraction
foliated by stable manifolds of individual model trajectories. In
addition, we require $\mathcal{M}(t)$ to be persistent (robust under
small perturbations) since mechanical models have inherent parameter
uncertainties and approximations, and a model reduction should be
robust with respect to these. 

Requirement (R2) ensures that full system trajectories not only approach
the set of model trajectories in the phase space, but also synchronize
with specific model trajectories. Consider, for example, a linear,
two-degree of freedom mechanical system with an asymptotically stable
fixed point at the origin. The fast stable manifold $\mathcal{M}$
of this fixed point (cf. Fig. \ref{fig:(R1)-(R2)}b) is invariant,
attracting and persistent, even unique (see, e.g., Cabre et al. \cite{cabre03}).
Yet, the dynamics on $\mathcal{M}$ fails to act as a faithful reduced-order
model for the typical near-equilibrium dynamics. Indeed, the flow
on $\mathcal{M}$ predicts a fast decay rate that is unobservable
along typical trajectories on their way to the fixed point. This is
because general trajectories first approach the $(x,\dot{x})=(0,0)$
subspace, then creep towards the origin more slowly, synchronizing
with motions along this subspace, rather than with those in $\mathcal{M}$. 

In contrast, the invariant manifold $\mathcal{\bar{M}}$ in Fig. \ref{fig:(R1)-(R2)}b
satisfies both (R1) and (R2) and is indeed a good choice for model
reduction. Tis is ensured by a dichotomy of time scales created by
the gap in the real part of the spectrum of the eigenvalues of the
fixed point. The larger this gap, the more efficient the reduced-order
model in predicting typical\emph{ }system behavior. 

In a nonlinear system, one generally loses the local slow-fast dichotomy
of time scales that may arise near fixed points, such as the one in
Fig. \ref{fig:(R1)-(R2)}b. A global reduced-order model with the
properties (R1)-(R2) will, therefore, not exist unless the slow-fast
timescale difference created locally by the fixed point extends to
a larger domain of the phase space. In more mechanical terms, a global
model reduction is only feasible when the $x$ variables stay globally
stiffer (i.e., faster) than the $y$ variables. 

Such a global slow-fast partition of coordinates has been assumed
in several case studies of mechanical systems, such as an undamped
spring coupled to a pendulum (Georgiou and Schwartz \cite{georgiou96})
and its extensions to higher or even infinitely many dimensions (Georgiou
and Schwartz \cite{georgiou99} and Georgiou and Vakakis \cite{georgiou96-1}).
These studies tacitly assume the existence of a slow manifold without
specific consideration to its stability and robustness. Due to a lack
of normal hyperbolicity for the limiting slow manifold (critical manifold),
well-defined invariant slow manifolds do not actually exist in these
mechanical models. Recent results guarantee only near-invariant surfaces
under certain conditions (MacKay \cite{mackay04}, Kristianssen and
Wullf \cite{kristiansen16}). These surfaces, however, do not attract
trajectories from an open neighborhood of the phase space. As a result,
their relevance for model reduction is a priori unclear, as they violate
the requirement (R1).

As a further case study, a forced and stiff linear oscillator coupled
to a soft nonlinear oscillator was considered by Georgiou et al. \cite{georgiou98,georgiou99a}.
As the authors observe, the existence of an attracting, two-dimensional
slow manifold in these two studies follows from a globalized version
of the center manifold theorem (Carr \cite{carr82}) and from the
geometric singular perturbation formulation of Fenichel \cite{fenichel79},
respectively. These approaches are similar in spirit to the work we
describe here, but pertain to specific, low-dimensional, soft-stiff
mechanical models without targeting model reduction issues per se. 

Related work also includes that of Lubich \cite{lubich93}, who developed
a numerical scheme for mechanical systems with very stiff potential
forces. In this context, all degrees of freedom are equally fast and
hence no oscillatory mode can be enslaved to the rest via model reduction.
An exceptional slow manifold (which involves coordinates from all
degrees of freedom) becomes attracting only under the numerical scheme.
A numerical procedure is introduced for approximating such slow manifolds
in more general but still uniformly stiff mechanical systems by Ariel
et al. \cite{ariel12}. In a more mathematical treatment, Stumpp \cite{stumpp08}
considered general mechanical systems with stiff damping forces and
showed the existence of an attracting slow manifold governing the
asymptotic behavior of the system. Again, all degrees of freedom are
assumed equally stiff and hence no modes can be eliminated via model
reduction. 

In contrast to these specific case studies and purely stiff reduction
procedures, we consider here general mechanical systems and establish
conditions under which stiffer degrees of freedom can be identified
and eliminated by reduction to an attracting slow manifold defined
over the remaining softer degrees of freedom. We do not assume any
specific force or inertia term to be large or small. Rather, we seek
the broadest set of conditions under which an exact slow-fast decomposition
emerges and yields a reduced-order mechanical system. The slow-fast
decomposition (SDF) procedure arising form our analysis satisfies
the key requirements (R1)-(R2) discussed above. 

We also establish the maximal domain of SFD, and give a specific upper
bound on the rate at which general solutions synchronize with those
of the reduced-order model. Our includes several classes of mechanical
systems and justifies earlier heuristic reduction schemes under certain
conditions. In particular, under the SFD conditions, the techniques
of static condensation and modal derivatives, respectively, can rigorously
be justified as first- and second-order local approximations to a
slow manifold near an equilibrium. At the same time, we give examples
of these reduction procedures fail when the SFD conditions are not
met.

We illustrate these results on simple mechanical systems, but our
formulas are explicit enough to be applied to higher-degree-of-freedom
problems. Importantly, determining the eigenvalues and modes shapes
is not a prerequisite for the application of SFD. Indeed, the stiffer
modes may be fully (both linearly and nonlinearly) coupled to the
rest of the modes

\section{Set-up }

\subsection{General form of the mechanical system}

Consider a $n$-degree of freedom, non-dimensionalized mechanical
system of the form 
\begin{equation}
M(q,t)\ddot{q}-F(q,\dot{q},t)=0,\label{eq:systtem00}
\end{equation}
where $M\in\mathbb{R}^{n\times n}$ is a nonsingular mass matrix that
may depend on the generalized coordinates $q$ and the time $t$ in
a smooth fashion (class $C^{r}$ for some $r\geq2$). The internal
and external forces acting on the system are contained in the term
$F\in\mathbb{R}^{n}$, which generally depends on $q$, $t$ and the
generalized velocities $\dot{q}\in\mathbb{R}^{n}$. 

\subsection{Classic model reduction by projection to a subspace}

As noted in the Introduction, model reduction for system \eqref{eq:systtem00}
is generally motivated by an assumed coordinate change
\begin{equation}
q=Ux,\label{eq:general linear model reduction}
\end{equation}
with a matrix $U\in\mathbb{R}^{n\times k}$ and a reduced coordinate
vector $x\in\mathbb{R}^{k}$ with $k<n$ (see, e.g., Geradin and Rixen
\cite{geradin15}). Substitution into \eqref{eq:systtem00}, followed
by a multiplication by $U^{T}$, then suggests the reduced equations
of motion
\begin{equation}
U^{T}M(Ux,t)U\ddot{x}-U^{T}F(Ux,U\dot{x},t)=0,\label{eq:projected eq}
\end{equation}
the projection of \eqref{eq:systtem00} from the full state space
$\mathbb{R}^{n}$ onto a $k$-dimensional subspace $E$, parametrized
by the variable $x$ (cf. Fig. \ref{fig:projection}). The main focus
of model reduction studies is then the most expedient choice of the
matrix $U$.

It is often forgotten, however, that for eq. \eqref{eq:projected eq}
to hold, one must have $q(t)=Ux(t)$ for all times, i.e., $E$ must
be an invariant plane for \eqref{eq:systtem00}. This assumption is
practically certain to be violated unless special symmetries are present.
Even for unforced and stable structural system (i.e., when \eqref{eq:systtem00}
is autonomous and $q=0$ is asymptotically stable), the invariance
of modal subspaces is violated when nonlinear terms are present. The
mismatch between modal subspaces and (nonlinear) invariant manifolds
emanating from the origin will only be small very close to the origin.
In addition, various choices of $U$ may render projected equations
that do not capture typical dynamics even close to $q=0$ (cf. Haller
and Ponsioen \cite{haller16} and Section \ref{sec:Approximate-SFD-near-equilbria}
below for examples). 

\subsection{Slow (flexible) and fast (stiff) variables }

If the system \eqref{eq:systtem00} is non-autonomous, we assume that
its explicit time-dependence in $M$ and $F$ is precisely one of
the following three types: (1) periodic (2) quasiperiodic with finitely
many rationally independent frequencies (3) aperiodic over a finite
time interval $[a,b]$. In the periodic and quasiperiodic cases, we
let $t\in\mathcal{T}=\mathbb{R}$, whereas in the aperiodic case,
we let $t\in\mathcal{T}=[a,b]$. 

Next, we split the generalized coordinate vector $q$ as 
\[
q=\left(\begin{array}{c}
x\\
y
\end{array}\right),\quad x\in\mathbb{R}^{s},\qquad y\in\mathbb{R}^{f},\quad s+f=n,
\]
into yet unspecified \emph{slow coordinates} $x$ and \emph{fast coordinates}
$y$. This slow-fast partition refers to the expected relative speed
of variation of the $x$ and $y$ variables. In mechanical terms,
we expect $x$ to label relatively flexible degrees of freedom as
opposed to the relatively stiff degrees of freedoms labeled by the
$y$ coordinates. 

We seek conditions under which a mathematically rigorous model reduction
process exists to express $y(t)$ uniquely as function of $x(t)$,
at least asymptotically in time, along general trajectories $q(t)$
of \eqref{eq:systtem00}. The $(x,y)$ partition of $q$ may be suggested
by a modal analysis of the linear system or simply by the physics
of a mechanical problem. Importantly, $q$ is not assumed to be a
set of linear modal coordinates, and hence our procedure does not
rely on an a priori identification of a linearized spectrum near an
equilibrium point.

To allow for a potentially stiff dependence of the system on $y$,
we introduce a small, non-dimensional parameter $\epsilon>0$ and
consider $M$ and $F$ as smooth functions of $y/\epsilon$ and $y$
for $\epsilon>0.$ At this point, this represents no loss of generality,
given that \emph{any} smooth function of $y$ and $\epsilon$ can
also be viewed as a smooth function of $y/\epsilon$ and $\epsilon$
for $\epsilon>0$ because $y=\epsilon\cdot(y/\epsilon).$ 

Using this notation, we split the mass matrices and forcing terms
in \eqref{eq:systtem00} by letting
\[
M(q,t)=\left(\begin{array}{cc}
M_{11}\left(x,\frac{y}{\epsilon},t;\epsilon\right) & M_{12}\left(x,\frac{y}{\epsilon},t;\epsilon\right)\\
M_{21}\left(x,\frac{y}{\epsilon},t;\epsilon\right) & M_{22}\left(x,\frac{y}{\epsilon},t;\epsilon\right)
\end{array}\right),\qquad F(q,\dot{q},t)=\left(\begin{array}{c}
F_{1}\left(x,\dot{x},\frac{y}{\epsilon},\dot{y},t;\epsilon\right)\\
F_{2}\left(x,\dot{x},\frac{y}{\epsilon},\dot{y},t;\epsilon\right)
\end{array}\right),\qquad\epsilon>0,
\]
where $M_{11}\in\mathbb{R}^{s\times s}$, $M_{12},M_{21}^{T}\in\mathbb{R}^{s\times f}$,
$M_{22}\in\mathbb{R}^{f\times f}$, $F_{1}\in\mathbb{R}^{s}$ and
$F_{2}\in\mathbb{R}^{f}$. Again, as notated above, this notation
is general enough to allow for cases in which $M$ or $F$ depends
purely on $y$, or depends both on $y$ and $y/\epsilon.$ The corresponding
equations of motion are
\begin{eqnarray}
M_{11}\ddot{x}+M_{12}\ddot{y}-F_{1} & = & 0,\nonumber \\
M_{22}\ddot{y}+M_{21}\ddot{x}-F_{2} & = & 0.\label{eq:eq mo intermediate}
\end{eqnarray}
 Taking appropriate linear combination of these equations, and introducing
the matrix $M_{i}$ and forces $Q_{i}$ via 
\begin{eqnarray}
M_{i}\left(x,\frac{y}{\epsilon},t;\epsilon\right) & = & M_{ii}-M_{ij}M_{jj}^{-1}M_{ji},\qquad i,j=1,2,\quad i\neq j,\nonumber \\
Q_{i}\left(x,\dot{x},\frac{y}{\epsilon},\dot{y},t;\epsilon\right) & = & F_{i}-M_{ij}M_{jj}^{-1}F_{j},\qquad i,j=1,2,\quad i\neq j,\label{eq:M_i and Q_i def}
\end{eqnarray}
we deduce from \eqref{eq:eq mo intermediate} the inertially decoupled
equations of motion: 
\begin{eqnarray}
M_{1}\left(x,\frac{y}{\epsilon},t;\epsilon\right)\ddot{x}-Q_{1}\left(x,\dot{x},\frac{y}{\epsilon},\dot{y},t;\epsilon\right) & = & 0,\nonumber \\
M_{2}\left(x,\frac{y}{\epsilon},t;\epsilon\right)\ddot{y}-Q_{2}\left(x,\dot{x},\frac{y}{\epsilon},\dot{y},t;\epsilon\right) & = & 0.\label{eq:inertially decoupled0}
\end{eqnarray}

Importantly, the equations \eqref{eq:inertially decoupled0} are fully
equivalent to \eqref{eq:systtem00} for \emph{any} choice of the partition
$q=(x,y)$ and for \emph{any} choice of a scalar parameter $\epsilon>0$.
In particular, $M_{1}\left(x,\frac{y}{\epsilon},t;\epsilon\right)\in\mathbb{R}^{s\times s}$
and $M_{2}\left(x,\frac{y}{\epsilon},t;\epsilon\right)\in\mathbb{R}^{f\times f}$
are nonsingular matrices for all $(x,y,t)$ and for all $\epsilon>0$.
As we shall see in later examples, the partition $q=(x,y)$ will need
to be selected in given problems in a way that further assumptions
detailed below are satisfied.

\subsection{Assumptions of the SFD and illustrating examples}

We now list assumptions that will be sufficient to guarantee the existence
of an exact reduced-order model satisfying the requirements (R1)-(R2).
First, using the new variable $\eta=y/\epsilon$, we define the mass-normalized
forcing terms
\begin{eqnarray*}
P_{1}\left(x,\dot{x},\eta,\dot{y},t;\epsilon\right) & = & M_{1}^{-1}\left(x,\eta,t;\epsilon\right)Q_{1}\left(x,\dot{x},\eta,\dot{y},t;\epsilon\right),\\
P_{2}\left(x,\dot{x},\eta,\dot{y},t;\epsilon\right) & = & \epsilon M_{2}^{-1}\left(x,\eta,t;\epsilon\right)Q_{2}\left(x,\dot{x},\eta,\dot{y},t;\epsilon\right),
\end{eqnarray*}
which are, by our assumptions, class $C^{r}$ in their arguments for
$\epsilon>0$. The following assumptions concern properties of $P_{i}$
in their $\epsilon=0$ limit.
\begin{description}
\item [{(A1)}] \textbf{Nonsingular extension to $\epsilon=0$:} The functions
$P_{1}$ and $P_{2}$ are at least of class $C^{2}$ in their arguments
at $\epsilon=0$. 
\end{description}

In other words, assumption (A1) requires continuous differentiability
of the transformed forcing terms $P_{i}$ also in the limit of $\epsilon=0$
, when the dummy variable $\eta=y/\epsilon$ is held fixed, independent
of $\epsilon$. 
\begin{description}
\item [{(A2)}] \textbf{Existence of a fast zero-acceleration set (critical
manifold):} The algebraic equation $Q_{2}\left(x,\dot{x},\eta,0,t;0\right)\equiv0$
can be solved for $\eta$ on an open, bounded domain $\mathcal{D}_{0}\subset\mathbb{R}^{s}\times\mathbb{R}^{s}\times\mathcal{T}$.
Specifically, there exists a $C^{1}$ function $G_{0}:\mathcal{D}_{0}\to\mathbb{R}^{s}$
such that
\begin{equation}
Q_{2}\left(x,\dot{x},G_{0}(x,\dot{x},t),0,t;0\right)\equiv0\label{eq:P_2=00003D0}
\end{equation}
 holds for all $(x,\dot{x},t)\in\mathcal{D}_{0}$. We refer to the
set $\mathcal{M}_{0}(t)$ defined by $\eta=G_{0}(x,\dot{x},t)$ as
a \emph{critical manifold.}
\end{description}

Assumption (A2) ensures the existence of a smooth set $\mathcal{M}_{0}(t)$
of instantaneous zero-acceleration states\emph{ (}critical manifold)
for the fast coordinates. Here the velocity variable $\dot{x}\in\mathbb{R}^{s}$
is viewed as arbitrary, and hence unrelated to the actual time derivative
of $x(t)$ along a trajectory $q(t)$. As a consequence, these instantaneous
zero-acceleration states are not equilibria and do not form an invariant
set for system \eqref{eq:inertially decoupled0}. Under assumption
(A3) below, however, $\mathcal{M}_{0}(t)$ will turn out to approximate
a slow invariant manifold that carries a reduced-order model satisfying
the requirements (R1)-(R2):
\begin{description}
\item [{(A3)}] \textbf{Formal asymptotic stability of the critical manifold}:
With the matrices 
\begin{equation}
A(x,\dot{x},t)=-\partial_{\dot{y}}P_{2}\left(x,\dot{x},G_{0}(x,\dot{x},t),0,t;0\right),\qquad B(x,\dot{x},t)=-\partial_{\eta}P_{2}\left(x,\dot{x},G_{0}(x,\dot{x},t),0,t;0\right),\label{eq:ABdef}
\end{equation}
the equilibrium solution $\eta\equiv0\in\mathbb{R}^{f}$ of the unforced,
constant-coefficient linear system
\begin{equation}
\eta^{\prime\prime}+A(x,\dot{x},t)\eta^{\prime}+B(x,\dot{x},t)\eta=0\label{eq:associated linear system}
\end{equation}
is asymptotically stable for all fixed parameter values $(x,\dot{x},t)\in\mathcal{D}_{0}$.
Here prime denotes differentiation with respect to an auxiliary time
$\tau$ that is independent of $t$. 
\end{description}
Note that assumption (A3) requires the linear unforced oscillatory
system \eqref{eq:associated linear system}, posed formally for the
dummy fast variable $\eta$, to be asymptotically stable. In this
context, $(x,\dot{x},t)$ play the role of constant parameters ranging
over $\mathcal{D}_{0}$. Assumption (A3) is satisfied, for instance,
when $A$ is symmetric, positive semi-definite and $B$ is symmetric,
positive definite over $\mathcal{D}_{0}$. In that case, $A$ represents
a damping matrix and $B$ represents a stiffness matrix for all parameter
values $(x,\dot{x},t)\in\mathcal{D}_{0}$. At this point, (A3) is
only a formal requirement with no immediately clear mathematical meaning.
This is because the critical manifold $\mathcal{M}_{0}(t)$ is not
invariant under equation \eqref{eq:inertially decoupled0} and hence
the arguments of $A$ and $B$ are, in fact, time-varying, and hence
do not determine the stability of \eqref{eq:associated linear system}.
\begin{example}
\label{ex: general oscillatory system}{[}\emph{Weakly nonlinear system
with parametric forcing}{]} Consider a typical multi-degree-of-freedom
mechanical system of the form
\begin{eqnarray}
M_{1}\ddot{x}+C_{1}\dot{x}+K_{1}x+S_{1}\left(x,y\right) & = & f_{1}(t),\nonumber \\
M_{2}\ddot{y}+C_{2}\dot{y}+K_{2}y+S_{2}\left(x,y\right) & = & f_{2}(t),\label{eq:oscillatory system}
\end{eqnarray}
with $x\in\mathbb{R}^{s}$ and $y\in\mathbb{R}^{f}$. Here the $M_{i}$
are symmetric and positive definite constant mass matrices; $C_{i}$
are constant symmetric damping matrices; $K_{i}$ are constant symmetric
stiffness matrices; and the functions 
\begin{equation}
S_{i}(x,y)=\mathcal{O}\left(\left|x\right|^{2},\text{\ensuremath{\left|x\right|}}\left|y\right|,\left|y\right|^{2}\right)\label{eq:quadraticS_i}
\end{equation}
model nonlinear coupling terms. By definition \eqref{eq:M_i and Q_i def},
for an arbitrary scalar parameter $\epsilon>0$ independent of $M_{i}$,
$C_{i}$, $D_{i}$ and $S_{i}$, we specifically have
\begin{eqnarray*}
Q_{1}\left(x,\dot{x},\frac{y}{\epsilon},\dot{y},t;\epsilon\right) & = & -\left[C_{1}\dot{x}+K_{1}x+S_{1}\left(x,\epsilon\frac{y}{\epsilon}\right)-f_{1}(t)\right],\\
Q_{2}\left(x,\dot{x},\frac{y}{\epsilon},\dot{y},t;\epsilon\right) & = & -\left[C_{2}\dot{y}+\epsilon K_{2}\left(\frac{y}{\epsilon}\right)+S_{2}\left(x,\epsilon\frac{y}{\epsilon}\right)-f_{2}(t)\right].
\end{eqnarray*}
 Therefore, the functions
\begin{eqnarray*}
P_{1}(x,\dot{x},\eta,\dot{y},t;\epsilon) & = & -M_{1}^{-1}\left[C_{1}\dot{x}+K_{1}x+S_{1}\left(x,\epsilon\eta\right)-f_{1}(t)\right],\\
P_{2}(x,\dot{x},\eta,\dot{y},t;\epsilon) & = & -\epsilon M_{2}^{-1}\left[C_{2}\dot{y}+\epsilon K_{2}\eta+S_{2}\left(x,\epsilon\eta\right)-f_{2}(t)\right],
\end{eqnarray*}
are differentiable in $\epsilon$ at the the $\epsilon=0$ limit,
satisfying assumption (A1). However, we have 
\[
P_{2}(x,\dot{x},\eta,\dot{y},t;0)\equiv0.
\]
Therefore, while any function $G_{0}(x,\dot{x},t)$ satisfies (A2),
both matrices $A$ and $B$ defined in \eqref{eq:ABdef} vanish, and
hence assumption (A3) never holds for system \eqref{eq:oscillatory system}.
For this assumption to hold, some of the system parameters must be
related to the small parameter $\epsilon$, as we shall see in the
next two examples.
\end{example}

\begin{example}
\label{ex: Partially stiff oscillatory system-1}{[}\emph{Partially
stiff weakly nonlinear system with very small stiff-inertia and parametric
forcing}{]} Consider now the slightly modified multi-degree-of-freedom
mechanical system
\begin{eqnarray}
M_{1}\ddot{x}+C_{1}\dot{x}+K_{1}x+S_{1}\left(x,y\right) & = & f_{1}(t),\nonumber \\
\epsilon^{2}M_{2}\ddot{y}+C_{2}\dot{y}+\frac{1}{\epsilon}K_{2}y+S_{2}\left(x,y\right) & = & f_{2}(t),\label{eq:oscillatory system-1-2}
\end{eqnarray}
with a non-dimensional small parameter $\epsilon\ll1$. All variables,
matrices and functions are the same as in Example \ref{ex: general oscillatory system},
but the $y$-component of this system generates very small inertial
forces and also has large linear stiffness. This time, we have
\begin{eqnarray*}
P_{1}(x,\dot{x},\eta,\dot{y},t;\epsilon) & = & -M_{1}^{-1}\left[C_{1}\dot{x}+K_{1}x+S_{1}\left(x,\epsilon\eta\right)-f_{1}(t)\right],\\
P_{2}(x,\dot{x},\eta,\dot{y},t;\epsilon) & = & -\frac{1}{\epsilon}M_{2}^{-1}\left[C_{2}\dot{y}+K_{2}\eta+S_{2}\left(x,\epsilon\eta\right)-f_{2}(t)\right],
\end{eqnarray*}
therefore assumption (A1) is not satisfied, given that $P_{2}$ is
not differentiable at $\epsilon=0$ . 
\end{example}

\begin{example}
\label{ex: Partially stiff oscillatory system-1-1}{[}\emph{Paradigm
for targeted energy transfer: Weakly nonlinear system with small inertia
in its essentially nonlinear component}{]} Consider the multi-degree-of-freedom
mechanical system
\begin{eqnarray}
M_{1}\ddot{x}+C_{1}\dot{x}+K_{1}x+S_{1}\left(x,y\right) & = & 0,\nonumber \\
\epsilon M_{2}\ddot{y}+C_{2}\dot{y}+S_{2}\left(x,y\right) & = & 0,\label{eq:oscillatory system-1-2-1}
\end{eqnarray}
with a non-dimensional small parameter $\epsilon\ll1$. Again, all
variables and matrices are the same as in Example \eqref{ex: general oscillatory system},
but the $y$-component of \eqref{eq:oscillatory system-1-2-1} generates
small inertial forces and no linear stiffness forces. This system
is noted as a prototype example of targeted energy transfer (cf. Vakakis
et al. \cite{vakakis08}) from the $x$ degrees of freedom to the
$y$ degrees of freedom. This energy transfer mechanism suggests the
lack of a reduced-order model over the $x$-degrees of freedom, given
that the $y$-variables display no long-term enslavement to the $x$-variables.
Calculating the quantities in our assumption (A1), we find 
\begin{eqnarray*}
P_{1}(x,\dot{x},\eta,\dot{y},t;\epsilon) & = & -M_{1}^{-1}\left[C_{1}\dot{x}+K_{1}x+S_{1}\left(x,\epsilon\eta\right)\right],\\
P_{2}(x,\dot{x},\eta,\dot{y},t;\epsilon) & = & -M_{2}^{-1}\left[C_{2}\dot{y}+S_{2}\left(x,\epsilon\eta\right)\right],
\end{eqnarray*}
are differentiable at $\epsilon=0$, and hence assumption (A1) holds.
However, the equation 
\[
Q_{2}(x,\dot{x},\eta,0,t;0)=-S_{2}\left(x,0\right)=0
\]
 cannot be solved for the variable $\eta$ at any point. As a consequence,
even though a set of zero acceleration states is defined by the equation
$S_{2}\left(x,0\right)=0,$ this set is not attracting. Indeed, the
matrix $B(x,v,t)$ defined in assumption (A3) vanishes identically
and hence the linear system \eqref{eq:associated linear system} is
not asymptotically stable.
\end{example}

\begin{example}
\label{ex: Partially-stiff-weakly}{[}\emph{Partially stiff weakly
nonlinear system with parametric forcing}{]} Consider now the multi-degree-of-freedom
mechanical system
\begin{eqnarray}
M_{1}\ddot{x}+C_{1}\dot{x}+K_{1}x+S_{1}\left(x,\frac{y}{\epsilon}\right) & = & f_{1}(t),\nonumber \\
\epsilon M_{2}\ddot{y}+C_{2}\dot{y}+\frac{1}{\epsilon}K_{2}y+S_{2}\left(x,y\right) & = & f_{2}(t),\label{eq:oscillatory system-1-1}
\end{eqnarray}
with the same quantities as in Example \eqref{ex: general oscillatory system}.
The difference here is that the mass matrix of the $y$ degrees of
freedom has small norm for $\epsilon\ll1$ and the stiffness matrix
is large in norm in the same equation. In addition, the nonlinear
coupling term in the $x$-equation is assumed to have a stiff dependence
on $y$. In this case, we have 
\begin{eqnarray*}
P_{1}(x,\dot{x},\eta,\dot{y},t;\epsilon) & = & -M_{1}^{-1}\left[C_{1}\dot{x}+K_{1}x+S_{1}\left(x,\eta\right)-f_{1}(t)\right],\\
P_{2}(x,\dot{x},\eta,\dot{y},t;\epsilon) & = & -M_{2}^{-1}\left[C_{2}\dot{y}+K_{2}\eta+S_{2}\left(x,\epsilon\eta\right)-f_{2}(t)\right],
\end{eqnarray*}
which satisfy assumption (A1). Solving the equation $Q_{2}(x,\dot{x},\eta,0,t;0)=0$
for $\eta,$ we find that assumption (A2) is satisfied by the function
\begin{equation}
G_{0}(x,\dot{x},t)=K_{2}^{-1}\left[f_{2}(t)-S_{2}(x,0)\right],\qquad(x,\dot{x},t)\in\mathcal{D}_{0}=\mathbb{R}^{s}\times\mathbb{R}^{s}\times\mathbb{R},\label{eq:G_0 in partially stiff system}
\end{equation}
provided that the stiffness matrix $K_{2}$ is invertible. In that
case, we obtain 
\[
A(x,\dot{x},t)=M_{2}^{-1}C_{2},\qquad B(x,\dot{x},t)=M_{2}^{-1}\left[K_{2}+\epsilon\partial_{y}S_{2}\left(x,\epsilon G_{0}(x,\dot{x},t)\right)\right]\vert_{\epsilon=0}=M_{2}^{-1}K_{2},
\]
and hence the homogeneous linear oscillatory system in assumption
(A3) becomes 
\[
\eta^{\prime\prime}+M_{2}^{-1}C_{2}\eta^{\prime}+M_{2}^{-1}K_{2}\eta=0
\]
or, equivalently, 
\begin{equation}
M_{2}\eta^{\prime\prime}+C_{2}\eta^{\prime}+K_{2}\eta=0.\label{eq:equivalent linear}
\end{equation}
The zero equilibrium of this system is asymptotically stable by our
assumptions on $M_{2}$, $C_{2}$ and $K_{2}$. Therefore, assumption
(A3) is also satisfied for system \eqref{eq:oscillatory system-1-1}. 
\end{example}

\section{Main result: Global existence of an exact reduced-order model }

\begin{figure}[H]
\centering{}\includegraphics[width=0.4\textwidth]{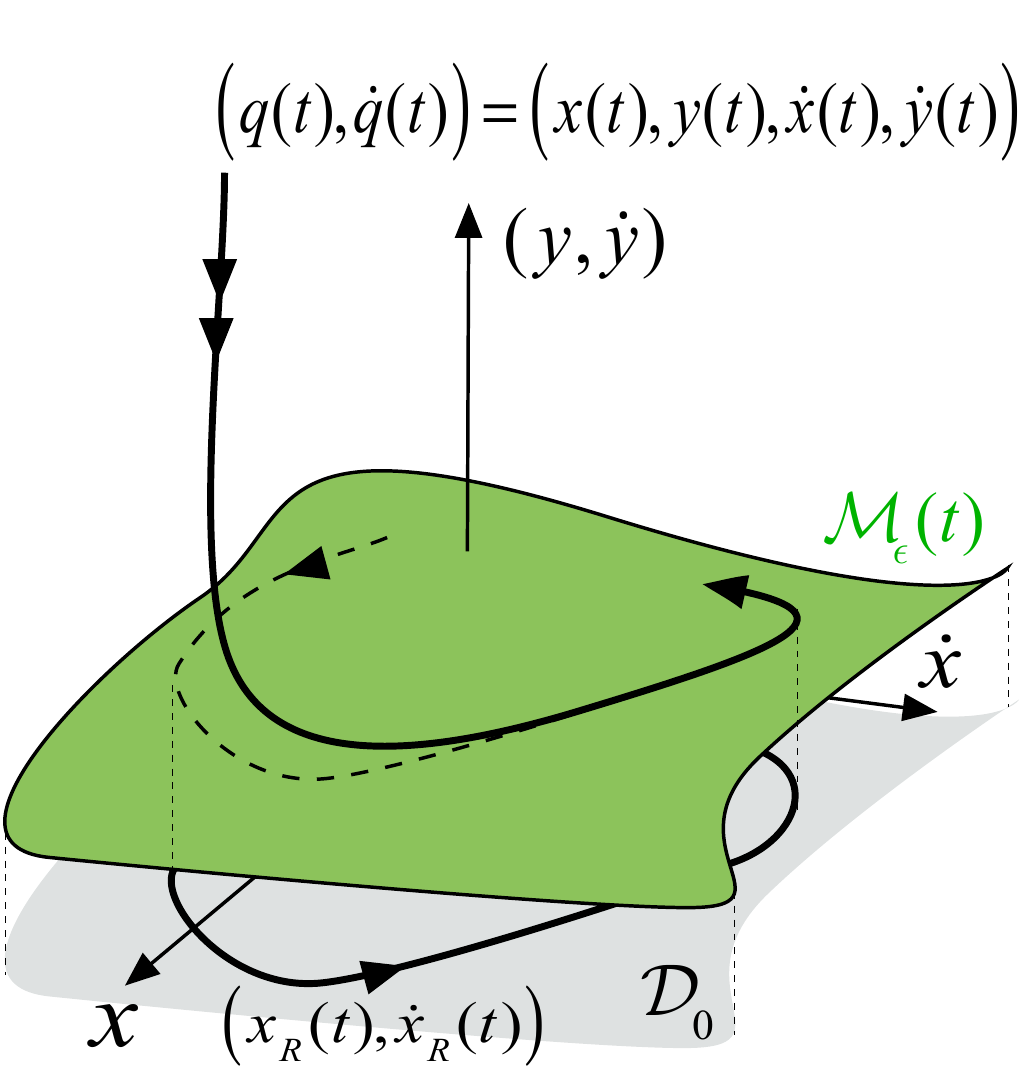}\caption{Reduced-order model trajectory $(x_{R}(t),\dot{x}_{R}(t))$ as a projection
from the slow manifold $\mathcal{M_{\epsilon}}(t)$ to the space of
the $(x,\dot{x})$ variables. Other nearby trajectories converge to
the slow manifold exponentially fast, and hence their projection on
the $(x,\dot{x})$ space synchronizes exponentially with trajectories
of the reduced-order model. \label{fig:SFD main} }
\end{figure}

To state our main result formally, we first define the following functions
for all $(x,\dot{x},t)\in\mathcal{D}_{0}:$ 
\begin{eqnarray}
H_{0}(x,\dot{x},t) & = & \partial_{x}G_{0}(x,\dot{x},t)\dot{x}+\partial_{\dot{x}}G_{0}(x,\dot{x},t)P_{1}\left(x,\dot{x},G_{0}(x,\dot{x},t),0,t;0\right)+\partial_{t}G_{0}(x,\dot{x},t),\nonumber \\
G_{1}(x,\dot{x},t) & = & -\left[D_{\eta}P_{2}\left(x,\dot{x},G_{0}(x,\dot{x},t),0,t;0\right)\right]^{-1}D_{\dot{y}}P_{2}\left(x,\dot{x},G_{0}(x,\dot{x},t),0,t;0\right)H_{0}(x,\dot{x},t)\nonumber \\
 &  & -\left[D_{\eta}P_{2}\left(x,\dot{x},G_{0}(x,\dot{x},t),0,t;0\right)\right]^{-1}D_{\epsilon}P_{2}\left(x,\dot{x},G_{0}(x,\dot{x},t),0,t;0\right),\label{eq:H0G1H1}\\
H_{1}(x,\dot{x},t) & = & \partial_{x}G_{1}(x,\dot{x},t)v+\partial_{\dot{x}}G_{1}(x,\dot{x},t)P_{1}\left(x,\dot{x},G_{0}(x,\dot{x},t),0,t;0\right)+\partial_{t}G_{1}(x,\dot{x},t).\nonumber 
\end{eqnarray}
With these quantitates, we have the following result:
\begin{thm}
\label{thm:main theorem}Under assumptions (A1)-(A3) and for $\epsilon>0$
small enough:
\end{thm}
\begin{description}
\item [{\emph{(i)}}] \emph{The mechanical system \eqref{eq:systtem00}
admits an exact reduced-order model satisfying the requirements (R1)-(R2). }
\item [{\emph{(ii)}}] \emph{The reduced-order model is given by 
\begin{eqnarray}
\ddot{x}-P_{1}\left(x,\dot{x},G_{0}(x,\dot{x},t),0,t;0\right) & = & \epsilon\left[D_{\eta}P_{1}\left(x,\dot{x},G_{0}(x,\dot{x},t),0,t;0\right)G_{1}(x,\dot{x},t)\right.\nonumber \\
 &  & \,\,\,\,\,\,+D_{\dot{y}}P_{1}\left(x,\dot{x},G_{0}(x,\dot{x},t),0,t;0\right)H_{0}(x,\dot{x},t)\label{eq:reduced model}\\
 &  & \,\,\,\,\,\,+\left.D_{\epsilon}P_{1}\left(x,\dot{x},G_{0}(x,\dot{x},t),0,t;0\right)\right]\nonumber \\
 &  & +\mathcal{O}(\epsilon^{2})\nonumber 
\end{eqnarray}
for all $(x,\dot{x},t)\in\mathcal{D}_{0}$.}
\item [{\emph{(iii)}}] \emph{If $M_{1}(x,\eta,\epsilon)$ is smooth in
$\epsilon$ at $\epsilon=0$, then the multiplication of \eqref{eq:reduced model}
by $M_{1}$ gives a form of the reduced-order model that does not
require the inversion of $M_{1}$: 
\begin{equation}
M_{1}\left(x,G_{0}(x,\dot{x},t),t;0\right)\ddot{x}-Q_{1}\left(x,v,G_{0}(x,\dot{x},t),0,t;0\right)=\mathcal{O}(\epsilon).\label{eq:reduced model 1}
\end{equation}
}
\item [{\emph{(iv)}}] \emph{The reduced-order models \eqref{eq:reduced model}-\eqref{eq:reduced model 1}
describe the reduced flow on a $2s$-dimensional invariant manifold
$\mathcal{M_{\epsilon}}(t)$ along which positions and velocities
in the stiff degrees of freedom are enslaved to those in the slow
degrees of freedom via 
\begin{eqnarray}
y & = & \epsilon G_{0}(x,\dot{x},t)+\epsilon^{2}G_{1}(x,\dot{x},t)+\mathcal{O}(\epsilon^{3}),\nonumber \\
\dot{y} & = & \epsilon H_{0}(x,\dot{x},t)+\epsilon^{2}H_{1}(x,\dot{x},t)+\mathcal{O}(\epsilon^{3}).\label{eq:slow man coords}
\end{eqnarray}
}
\item [{\emph{(v)}}] \emph{The $x(t)$ components of the trajectories of
system \eqref{eq:systtem00} synchronize with appropriate model trajectories
$x_{R}(t)$ of \eqref{eq:reduced model} or \eqref{eq:reduced model 1}
at an exponential rate. Specifically, let $q(t)=\left(x(t),y(t)\right)$
be a full trajectory of system \eqref{eq:systtem00} such that at
a time $t_{0}$, the initial position $q(t_{0})$ is close enough
to the slow manifold carrying the reduced order model. Then there
exists a trajectory $x_{R}(t)$ of the reduced-order model \eqref{eq:reduced model}
or \eqref{eq:reduced model 1} such that
\begin{equation}
\left|\left(\begin{array}{c}
x(t)-x_{R}(t)\\
\dot{x}(t)-\dot{x}_{R}(t)
\end{array}\right)\right|\leq C\left|\left(\begin{array}{c}
x(t_{0})-x_{R}(t_{0})\\
\dot{x}(t_{0})-\dot{x}_{R}(t_{0})\\
\frac{1}{\epsilon}y(t_{0})-G_{0}\left(x_{R}(t_{0}),\dot{x}_{R}(t_{0}),t\right)+\mathcal{O}(\epsilon)\\
\dot{y}(t_{0})-\epsilon H_{0}\left(x_{R}(t_{0}),\dot{x}_{R}(t_{0}),t\right)+\mathcal{O}(\epsilon^{2})
\end{array}\right)\right|e^{-\frac{\Lambda}{\epsilon}(t-t_{0})},\quad t>t_{0},\label{eq:synchronization inequaltiy}
\end{equation}
 where $\Lambda>0$ can be selected as any constant satisfying 
\[
\max_{j\in[1,f],\,(x,\dot{x},t)\in\mathcal{D}_{0}}\mathrm{Re}\,\lambda_{j}(x,\dot{x},t)<-\Lambda<0,
\]
with $\lambda_{j}(x,\dot{x},t),$ $j=1,\ldots,f$, denoting the eigenvalues
of the associated linear system \eqref{eq:associated linear problem}.
The constant $C>0$ generally depends on the choice of $\Lambda$
but is independent of the choice of the initial conditions $q(t_{0})$
and $\dot{q}(t_{0})$.}
\end{description}
\begin{proof}
See Appendix \eqref{app: proof of main}.
\end{proof}
In Fig. \ref{fig:SFD main}, we illustrate the geometric relation
between the reduced model flow on the slow manifold to general trajectories
of the full system, as described by Theorem \ref{thm:main theorem}.

\begin{rem}
If the left-hand side of the reduced-order model \eqref{eq:reduced model}
has structurally stable features (cf. Guckenheimer and Holmes \cite{guckenheimer83}),
then, for $\epsilon>0$ small enough, these features persist smoothly
under the addition of the $\mathcal{O}(\epsilon)$ terms of the right-hand
side, and hence an explicit computation of these terms is not necessary.
For instance, if system \eqref{eq:reduced model} has a single attracting
fixed point or periodic orbit over the compact domain $\mathcal{D}_{0}$,
then wither of these features is robust without the explicit inclusion
of the $\mathcal{O}(\epsilon)$ and higher-order terms on its right-hand
side. If, however, system \eqref{eq:reduced model} is conservative,
then the inclusion of $\mathcal{O}(\epsilon)$ terms is necessary
to obtain a robust, dissipative reduced-order model. If the $\mathcal{O}(\epsilon)$
terms are also conservative, then explicit evaluation of the $\mathcal{O}(\epsilon^{2})$
is required following the expansion scheme used in the proof of Theorem
\ref{thm:main theorem}.
\end{rem}

\begin{rem}
The synchronization expressed by \eqref{eq:synchronization inequaltiy}
means that both positions and velocities predicted by the reduced-order
model \eqref{eq:reduced model} are relevant for the observed system
dynamics as long the time $t>t_{0}$ is selected from the domain $\mathcal{T}$.
This time-domain is unbounded (i.e., $\mathcal{T}=\mathbb{R}$) for
mechanical systems with explicit periodic and quasiperiodic time dependence.
For the case of temporally aperiodic time dependence, the times allowed
in \eqref{eq:synchronization inequaltiy} are restricted to the finite
interval $\mathcal{T}=[a,b]$. 
\end{rem}
\begin{example}
{[}\emph{Partially stiff weakly nonlinear system with parametric forcing}{]}
We recall that the partially stiff system \eqref{eq:oscillatory system-1-1}
in Example \ref{ex: Partially-stiff-weakly} satisfies assumptions
(A1)-(A3) and hence admits an exact, global reduced-order model. The
form of the function $G_{0}$ from \eqref{eq:G_0 in partially stiff system}
is 
\begin{equation}
G_{0}(x,\dot{x},t)=K_{2}^{-1}\left[f_{2}(t)-S_{2}(x,0)\right].\label{eq:G_0 in partially stiff system-1}
\end{equation}
The mass matrix $M_{1}$ is independent of $\epsilon$, and hence
the equivalent form \eqref{eq:reduced model 1} of the reduced-order
model applies and gives 
\[
M_{1}\ddot{x}+C_{1}\dot{x}+K_{1}x+S_{1}\left(x,K_{2}^{-1}\left[f_{2}(t)-S_{2}(x,0)\right]\right)=f_{1}(t)+O(\epsilon).
\]
The leading-order terms in the expressions \eqref{eq:slow man coords}
for the slow manifold are
\begin{eqnarray*}
y & = & \epsilon G_{0}(x,\dot{x},t)+\mathcal{O}(\epsilon^{2})=\epsilon K_{2}^{-1}\left[f_{2}(t)-S_{2}(x,0)\right]+\mathcal{O}(\epsilon^{2}),\\
\dot{y} & = & \epsilon H_{0}(x,\dot{x},t)+\mathcal{O}(\epsilon^{2})=\epsilon K_{2}^{-1}\left[\dot{f}_{2}(t)-\partial_{x}S_{2}(x,0)\dot{x}\right]+\mathcal{O}(\epsilon^{2}).
\end{eqnarray*}
 If $f(t)$ is periodic or quasiperiodic in time, then we have the
synchronization estimate \eqref{eq:synchronization inequaltiy} for
all times $t>t_{0}$. Specifically, any $\Lambda>0$ can be selected
such that $-\Lambda<0$ is a strict upper bound on the real part of
the spectrum of the oscillatory system \eqref{eq:equivalent linear}.
We note that if we had assumed a non-stiff coupling of the form $S_{1}\left(x,y\right)$
in Example \ref{ex: Partially-stiff-weakly}, then assumptions (A1)-(A4)
would still have been satisfied, but the reduced model would simplify
to 
\[
M_{1}\ddot{x}+C_{1}\dot{x}+K_{1}x+S_{1}\left(x,0\right)=f_{1}(t)+O(\epsilon),
\]
uncoupling completely from the stiff modes at leading order. The convergence
estimate \eqref{eq:synchronization inequaltiy} would remain valid
in this case, too.
\end{example}

\section{The boundary of the domain of model reduction}

In the examples we have discussed so far, the domain $\mathcal{D}_{0}$
could be selected arbitrarily large. Thus, a reduced-order model exists
over arbitrarily large $(x,\dot{x},t)$ values in these problems,
as long as $\epsilon$ is kept small enough. In general, however,
$\mathcal{D}_{0}$ will have a nonempty boundary $\partial\mathcal{D}_{0}$
over which the reduced-order model \eqref{eq:reduced model}-\eqref{eq:reduced model 1}
cannot be further extended. 

Such non-extendibility of the reduced-order model domain arises from
a break-down in the solvability of the algebraic equation \eqref{eq:P_2=00003D0}
for the critical manifold. By the implicit function theorem, this
occurs along points satisfying 
\begin{equation}
\det\left[\partial_{\eta}P_{2}(x,\dot{x},G_{0}(x,\dot{x},t),0,t;0)\right]=0,\quad(x,\dot{x},t)\in\partial\mathcal{D}_{0}.\label{eq:detcond}
\end{equation}
In the generic case, this determinant becomes zero at points where
$\partial_{\eta}P_{2}$ has a single zero eigenvalue, i.e., 
\begin{equation}
\mathrm{rank}\left[\partial_{\eta}P_{2}(x,\dot{x},G_{0}(x,\dot{x},t),0,t;0)\right]=f-1,\quad(x,\dot{x},t)\in\partial\mathcal{D}_{0}.\label{eq:detcond2}
\end{equation}

Under further nondegeneracy conditions (see., e..g., Arnold \cite{arnold92}),
a fold develops in the critical manifold along $\partial\mathcal{D}_{0}$,
i.e., $\mathcal{M}_{0}\equiv\mathcal{M}_{0}^{+}$ ceases to be a locally
unique graph over the $(x,\dot{x},t)$ variables. As we pass from
$\mathcal{M}_{0}$ to the newly bifurcating critical manifold branch
$\mathcal{\mathcal{M}}_{0}^{-}$, the matrices $A(x,\dot{x},t)$ and
$B(x,\dot{x},t)$ vary smoothly in their arguments, given that one
manifold branch is smoothly connected to the other one along a fold.
Under the nondegeneracy condition \eqref{eq:detcond2}, precisely
one eigenvalue of the matrix $B(x,\dot{x},t)$ will cross zero in
the passage from $\mathcal{M}_{0}^{+}$ to $\mathcal{M}_{0}^{-}$
along the critical manifold. In this case, the graph segment $\eta=G_{0}^{_{-}}(x,\dot{x},t)$
describing the bifurcating branch $\mathcal{M}_{0}^{-}$ \eqref{eq:associated linear system}
violates assumption (A2). As a consequence, the folded slow manifold
branch $\mathcal{M}_{\epsilon}^{-}$ perturbing from $\mathcal{M}_{0}^{-}$
is unstable and hence irrelevant for reduced-order modeling. 

In summary, unlike in the setting of the local construction of spectral
submanifolds near equilibria (cf. Haller and Ponsioen \cite{haller16}),
a folding invariant manifold arising in SFD is not a technical limitation
to overcome when one is in pursuit of a more global reduced-order
model. Rather, a fold in the slow manifold over the plane of slow
variables signals precisely the limit beyond which no reduced-order
model satisfying (R1)-(R2) exists in a given part of the phase space.

\begin{example}
\label{ex: Partially-stiff-weakly-1}{[}\emph{Partially stiff weakly
nonlinear system with parametric forcing}{]} Consider now the multi-degree-of-freedom
mechanical system
\begin{eqnarray}
M_{1}\ddot{x}+C_{1}\dot{x}+K_{1}x+S_{1}\left(x,\frac{y}{\epsilon}\right) & = & f_{1}(t),\nonumber \\
\epsilon M_{2}\ddot{y}+C_{2}\dot{y}+\frac{1}{\epsilon}K_{2}y+S_{2}\left(x,\frac{y}{\epsilon}\right) & = & f_{2}(t),\label{eq:oscillatory system-1-1-1}
\end{eqnarray}
with the same variables, matrices and functions used in Example \ref{ex: Partially-stiff-weakly},
except that here the coupling function $S_{2}$ also has a stiff dependence
on the $y$ variables. We then obtain
\begin{eqnarray*}
P_{1}(x,\dot{x},\eta,\dot{y},t;\epsilon) & = & -M_{1}^{-1}\left[C_{1}\dot{x}+K_{1}x+S_{1}\left(x,\eta\right)-f_{1}(t)\right],\\
P_{2}(x,\dot{x},\eta,\dot{y},t;\epsilon) & = & -M_{2}^{-1}\left[C_{2}\dot{y}+K_{2}\eta+S_{2}\left(x,\eta\right)-f_{2}(t)\right],
\end{eqnarray*}
thus assumption (A1) is satisfied again. The condition \eqref{eq:detcond}
in this case gives 
\[
\det\left[M_{2}^{-1}\left(K_{2}+\partial_{y}S_{2}\left(x,\eta\right)\right)\right]\neq0.
\]
 By the non-singularity of $M_{2}$, this latter condition is equivalent
to 
\[
\det\left[K_{2}+\partial_{y}S_{2}\left(x,\eta\right)\right]\neq0.
\]
 For instance, when $S_{2}$ has only quadratic terms, then this last
condition can always be written as
\begin{equation}
\det\left[K_{2}+\Pi x+\Phi\eta\right]\neq0,\label{eq:detnonzero-1}
\end{equation}
where $\Pi$ and $\Phi$ are 3-tensors of appropriate dimensions.
Suppose now, for simplicity, that $\Phi\equiv0$, the master variable
$x$ is a scalar $(m=1)$, and $\Pi\in\mathbb{R}^{f\times f}$ is
nonsingular. The requirement \eqref{eq:detnonzero-1} then becomes
\begin{equation}
\det\left[\Pi^{-1}K_{2}-(-x)I\right]\neq0,\label{eq:eig sing cond}
\end{equation}
which implies that $-x$ cannot be an eigenvalue of $\Pi^{-1}K_{2}$.
Consequently, condition \eqref{eq:eig sing cond} fails along the
domain boundary
\[
\partial\mathcal{D}_{0}=\left\{ (x,\dot{x},t)\,:\,\exists j:\,\,\,x=-\lambda_{j}(\Pi^{-1}K_{2}),\right\} ,
\]
with $\lambda_{j}(\Pi^{-1}K_{2})$ denoting the $j^{th}$ real eigenvalue
of the matrix $\Pi^{-1}K_{2}$. 

To illustrate the geometry of the critical manifold in a simple case,
we let $s=f=1$ and select the parameters, the coupling and the forcing
terms as 
\[
K_{2}=4,\quad S_{2}(x,\eta)=x^{2}+4\eta^{2},\quad f_{2}(t)=\sin t,
\]
so that the equation $P_{2}(x,\dot{x},\eta,0,t;0)=0$ takes the form
\begin{equation}
4\eta+x^{2}+4\eta^{2}-\sin t=0.\label{eq:p2=00003D0_ex}
\end{equation}
This equation is solved by $\eta=x=t=0$ and hence the set $\mathcal{D}_{0}$
is nonempty. The boundary $\partial\mathcal{D}_{0}$, defined by conditions
\eqref{eq:detcond}-\eqref{eq:detcond2}, satisfies
\begin{eqnarray*}
\det\left[\partial_{\eta}P_{2}(x,\dot{x},\eta,0,t;0)\right] & = & \partial_{\eta}P_{2}(x,\dot{x},\eta,0,t;0)=4+8\eta=0\quad\iff\quad\eta=-\frac{1}{2},\\
\partial_{\eta}\det\left[\partial_{\eta}P_{2}\left(x,\dot{x},-\frac{1}{2},0,t;0\right)\right] & = & 8\neq0,
\end{eqnarray*}
where the second condition here is the classic nondegeneracy condition
for fold bifurcations in the one-dimensional case (cf. Arnold \cite{arnold92}).
Substitution of $\eta=-\frac{1}{2}$ into the equation\eqref{eq:p2=00003D0_ex}
gives an explicit definition for $\partial\mathcal{D}_{0}$ in the
$(x,\dot{x},t)$ space as 
\begin{equation}
\partial\mathcal{D}_{0}=\left\{ (x,\dot{x},t)\,:\,x^{2}=1+\sin t\,\right\} .\label{eq:partial D0 ex}
\end{equation}
A direct solution of equation \eqref{eq:p2=00003D0_ex} through the
quadratic formula confirms that the zero set 
\[
\eta=G_{0}^{\pm}(x,\dot{x},t)=\frac{-1\pm\sqrt{1-(x^{2}-\sin t)}}{2},\quad(x,\dot{x},t)\in\mathcal{D}_{0}=\left\{ (x,\dot{x},t)\,:\,x^{2}<1+\sin t\,\right\} 
\]
indeed ceases to be a graph and develops a fold singularity along
$\partial\mathcal{D}_{0}$. The stability of the two branches of $G_{0}^{\pm}(x,\dot{x},t)$
can be determined by calculating \eqref{eq:associated linear system}
along both branches:
\begin{eqnarray*}
A^{\pm}(x,\dot{x},t) & = & M_{2}^{-1}C_{2},\\
B^{\pm}(x,\dot{x},t) & = & -M_{2}^{-1}\left(4+8G_{0}^{\pm}(x,\dot{x},t)\right)=\mp4\sqrt{1-(x^{2}-\sin t)}.
\end{eqnarray*}
Therefore, the critical manifold 
\[
\mathcal{M}_{0}^{+}=\left\{ (x,\dot{x},t)\in\mathcal{D}_{0}\,:\,\eta=G_{0}^{+}(x,\dot{x},t)\:\right\} 
\]
satisfies assumption (A1)-(A3) but develops a fold over $\mathcal{D}_{0}$
along the boundary curve $\partial\mathcal{D}_{0}$ defined in \eqref{eq:partial D0 ex}.
The additional branch 
\[
\mathcal{M}_{0}^{-}=\left\{ (x,\dot{x},t)\in\mathcal{D}_{0}\,:\,\eta=G_{0}^{-}(x,\dot{x},t)\;\right\} 
\]
emanating from the domain boundary $\partial\mathcal{D}_{0}$ is unstable,
as its associated constant-coefficient linear system (cf. assumption
(A2)), given by 
\[
M_{2}u^{\prime\prime}+C_{2}u^{\prime}-4\sqrt{1-(x^{2}-\sin t)}u=0,
\]
is unstable. We show the stable and unstable critical manifolds, as
well as the domain boundary $\partial\mathcal{D}_{0}$, in Fig. \ref{fig:fold}.
\begin{figure}[H]
\begin{centering}
\includegraphics[width=0.7\textwidth]{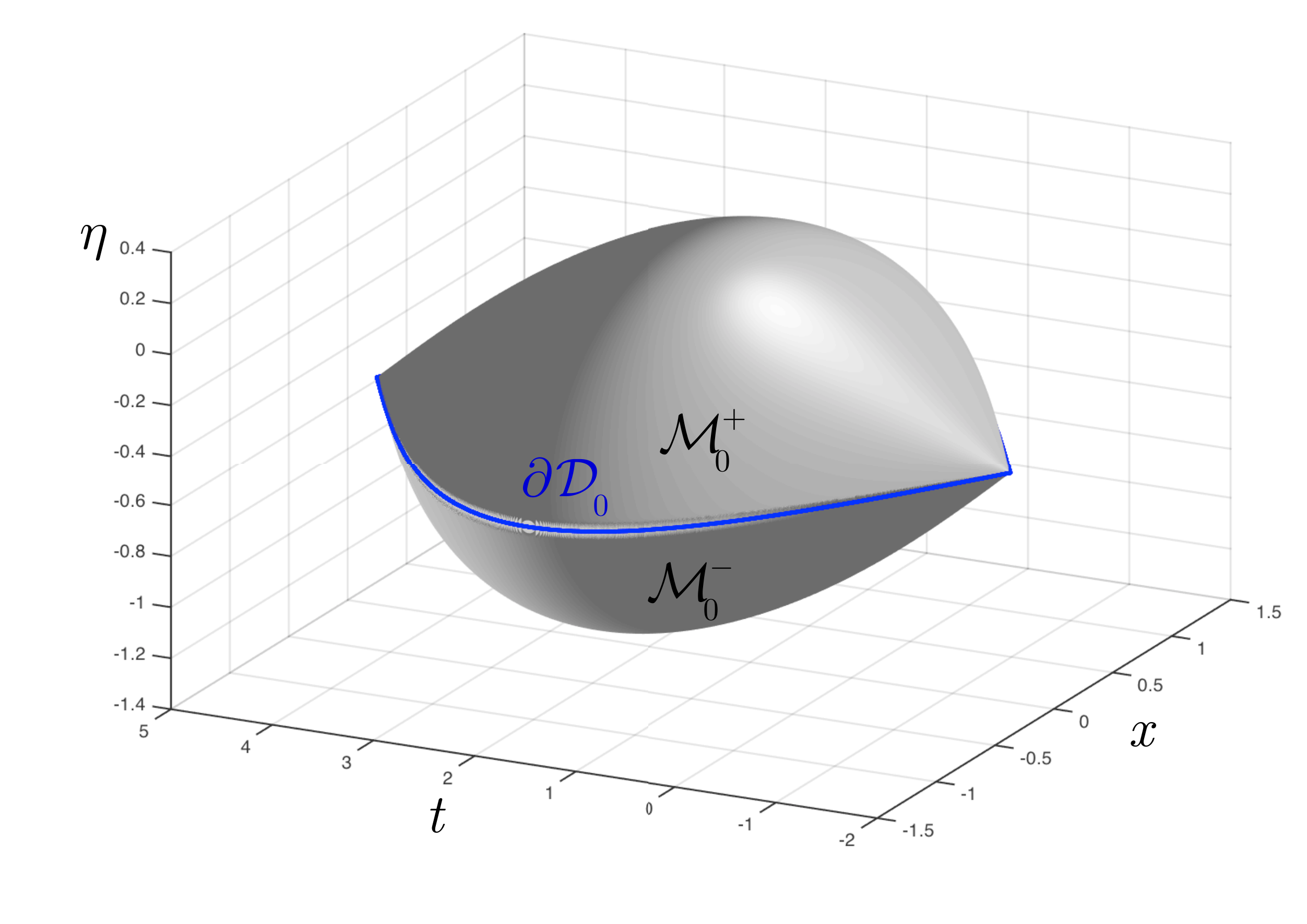}
\par\end{centering}
\caption{The stable critical manifold branch $\mathcal{M}_{0}^{+}$ and the
unstable branch $\mathcal{M}_{0}^{-}$ for the nonlinear mechanical
system \eqref{eq:oscillatory system-1-1-1} with $s=f=1$. Also shown
is the domain boundary $\partial\mathcal{D}_{0}$ along which the
fold in the critical manifold $\mathcal{M}_{0}$ develops.\label{fig:fold}}

\end{figure}
\end{example}

\section{Approximate SFD near equilibria: Static condensation and modal derivatives\label{sec:Approximate-SFD-near-equilbria}}

Here we show that at least two formal reduction procedures used in
structural dynamics, modal condensation and the method of modal derivatives,
can be mathematically justified when the conditions (A1)-(A3) of the
SFD are satisfied. In this case, these two procedures turn out to
provide local first- and second-order approximations, respectively,
to a slow manifold $\mathcal{M}_{\epsilon}$ emanating from an equilibrium
point of the unforced mechanical system

To show this, we also assume the following:
\begin{description}
\item [{(A4)}] \textbf{Independence of critical manifold of the slow velocities:
}The relation
\end{description}
\begin{equation}
\partial_{\dot{x}}Q_{2}(x,\dot{x},\eta,0,t;0)\equiv0,\label{eq:P_2-mod-der}
\end{equation}
holds, i..e, the function $P_{2}$ has not explicit dependence on
the slow velocities $\dot{x}$ for $\epsilon=0$ and $\dot{y}=0$. 

\begin{figure}[h]

\begin{centering}
\includegraphics[width=0.5\textwidth]{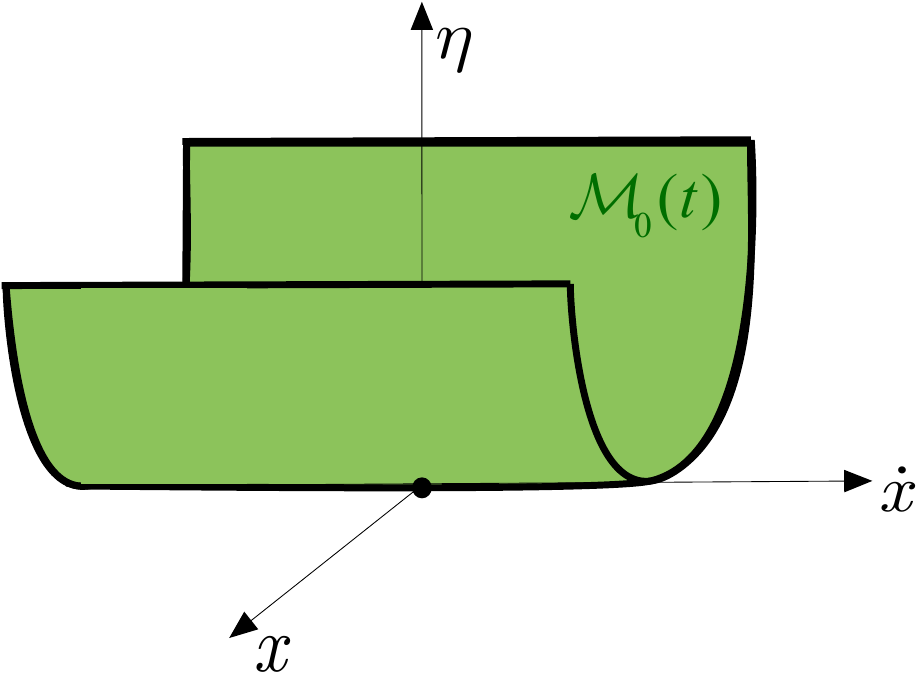}\caption{The geometry of the critical manifold $\mathcal{M}_{0}(t)$ under
assumptions (A3)-(A4) at an arbitrary time $t$.}
\label{fig:flat critical manifold}
\par\end{centering}
\end{figure}

We further assume that the domain $\mathcal{D}_{0}$, over which the
graph $\eta=G_{0}(x,\dot{x},t)$ is defined, contains the line $x=0$
of the $(x,\dot{x},t)$ parameter space, i..e,
\begin{description}
\item [{(A5)}] \textbf{Critical manifold contains an unforced fixed point:}
We assume
\end{description}
\begin{equation}
\left\{ (x,\dot{x},t):\,x=0,\,\,\,\dot{x}=0\right\} \subset\mathcal{D}_{0}.\label{eq:line mod der}
\end{equation}
This condition is satisfied, for instance, when \eqref{eq:systtem00}
is a weakly nonlinear system whose unforced part admits a fixed point
at $q=(x,y)=0.$ The implications of assumptions (A4)-(A5) for the
geometry of the critical manifold are illustrated in Fig. \ref{fig:flat critical manifold}.

Static condensation (Geradin and Rixen \cite{geradin15}) is a linear
reduction procedure applied to a $q=(x,y)$ partition of the degrees
of freedom in system \eqref{eq:systtem00} near an equilibrium point.
In this reduction method, the inertial terms and velocities are simply
ignored in the linearized equation for the $y$ degrees of freedom.
The resulting linear algebraic equation is solved for $y$, and the
result is substituted for $y$ in the $x$ equations, yielding a single
second-order differential equation in the $x$ variables. 

The method of modal derivatives (Idelsohn and Cardona \cite{idelsohn85},
Rutzmoser et al. \cite{rutzmoser14}, Wu and Tiso \cite{wu16}) considers
a similar $q=(x,y)$ partition of coordinates near an unforced equilibrium
and seeks a quadratic invariant manifold tangent to an eigenspace
of the linearized system. The main assumption is that along this quadratic
manifold, the $y$ coordinates can be written as purely quadratic
functions of the $x$ coordinates, with the coefficients of this quadratic
forms collected in an appropriate \emph{modal derivative tensor.}

The above two reduction methods can be justified in our present setting
as follows:
\begin{prop}
\emph{\label{prop:Guyan and modal derivatives}}Under assumptions
(A1)-(A5):
\end{prop}
\begin{description}
\item [{\emph{(i)}}] \emph{The expressions derived for the slow manifold
$\mathcal{M_{\epsilon}}(t)$ in \eqref{eq:slow man coords} satisfy
\begin{equation}
G_{0}(x,\dot{x},t)=\Gamma(t)+\Phi(t)x+\left(\Theta(t)x\right)x+\mathcal{O}\left(\left|x\right|^{3}\right),\label{eq:G0 quadratic}
\end{equation}
where the function $\Gamma(t)$ is the solution of the equation $P_{2}(0,\Gamma(t),0,t;0)=0,$
and the two-tensor $\Phi(t)$ and the three-tensor $\Theta(t)$ satisfy
\begin{eqnarray}
\Phi(t) & = & -\left.\left[\partial_{\eta}P_{2}\right]^{-1}\partial_{x}P_{2}\right|_{x=0,\eta=\Gamma(t),\dot{y}=0,\epsilon=0},\nonumber \\
\Theta(t) & = & -\left.\frac{1}{2}\left[\partial_{\eta}P_{2}\right]^{-1}\left[\partial_{xx}^{2}P_{2}+\left(2\partial_{x\eta}^{2}P_{2}+\partial_{\eta\eta}^{2}P_{2}\Phi(t)\right)\Phi(t)\right]\right|_{x=0,\eta=\Gamma(t),\dot{y}=0,\epsilon=0}.\label{eq:MD-general_form-1}
\end{eqnarray}
}
\item [{\emph{(ii)}}] \emph{Assume that $P_{2}(x,\dot{x},\eta,0,t;0)$
has no explicit time dependence and $(x,y)$ are modal coordinates
for the linearized system at $(x,y)=0$, i.e.,
\[
\partial_{t}P_{2}(x,\dot{x},\eta,0,t;0)\equiv0,\quad\partial_{x}P_{2}(0,0,0,0,t;0)=0,\quad\partial_{\eta}P_{1}(0,0,0,0,t;0)=0.
\]
We then obtain
\begin{equation}
\Gamma=0,\qquad\Phi=0,\qquad\Theta=-\frac{1}{2}\left[\partial_{\eta}P_{2}(0,0,0,0,t;0)\right]^{-1}\partial_{xx}^{2}P_{2}(0,0,0,0,t;0).\label{eq:simple_MD_formulas-2}
\end{equation}
}
\item [{\emph{(iii)}}] \emph{Under the conditions of statement (ii), a
linear-in-$x$ and zeroth-order-in-$\epsilon$ approximation to $\mathcal{M_{\epsilon}}(t)$
yields the modal-condensation-based reduced model
\begin{equation}
\ddot{x}-P_{1}\left(x,\dot{x},0,0,t;0\right)+\mathcal{O}(\epsilon,\left|x\right|^{3})=0\label{eq:Guyan reduced}
\end{equation}
 for the dynamics on $\mathcal{M_{\epsilon}}(t)$ }
\item [{\emph{(iv)}}] \emph{Under the conditions of statement (ii), a quadratic-in-$x$
and zeroth-order-in-$\epsilon$ approximation to $\mathcal{M_{\epsilon}}(t)$
yields the modal-derivatives-based reduced-order model
\begin{equation}
\ddot{x}-P_{1}\left(x,\dot{x},\left(\Theta(t)x\right)x,0,t;0\right)+\mathcal{O}(\epsilon,\left|x\right|^{4})=0\label{eq:modal reduced}
\end{equation}
 for the dynamics on $\mathcal{M_{\epsilon}}(t)$, with $\Theta(t)$
generally referred to as the modal derivative tensor.}
\end{description}
\begin{proof}
See Appendix \ref{app:Guyan and modal derivatives}.
\end{proof}
\begin{rem}
Combining statement (iii) of Theorem \ref{thm:main theorem} with
Proposition \ref{prop:Guyan and modal derivatives} gives that if
$M_{1}(x,\eta,\epsilon)$ is smooth in $\epsilon$ at $\epsilon=0$,
then the static-condensation-based model \eqref{eq:Guyan reduced}
is equivalent to\emph{
\begin{equation}
M_{1}\left(x,0;0\right)\ddot{x}-Q_{1}\left(x,\dot{x},0,0,t;0\right)+\mathcal{O}(\epsilon,\left|x\right|^{3})=0,\label{eq:Guyan reduced-1}
\end{equation}
}and\emph{ }the modal-derivatives-based reduced model \eqref{eq:modal reduced}
is equivalent to 
\begin{equation}
M_{1}\left(x,\left(\Theta x\right)x;0\right)\ddot{x}-Q_{1}\left(x,\dot{x},\left(\Theta x\right)x,0,t;0\right)+\mathcal{O}(\epsilon,\left|x\right|^{4})=0.\label{eq:reduced model 1-1}
\end{equation}
\end{rem}

\begin{rem}
The unevaluated higher-order $\mathcal{O}(\left|x\right|^{3})$ and
$\mathcal{O}(\left|x\right|^{4})$ terms in eqs. \eqref{eq:Guyan reduced-1}
and \eqref{eq:reduced model 1-1} generally do not remain uniformly
small over the full model-reduction domain $\mathcal{D}_{0}.$ Rather,
one can only use the leading-order model terms in these equations
reliably as long as the slow coordinates are rescaled as $x=\sqrt[3]{\epsilon}\xi$
and $x=\sqrt[4]{\epsilon}\xi$, respectively. In that case, \eqref{eq:Guyan reduced-1}
and \eqref{eq:reduced model 1-1} can be re-written as 
\begin{eqnarray}
M_{1}\left(\xi,0;0\right)\ddot{\xi}-Q_{1}\left(\xi,\dot{\xi},0,0,t;0\right)+\mathcal{O}(\epsilon) & = & 0,\label{eq:reduced model 1-1-1}\\
M_{1}\left(\xi,\left(\Theta\xi\right)\xi;0\right)\ddot{\xi}-Q_{1}\left(\xi,\dot{\xi},\left(\Theta\xi\right)\xi,0,t;0\right)+\mathcal{O}(\epsilon) & = & 0,
\end{eqnarray}
respectively. One can then arguably focus on the $\epsilon$-independent
leading order terms for $\epsilon>0$ small enough. The static-condensation-
and model-derivative-based reductions are, therefore, justified in
order $\mathcal{O}(\sqrt[3]{\epsilon})$ and $\mathcal{O}(\sqrt[4]{\epsilon})$
neighborhoods of the $x=0$ equilibrium, respectively, provided that
the assumptions of Proposition \ref{prop:Guyan and modal derivatives}
are satisfied.
\end{rem}
\begin{example}
{[}\emph{Localized reduced-order model for a stiff, weakly nonlinear
system with parametric forcing}{]} We reconsider now the multi-degree-of-freedom
mechanical system \eqref{eq:oscillatory system-1-1} and assume that
\begin{equation}
f_{2}(t)\equiv0,\label{eq:f2 zero assumption}
\end{equation}
i.e., that the external forcing on the stiff degrees of freedom vanishes.
Using the results from Example \ref{ex: Partially-stiff-weakly-1},
we have
\begin{eqnarray*}
P_{1}(x,\dot{x},\eta,\dot{y},t;\epsilon) & = & -M_{1}^{-1}\left[C_{1}\dot{x}+K_{1}x+S_{1}\left(x,\eta\right)-f_{1}(t)\right],\\
P_{2}(x,\dot{x},\eta,\dot{y};\epsilon) & = & -M_{2}^{-1}\left[C_{2}\dot{y}+K_{2}\eta+S_{2}\left(x,\eta\right)\right].
\end{eqnarray*}
As already discussed in Example \ref{ex: Partially-stiff-weakly-1},
conditions (A1)-(A3) are satisfied and hence Theorem \ref{thm:main theorem}
guarantees a slow manifold and determines its reduced dynamics. Assumption
\eqref{eq:P_2-mod-der} is clearly satisfied, as $P_{2}$ does not
depend on $\dot{x}$. Assumption \eqref{eq:line mod der} also holds,
as one sees from the expression for $G_{0}$ in \eqref{eq:G_0 in partially stiff system-1}.
Since $P_{2}$ has no explicit time dependence, the static condensation
and modal derivative formulas in \eqref{eq:simple_MD_formulas} apply
and take the specific form
\begin{equation}
\Gamma\equiv0,\qquad\Phi\equiv0,\qquad\Theta=-\left[-M_{2}^{-1}K_{2}\right]^{-1}\left[-M_{2}^{-1}\partial_{xx}^{2}S_{2}(0,0)\right]=-K_{2}^{-1}\partial_{xx}^{2}S_{2}(0,0).\label{eq:simple_MD_formulas-1}
\end{equation}
Therefore, in a neighborhood of the origin, the reduced-order formulation
\eqref{eq:reduced model 1-1-1} applies and statement (iii) of Proposition
\ref{prop:Guyan and modal derivatives} justifies the static-condensation-based
reduced model
\[
M_{1}\ddot{\xi}+C_{1}\dot{\xi}+K_{1}\xi+S_{1}\left(x,0\right)-f_{1}(t)+\mathcal{O}(\epsilon)=0
\]
as a leading-order reduced model for the dynamics on $\mathcal{M}_{\epsilon}(t)$
in an order $\mathcal{O}(\sqrt[3]{\epsilon})$ neighborhood of the
unforced equilibrium $x=0$. Similarly, statement (iv) of Proposition
\eqref{prop:Guyan and modal derivatives} justifies the modal-derivatives-based
reduced-order model 
\begin{equation}
M_{1}\ddot{\xi}+C_{1}\dot{\xi}+K_{1}\xi+S_{1}\left(x,-\left[K_{2}^{-1}\partial_{xx}^{2}S_{2}(0,0)x\right]x\right)-f_{1}(t)+\mathcal{O}(\epsilon)=0\label{eq:reduced model 1-1-1-1}
\end{equation}
in an order $\mathcal{O}(\sqrt[4]{\epsilon})$ neighborhood of the
unforced equilibrium $x=0$.
\end{example}
The above example illustrates how Proposition \ref{prop:Guyan and modal derivatives}
puts static condensation and modal derivatives in a rigorous context
under appropriate assumptions. We now also illustrate, however, that
these two intuitive reduction methods give incorrect results when
the assumptions of Proposition \ref{prop:Guyan and modal derivatives}
are not satisfied. 
\begin{example}
\emph{\label{ex:reduction failure}{[}Failure of static modal condensation
and model-derivative-based reduction{]} }Consider a two-degree-of-freedom
nonlinear, coupled oscillator system with amplitude-dependent damping
in the first mode, given by the equations 
\begin{eqnarray}
\ddot{x}+\left(c_{1}+\mu_{1}x^{2}\right)\dot{x}+k_{1}x+axy+bx^{3} & = & 0,\qquad x\in\mathbb{R},\nonumber \\
\ddot{y}+c_{2}\dot{y}+k_{2}y+cx^{2} & = & 0,\qquad y\in\mathbb{R}.\label{eq:2DOF example}
\end{eqnarray}
Note that the linearized system at the $(x,y)=(0,0)$ equilibrium
 is in modal coordinates. For $c_{2}>c_{1},$ we obtain slower linear
amplitude decay in the two-dimensional modal subspace of the $x$
variable than in the modal subspace of the $y$ variable. This suggests
a reduction to a model involving only the slower $x$ variables. The
argument used in Example \ref{ex: general oscillatory system}, however,
shows that \eqref{eq:2DOF example} violates assumption (A3) and hence
Proposition \ref{prop:Guyan and modal derivatives} does not apply.
The static condensation procedure nevertheless gives the formal reduced-order
model 
\begin{equation}
\ddot{x}+\left(c_{1}+\mu_{1}x^{2}\right)\dot{x}+k_{1}x+bx^{3}=0,\label{eq:ex guyan reduced}
\end{equation}
and the method of modal-derivates formally gives the formal reduced
model 
\begin{equation}
\ddot{x}+\left(c_{1}+\mu_{1}x^{2}\right)\dot{x}+k_{1}x+\left(b-\frac{ac}{k_{2}}\right)x^{3}=0,\label{eq:ex md reduced}
\end{equation}
modifying \eqref{eq:ex guyan reduced} at cubic order only. While
a global slow manifold is not guaranteed to exist in this example,
a unique, two-dimensional analytic invariant manifold tangent to the
subspace of the $x$ variables at the origin does exist (cf. Haller
and Ponsioen \cite{haller16}). This spectral submanifold (SSM) offers
a mathematically rigorous process for model reduction in system \eqref{eq:2DOF example},
providing an exact reduced flow to which \eqref{eq:ex guyan reduced}
and \eqref{eq:ex md reduced} can be compared. As we show in Appendix
\ref{sec:Appendix:-Details-for-failure}, the reduced model on the
slow SSM is of the form

\noindent 
\begin{eqnarray}
 &  & \ddot{x}+\left[c_{1}+\left(\mu_{1}-\frac{2ac\left(4c_{1}k_{1}+k_{2}\left(c_{1}-c_{2}\right)+2c_{1}c_{2}^{2}-6c_{1}^{2}c_{2}+4c_{1}^{3}\right)}{D}\right)x^{2}\right]\dot{x}\nonumber \\
 &  & +\left[k_{1}-\frac{2ac\left(2c_{1}^{2}-3c_{1}c_{2}+c_{2}^{2}+4k_{1}-k_{2}\right)}{D}\dot{x}^{2}\right]x\label{eq:exact reduced example}\\
 &  & +\left[b-\frac{ac\left(4c_{1}^{4}-6c{}_{1}^{3}c_{2}+2c{}_{1}^{2}c{}_{2}^{2}+5c{}_{1}^{2}k_{2}-c_{1}c_{2}\left(2k_{1}+3k_{2}\right)+2c{}_{2}^{2}k_{1}+8k_{1}^{2}-6k_{1}k_{2}+k{}_{2}^{2}\right)}{D}\right]x^{3}\nonumber \\
 &  & +\,\mathcal{O}(4)=0,\nonumber 
\end{eqnarray}
where
\begin{equation}
D=\left(c_{1}^{2}-c_{1}c_{2}+k_{2}\right)\left(4c_{1}^{2}k_{2}-8c_{1}c_{2}k_{1}-2c_{1}c_{2}k_{2}+4c_{2}^{2}k_{1}+16k_{1}^{2}-8k_{1}k_{2}+k_{2}^{2}\right).\label{eq:Ddef}
\end{equation}
A comparison of the exact reduced model \eqref{eq:exact reduced example}
with the statically condensed version \eqref{eq:ex guyan reduced}
and with the modal-derivatives-based version \eqref{eq:ex md reduced}
shows that the latter two heuristic reduction methods miss most terms
already in the leading-order (cubic) nonlinearities. Depending on
the specific value of the parameters, the missing terms can significantly
impact the nature of the reduced dynamics and hence cannot be omitted.
We note that in the slow-fast limit expressed by the scaling
\begin{equation}
c_{2}\to c_{2}/\epsilon,\qquad k_{2}\to k_{2}/\epsilon^{2},\label{eq:slow fast scaling}
\end{equation}
system \eqref{eq:2DOF example} satisfies the assumptions of Proposition
\ref{prop:Guyan and modal derivatives} and hence the approximation
to the slow SSM should coincide with the approximation to the global
slow manifold $\mathcal{M}_{\epsilon}$ in this case. Indeed, in this
scaling, formulas \eqref{eq:alpha/beta/gamma-1} for the constants
$\alpha,$ $\beta$ and $\gamma$ in Appendix \emph{\ref{sec:Appendix:-Details-for-failure}}
simplify to 
\[
\alpha=-\frac{c}{k_{2}}\epsilon^{2}+\mathcal{O}\left(\epsilon^{3}\right),\qquad\beta=\mathcal{O}\left(\epsilon^{3}\right),\qquad\gamma=\mathcal{O}\left(\epsilon^{3}\right),
\]
and hence the exact reduced model \eqref{eq:exact reduced example}
simplifies to 
\begin{equation}
\ddot{x}+\left(c_{1}+\mu_{1}x^{2}\right)\dot{x}+k_{1}x+\left(b-\frac{ac}{k_{2}/\epsilon^{2}}\right)x^{3}+\mathcal{O}(\epsilon^{3})=0,\label{eq:exact reduced example-1}
\end{equation}
coinciding with the modal-derivatives-based reduced-order model \eqref{eq:ex md reduced}.
This agreement, however, only holds in the slow-fast setting \eqref{eq:slow fast scaling}. 

Even in the conservative limit, when the SSM is replaced by a unique,
analytic Lyapunov-subcenter manifold (Kelley \cite{kelley69}), we
obtain a conservative limit of the exact reduced-order model \eqref{eq:exact reduced example}
in the form
\begin{equation}
\ddot{x}+\left[k_{1}-\frac{2ac}{k_{2}\left(4k_{1}-k_{2}\right)}\dot{x}^{2}\right]x+\left[b-\frac{ac\left(2k_{1}-k_{2}\right)}{k_{2}\left(4k_{1}-k_{2}\right)}\right]x^{3}+\,\mathcal{O}(x^{4})=0,\label{eq:conservative exact}
\end{equation}
filled with nonlinear normal modes (periodic orbits). At the same
time, the conservative limit of the static condensation procedure
gives
\begin{equation}
\ddot{x}+k_{1}x+bx^{3}=0,\label{eq:ex guyan reduced-1}
\end{equation}
while the modal derivatives-based reduction \eqref{eq:ex md reduced}
gives 
\begin{equation}
\ddot{x}+k_{1}x+\left(b-\frac{ac}{k_{2}}\right)x^{3}+\mathcal{O}(x^{4})=0.\label{eq:conservative md}
\end{equation}
Comparing \eqref{eq:conservative exact} and \eqref{eq:conservative md}
shows that the method of modal derivatives gives an incorrect reduced-order
model up to cubic order, unless we have either $a=0$ or $c=0$. As
shown in Fig. \ref{ fig:SSM vs MD near resonance}, the error between
the actual reduced flow \eqref{eq:conservative exact} and \eqref{eq:conservative md}
grows unbounded in the vicinity of the $2:1$ resonance (represented
by $k_{2}=4k_{1})$ between the two natural frequencies of the undamped
limit of system \eqref{eq:2DOF example}. In the limit of an exact
$2:1$ resonance, no invariant manifold tangent to the $x$-subspace
exists, even though the modal derivative approach still suggests the
existence a bounded reduced flow on such a manifold.

\begin{figure}[H]
\begin{centering}
\includegraphics[width=0.8\textwidth]{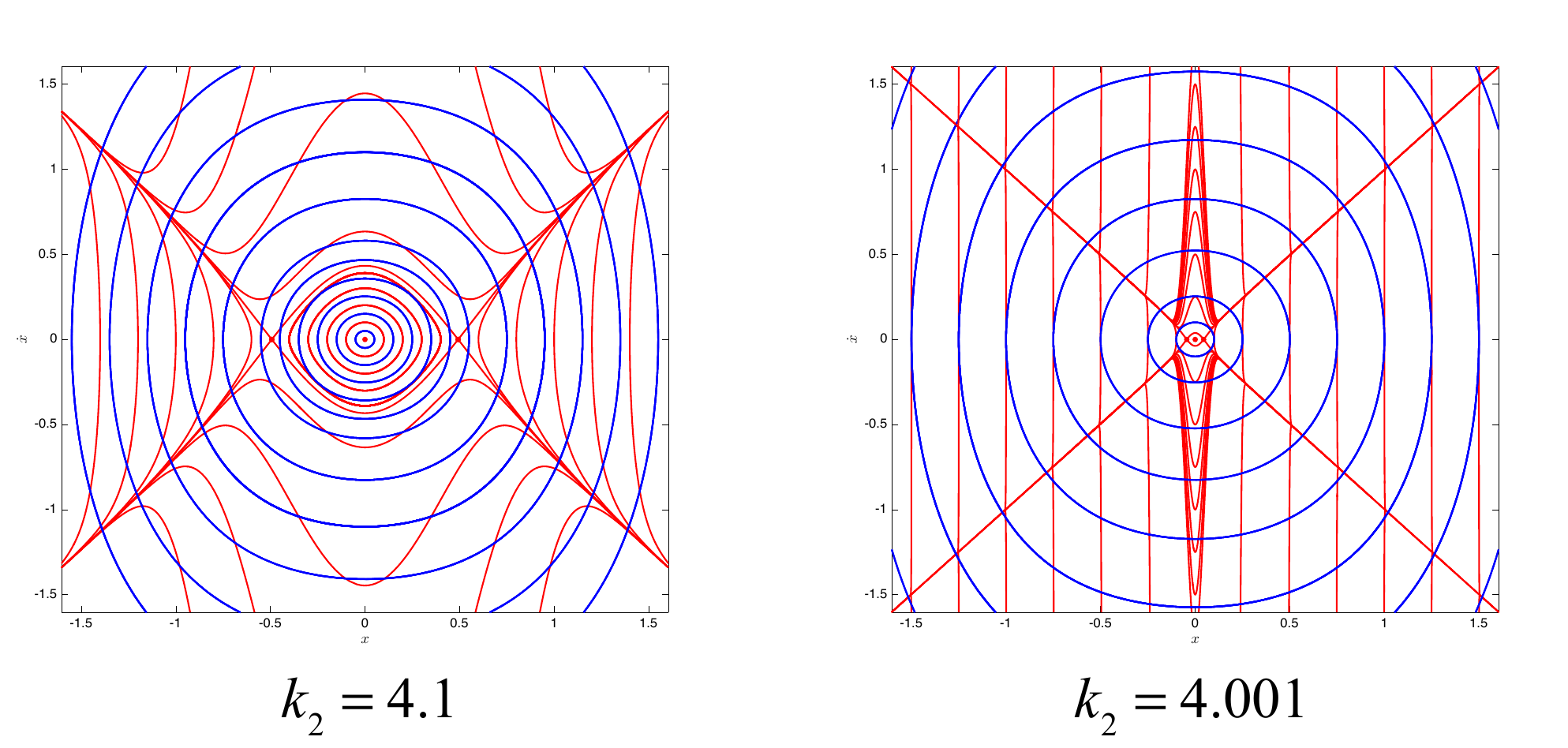}
\par\end{centering}
\caption{Trajectories of the cubic modal-derivatives-based reduction \eqref{eq:conservative md}
(blue) and those of the exact cubic reduction \eqref{eq:conservative exact}
(red) to the unique, 2D analytic invariant manifold over the $(x,\dot{x})$
variables. The remaining parameters are set as $k_{1}=a=b=c=1.$ \label{ fig:SSM vs MD near resonance}}
\end{figure}
\end{example}

\section{A detailed example: Three-degree-of-freedom system with a pendulum
damper}

\begin{figure}[H]
\begin{centering}
\includegraphics[scale=0.3]{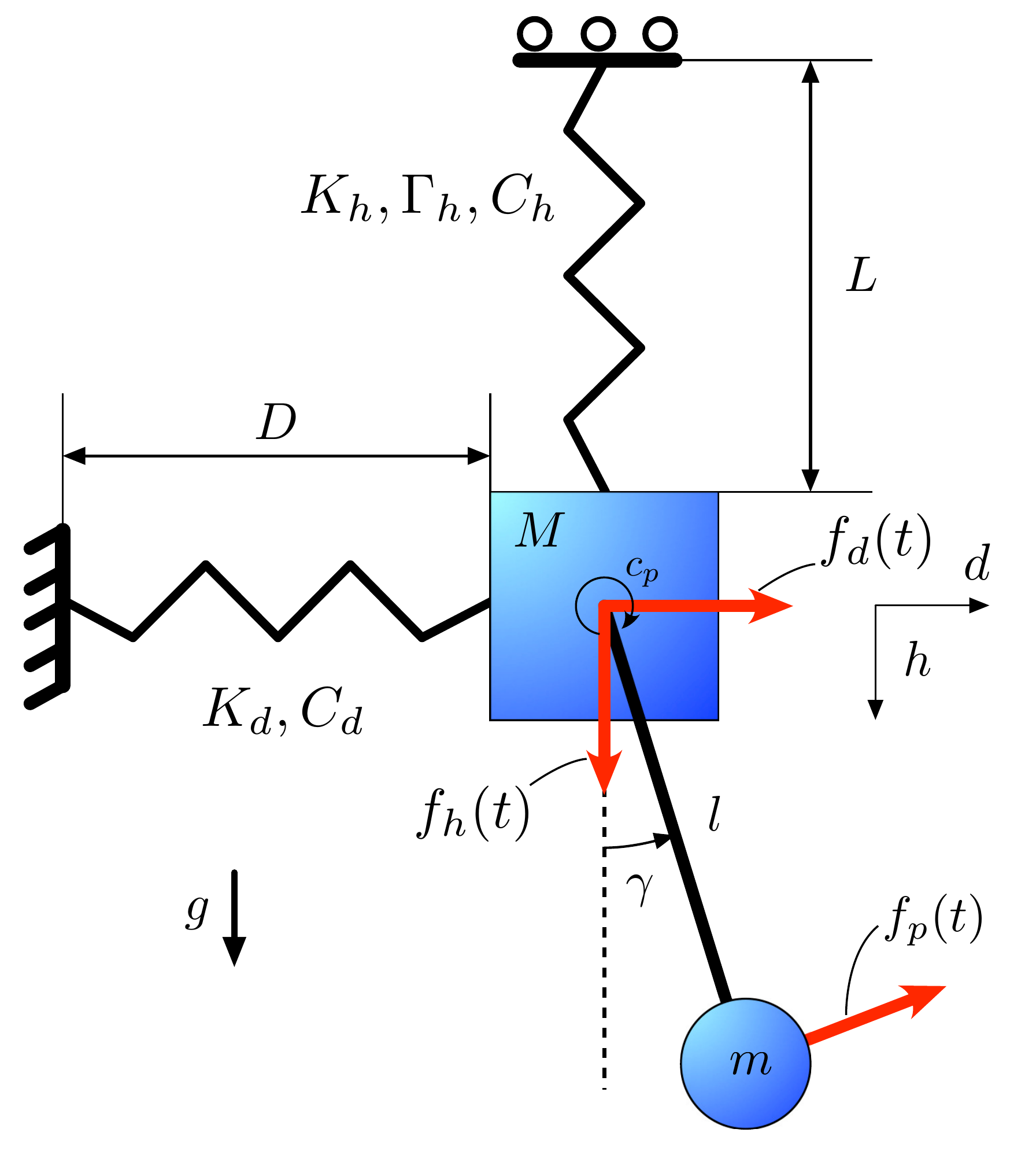}
\par\end{centering}
\caption{Three degree-of-freedom coupled pendulum\label{fig:Three-degree-of-freedom-coupled-1}}
\end{figure}

We consider a the system depicted in Fig. \ref{fig:Three-degree-of-freedom-coupled-1},
with a mass $M$ hanging on a vertical spring of unstretched length
$L$ and linear viscous damping $C_{h}$. The spring is hardening,
with linear stiffness coefficient $K_{h}$ and cubic stiffness coefficient
$\Gamma_{h}>0$. The mass is subject to downward external periodic
forcing of the form $f_{h}(t)=f_{h_{0}}\sin\omega_{1}t$, as well
as to gravity whose constant is $g$. The downward position of the
mass from the unstretched spring position is measured by the coordinate
$h$. The horizontal spring with linear stiffness coefficient $K_{d}$
and natural length $D$ is fixed to the surroundings, thereby introducing
geometric nonlinearities. Added in this direction is a viscous damper
with damping coefficient $C_{d}$ and an external periodic force $f_{d}(t)=f_{d_{0}}\sin\omega_{1}t$,
both acting in the horizontal direction. 

As indicated in Fig. \ref{fig:Three-degree-of-freedom-coupled-1},
a pendulum of mass $m$ and length $l$ is attached to the mass $M$.
The angle of the pendulum from the vertical is denoted by $\gamma$.
The pendulum is also subject to angular viscous damping with coefficient
$c_{p}$, and to an external periodic force $f_{p}(t)=f_{p_{0}}\sin\omega_{2}t$
acting on $m$ in a direction normal to the pendulum . 

The equations of motion for this system are

\begin{eqnarray}
ml^{2}\ddot{\gamma}-ml\sin\gamma\ddot{h}+ml\cos\gamma\ddot{d}+c_{p}\dot{\gamma}+mgl\sin\gamma & = & f_{p}(t)l,\nonumber \\
(M+m)\ddot{h}-ml\sin\gamma\ddot{\gamma}-ml\cos\gamma\dot{\gamma}^{2}+C_{h}\dot{h}+K_{h}h+K_{d}Q(d,h)h+\Gamma_{h}h^{3} & = & (M+m)g+f_{h}(t)-f_{p}(t)\sin\gamma,\nonumber \\
(M+m)\ddot{d}+ml\cos\gamma\ddot{\gamma}-ml\sin\gamma\dot{\gamma}^{2}+C_{d}\dot{d}+K_{d}\left(D+d\right)Q(d,h) & = & f_{d}(t)+f_{p}(t)\cos\gamma,\label{eq:spring-pendulum}
\end{eqnarray}
with the geometric nonlinear term $Q(d,h)$ 

\[
Q(d,h)=\left(1-\frac{D}{\sqrt{\left(D+d\right)^{2}+h^{2}}}\right).
\]
The linearized oscillation frequencies of the uncoupled springs and
of the pendulum are
\begin{equation}
\omega_{h}=\sqrt{\frac{K_{h}}{M}},\qquad\omega_{d}=\sqrt{\frac{K_{d}}{M}},\qquad\omega_{p}=\sqrt{\frac{g}{l},}\label{eq:lineigen-1-1}
\end{equation}
respectively. With the help of these frequencies, we non-dimensionalize
the $h$ and $d$ coordinates, the time $t$, and all system parameters
by letting 

\[
\tilde{h}=\frac{h}{L},\quad\tilde{d}=\frac{d}{D},\quad\tilde{t}=\omega_{p}t,
\]

\[
\Delta=\frac{l}{L},\quad\rho=\frac{D}{L},\quad\beta=\frac{m}{M},\quad F_{h}(t)=\frac{f_{h}(t)}{Mg},\quad F_{p}(t)=\frac{f_{p}(t)}{Mg},\quad F_{d}(t)=\frac{f_{d}(t)}{Mg},\quad G_{p}(t)=\frac{f_{p}(t)}{mg},
\]

\[
\pi_{h}=\frac{C_{h}}{\omega_{p}M},\quad\pi_{d}=\frac{C_{d}}{\omega_{p}M},\quad\pi_{p}=\frac{c_{p}}{\omega_{p}mL^{2}},\quad q_{h}=\frac{\omega_{h}^{2}}{\omega_{p}^{2}},\quad q_{d}=\frac{\omega_{d}^{2}}{\omega_{p}^{2}},\quad a_{h}=\frac{\Gamma_{h}L^{2}}{M\omega_{p}^{2}},
\]
which leads to the following definition for the scaled version of
$Q(d,h)$

\[
\tilde{Q}(\tilde{d},\tilde{h})=\left(1-\frac{\rho}{\sqrt{\rho^{2}\left(1+\tilde{d}\right)^{2}+\tilde{h}^{2}}}\right).
\]
Denoting differentiation with respect to the new time $\tilde{t}$
still by a dot, then dropping all the tildes, we obtain the non-dimensionalized
equations of motions

\begin{eqnarray}
\Delta^{2}\ddot{\gamma}-\Delta\sin\gamma\ddot{h}+\rho\Delta\cos\gamma\ddot{d}+\pi_{p}\dot{\gamma}+\Delta^{2}\sin\gamma & = & \Delta^{2}G_{p}(t),\nonumber \\
(1+\beta)\ddot{h}-\beta\Delta\sin\gamma\ddot{\gamma}-\beta\Delta\cos\gamma\dot{\gamma}^{2}+\pi_{h}\dot{h}+q_{h}h+q_{d}hQ(d,h)+a_{h}h^{3} & = & (1+\beta)\Delta+F_{h}(t)\Delta-F_{p}(t)\Delta\sin\gamma,\nonumber \\
(1+\beta)\ddot{d}+\beta\frac{\Delta}{\rho}\cos\gamma\ddot{\gamma}-\beta\frac{\Delta}{\rho}\sin\gamma\dot{\gamma}^{2}+\pi_{d}\dot{d}+q_{d}\left(1+d\right)Q(d,h) & = & F_{d}(t)\frac{\Delta}{\rho}+F_{p}(t)\frac{\Delta}{\rho}\cos\gamma.\nonumber \\
\label{eq:nondim_pendeq}
\end{eqnarray}

\subsection{Two soft degrees of freedom\label{subsec:Two-slow-degrees}}

We are interested in applying the SFD procedure to system \eqref{eq:nondim_pendeq}
to obtain an exact reduced-order model for the dynamics. First, we
assume that $h$ is a stiff degree of freedom and $(\gamma,d)$ represent
soft degrees of freedom. In that case, using the notation from system
\eqref{eq:systtem00}, we can write the mass matrix $M(q,t;\epsilon)$
and the forcing term $F(q,\dot{q},t;\epsilon)$ as

\begin{eqnarray*}
M(q,t;\epsilon) & = & \left(\begin{array}{ccc}
\Delta^{2} & \rho\Delta\cos x_{\gamma} & -\Delta\sin x_{\gamma}\\
\beta\frac{\Delta}{\rho}\cos x_{\gamma} & 1+\beta & 0\\
-\beta\Delta\sin x_{\gamma} & 0 & 1+\beta
\end{array}\right),\\
F(q,\dot{q},t;\epsilon) & = & \left(\begin{array}{l}
-\pi_{p}\dot{x}_{\gamma}-\Delta^{2}\sin x_{\gamma}+\Delta^{2}G_{p}(t)\vspace{9pt}\\
\beta\frac{\Delta}{\rho}\sin x_{\gamma}\dot{x}_{\gamma}^{2}-\pi_{d}\dot{x}_{d}-q_{d}\left(1+x_{d}\right)Q(x_{d},\frac{y}{\epsilon})+F_{d}(t)\frac{\Delta}{\rho}+F_{p}(t)\frac{\Delta}{\rho}\cos x_{\gamma}\vspace{5pt}\\
\beta\Delta\cos x_{\gamma}\dot{x}_{\gamma}^{2}-\pi_{h}\dot{y}_{h}-q_{h}\epsilon\frac{y}{\epsilon}-q_{d}\epsilon\frac{y}{\epsilon}Q(x_{d},\frac{y}{\epsilon})-a_{h}\epsilon^{3}\left(\frac{y}{\epsilon}\right)^{3}\\
+(1+\beta)\Delta+F_{h}(t)\Delta-F_{p}(t)\Delta\sin x_{\gamma}
\end{array}\right).
\end{eqnarray*}
Here we have introduced the coordinates $(x,y)$ by letting

\[
x_{\gamma}=\gamma,\quad x_{d}=d,\quad y=h.
\]
The modified mass matrices $M_{i}$ and the forcing terms $Q_{i}$
defined in \eqref{eq:M_i and Q_i def} take the specific form
\begin{eqnarray*}
M_{1} & = & M_{11}-M_{12}M_{22}^{-1}M_{21}=\left(\begin{array}{cc}
\frac{\Delta^{2}}{1+\beta}\left(1+\beta\cos^{2}x_{\gamma}\right) & \rho\Delta\cos x_{\gamma}\\
\beta\frac{\Delta}{\rho}\cos x_{\gamma} & 1+\beta
\end{array}\right),\\
M_{2} & = & M_{22}-M_{21}M_{11}^{-1}M_{12}=\frac{1+\beta}{1+\beta\sin^{2}x_{\gamma}},\\
Q_{1} & = & F_{1}-M_{12}M_{22}^{-1}F_{2}=\left[\begin{array}{c}
q_{1}\\
q_{2}
\end{array}\right]\\
q_{1} & = & -\pi_{p}\dot{x}_{\gamma}-\Delta^{2}\sin x_{\gamma}+\Delta^{2}G_{p}(t)+\frac{\Delta\sin x_{\gamma}}{1+\beta}\bigl[\beta\Delta\cos x_{\gamma}\dot{x}_{\gamma}^{2}-\pi_{h}\dot{y}-q_{h}\epsilon\frac{y}{\epsilon}\\
 &  & -q_{d}\epsilon\frac{y}{\epsilon}Q(x_{d},\frac{y}{\epsilon})-a_{h}\epsilon^{3}\left(\frac{y}{\epsilon}\right)^{3}+(1+\beta)\Delta+F_{h}(t)\Delta-F_{p}(t)\Delta\sin x_{\gamma}\bigr]\\
q_{2} & = & \beta\frac{\Delta}{\rho}\sin x_{\gamma}\dot{x}_{\gamma}^{2}-\pi_{d}\dot{x}_{d}-q_{d}\left(1+x_{d}\right)Q(x_{d},\frac{y}{\epsilon})+F_{d}(t)\frac{\Delta}{\rho}+F_{p}(t)\frac{\Delta}{\rho}\cos x_{\gamma}\\
Q_{2} & = & F_{2}-M_{21}M_{11}^{-1}F_{1}\\
 & = & \beta\Delta\cos x_{\gamma}\dot{x}_{\gamma}^{2}-\pi_{h}\dot{y}-q_{h}\epsilon\frac{y}{\epsilon}-q_{d}\epsilon\frac{y}{\epsilon}Q(x_{d},\frac{y}{\epsilon})-a_{h}\epsilon^{3}\left(\frac{y}{\epsilon}\right)^{3}+(1+\beta)\Delta+F_{h}(t)\Delta-F_{p}(t)\Delta\sin x_{\gamma}\\
 &  & +\frac{\left(1+\beta\right)\beta\sin x_{\gamma}}{\Delta\left(1+\beta\sin^{2}x_{\gamma}\right)}\left[-\pi_{p}\dot{x}_{\gamma}-\Delta^{2}\sin x_{\gamma}+\Delta^{2}G_{p}(t)\right]\\
 &  & -\frac{\beta\rho\sin x_{\gamma}\cos x_{\gamma}}{1+\beta\sin^{2}x_{\gamma}}\left[\beta\frac{\Delta}{\rho}\sin x_{\gamma}\dot{x}_{\gamma}^{2}-\pi_{d}\dot{x}_{d}-q_{d}\left(1+x_{d}\right)Q(x_{d},\frac{y}{\epsilon})+F_{d}(t)\frac{\Delta}{\rho}+F_{p}(t)\frac{\Delta}{\rho}\cos x_{\gamma}\right].
\end{eqnarray*}
These give the following expression for the function $P_{1}$ 

\begin{align}
 & P_{1}\left(x,v,\eta,w,t;\epsilon\right)=M_{1}^{-1}\left[\begin{array}{l}
q_{1s}\\
q_{2s}
\end{array}\right],\label{eq:P1_slow_slow}
\end{align}

\begin{align*}
q_{1s} & =-\pi_{p}v_{\gamma}-\Delta^{2}\sin x_{\gamma}+\Delta^{2}G_{p}(t)+\frac{\Delta\sin x_{\gamma}}{1+\beta}\bigl[\beta\Delta\cos x_{\gamma}v_{\gamma}^{2}-\pi_{h}w_{h}-q_{h}\epsilon\eta\\
 & -q_{d}\epsilon\eta Q(x_{d},\eta)-a_{h}\epsilon^{3}\eta^{3}+(1+\beta)\Delta+F_{h}(t)\Delta-F_{p}(t)\Delta\sin x_{\gamma}\bigr],\\
q_{2s} & =\beta\frac{\Delta}{\rho}\sin x_{\gamma}v_{\gamma}^{2}-\pi_{d}v_{d}-q_{d}\left(1+x_{d}\right)Q(x_{d},\eta)+F_{d}(t)\frac{\Delta}{\rho}+F_{p}(t)\frac{\Delta}{\rho}\cos x_{\gamma},
\end{align*}
with the inverse of $M_{1}$ given by

\begin{eqnarray*}
M_{1}^{-1} & = & \frac{1}{\Delta^{2}}\left[\begin{array}{cc}
1+\beta & -\rho\Delta\cos x_{\gamma}\\
-\beta\frac{\Delta}{\rho}\cos x_{\gamma} & \frac{\Delta^{2}}{1+\beta}\left(1+\beta\cos^{2}x_{\gamma}\right)
\end{array}\right].
\end{eqnarray*}
The function $P_{2}$ takes the specific form

\begin{eqnarray}
P_{2}\left(x,v,\eta,w,t;\epsilon\right) & =\epsilon & \left(\frac{1+\beta\sin^{2}x_{\gamma}}{1+\beta}\right)\biggl(\beta\Delta\cos x_{\gamma}v_{\gamma}^{2}-\pi_{h}w_{h}-q_{h}\epsilon\eta\nonumber \\
 &  & -q_{d}\epsilon\eta Q(x_{d},\eta)-a_{h}\epsilon^{3}\eta^{3}+(1+\beta)\Delta+F_{h}(t)\Delta-F_{p}(t)\Delta\sin x_{\gamma}\nonumber \\
 &  & +\frac{\left(1+\beta\right)\beta\sin x_{\gamma}}{\Delta\left(1+\beta\sin^{2}x_{\gamma}\right)}\left[-\pi_{p}v_{\gamma}-\Delta^{2}\sin x_{\gamma}+\Delta^{2}G_{p}(t)\right]\nonumber \\
 &  & -\frac{\beta\rho\sin x_{\gamma}\cos x_{\gamma}}{1+\beta\sin^{2}x_{\gamma}}\biggl[\beta\frac{\Delta}{\rho}\sin x_{\gamma}v_{\gamma}^{2}-\pi_{d}v_{d}-q_{d}\left(1+x_{d}\right)Q(x_{d},\eta)+F_{d}(t)\frac{\Delta}{\rho}+F_{p}(t)\frac{\Delta}{\rho}\cos x_{\gamma}\biggr]\biggr).\nonumber \\
\label{eq:P2_slow_slow}
\end{eqnarray}

We observe that $\lim_{\epsilon\to0}P_{2}\left(x,v,\eta,w,t;\epsilon\right)\equiv0,$
and therefore assumptions (A1)-(A3) are not satisfied without further
assumptions on the parameters that ensure the stiff-soft partition
of the coordinates. To this end we express the stiffness of the $y$
degree of freedom by letting

\[
\Delta=\frac{l}{L}=\frac{\delta}{\epsilon},\quad\rho=\frac{D}{L}=\frac{\phi}{\epsilon},\quad q_{d}=\frac{\omega_{d}^{2}}{\omega_{p}^{2}}=\Omega_{d}^{2},\quad q_{h}=\frac{\omega_{h}^{2}}{\omega_{p}^{2}}=\frac{\Omega_{h}^{2}}{\epsilon^{2}},\quad a_{h}=\frac{\alpha_{h}}{\epsilon^{4}},
\]
\begin{align}
\pi_{h} & =\frac{C_{h}}{\omega_{p}M}=\frac{\mu_{h}}{\epsilon},\quad\pi_{d}=\frac{C_{d}}{\omega_{p}M}=\mu_{d},\quad\pi_{p}=\frac{c_{p}}{\omega_{p}mL^{2}}=\frac{c_{p}}{\omega_{p}m\left(\frac{\epsilon}{\delta}l\right)^{2}}=\frac{c_{p}\delta^{2}}{\omega_{p}m\epsilon^{2}l^{2}}=\frac{\mu_{p}}{\epsilon^{2}}.\nonumber \\
\label{eq:rescalings two slow}
\end{align}
In this parameter range, assumptions (A1)-(A3) are satisfied, as we
show in Appendix \ref{sec:Appendix:-Details-for-fast-slow-slow system}.
The reduced model arising from these calculations is of the form

\begin{align}
\ddot{\gamma} & =\frac{\omega_{p}^{2}\left(M+m\right)}{M+m\sin^{2}\gamma}\left(-\frac{c_{p}}{\omega_{p}^{2}ml^{2}}\dot{\gamma}-\sin\gamma+\frac{f_{p}(t)}{mg}\right)\label{eq:red_mod_1}\\
 & -\frac{\omega_{p}^{2}M\cos\gamma}{M+m\sin^{2}\gamma}\left(\frac{m}{M\omega_{p}^{2}}\sin\gamma\dot{\gamma}^{2}-\frac{C_{d}}{\omega_{p}^{2}Ml}\dot{d}-\frac{K_{d}}{Mg}d+\frac{f_{d}(t)}{Mg}+\frac{f_{p}(t)}{Mg}\cos\gamma\right)+\mathcal{O}(\epsilon),\nonumber 
\end{align}

\begin{align}
\ddot{d} & =\frac{\omega_{p}^{2}DM}{M+m\sin^{2}\gamma}\left(\frac{ml}{MD\omega_{p}^{2}}\sin\gamma\dot{\gamma}^{2}-\frac{C_{d}}{MD\omega_{P}^{2}}\dot{d}-\frac{K_{d}l}{MgD}d+\frac{f_{d}(t)l}{MgD}+\frac{f_{p}(t)l}{MgD}\cos\gamma\right)\label{eq:red_mod_2}\\
 & -\frac{\omega_{p}^{2}Dm\cos\gamma}{M+m\sin^{2}\gamma}\left(-\frac{c_{p}}{\omega_{p}^{2}mDl}\dot{\gamma}-\frac{l}{D}\sin\gamma+\frac{l}{D}\frac{f_{p}(t)}{mg}\right)+\mathcal{O}(\epsilon).\nonumber 
\end{align}

\begin{figure}[H]
\centering{}\includegraphics[width=0.4\textwidth]{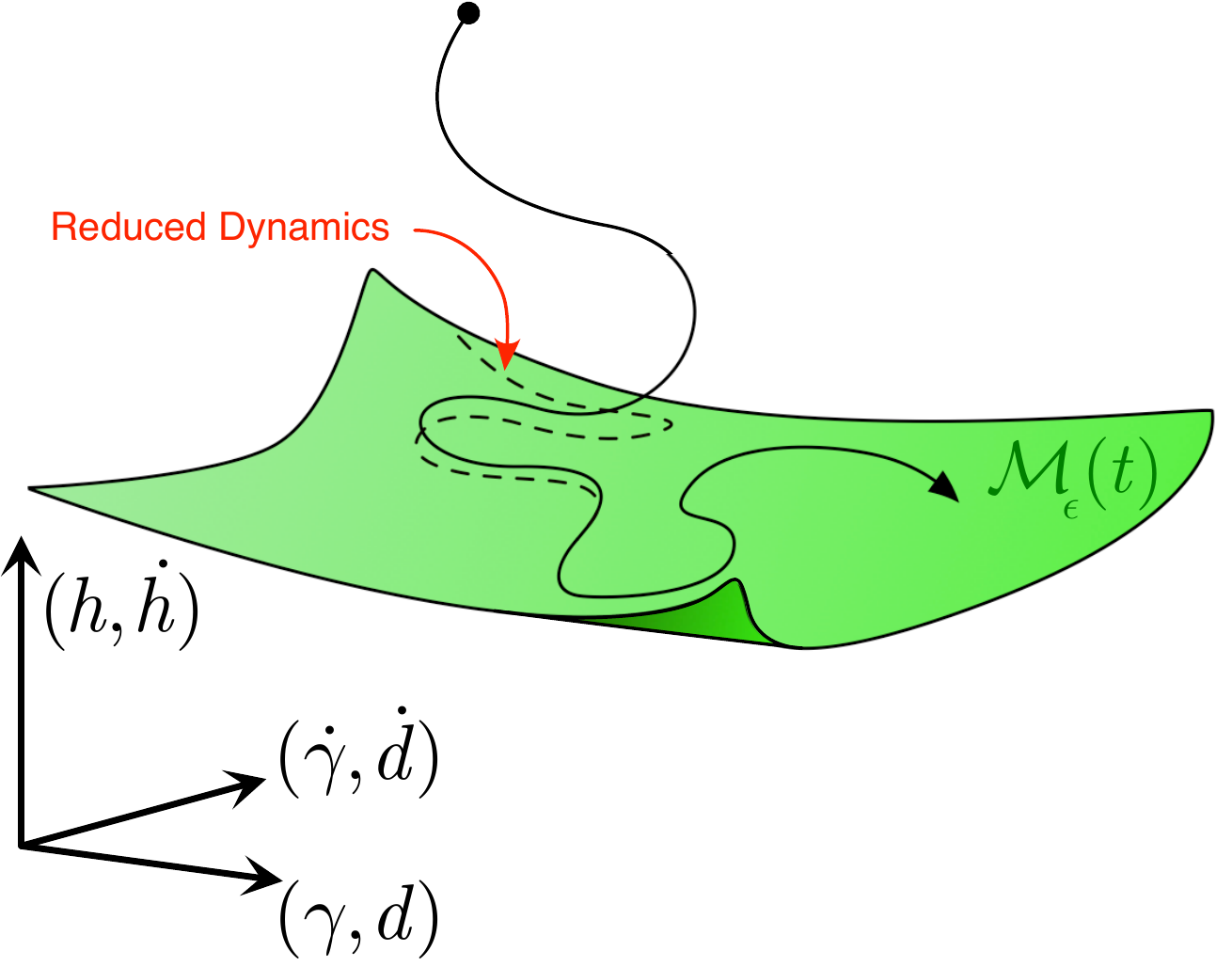}\caption{Illustration of the attracting slow manifold $\mathcal{M_{\epsilon}}$
for the mechanical system \eqref{eq:nondim_pendeq}, graphed over
the two slow degrees of freedom $x_{\gamma}$ and $x_{d}$ and their
corresponding velocities. A general trajectory $q(t)$ is attracted
to the slow manifold, synchronizing exponentially fast with a trajectory
of reduced dynamics (dashed line). \label{fig:slow_manifold}}
\end{figure}

We have implemented this model in \textit{Mathematica} to show how
a general trajectory $x(t)$ of the full system is attracted to reduced
model-trajectories the slow manifold $\mathcal{M}_{\epsilon}$. A
graphical illustration of this behavior is shown in Fig. \ref{fig:slow_manifold}.

\begin{figure}[H]
\begin{centering}
\subfloat[\label{fig:gamma_t}]{\begin{centering}
\includegraphics[scale=0.35]{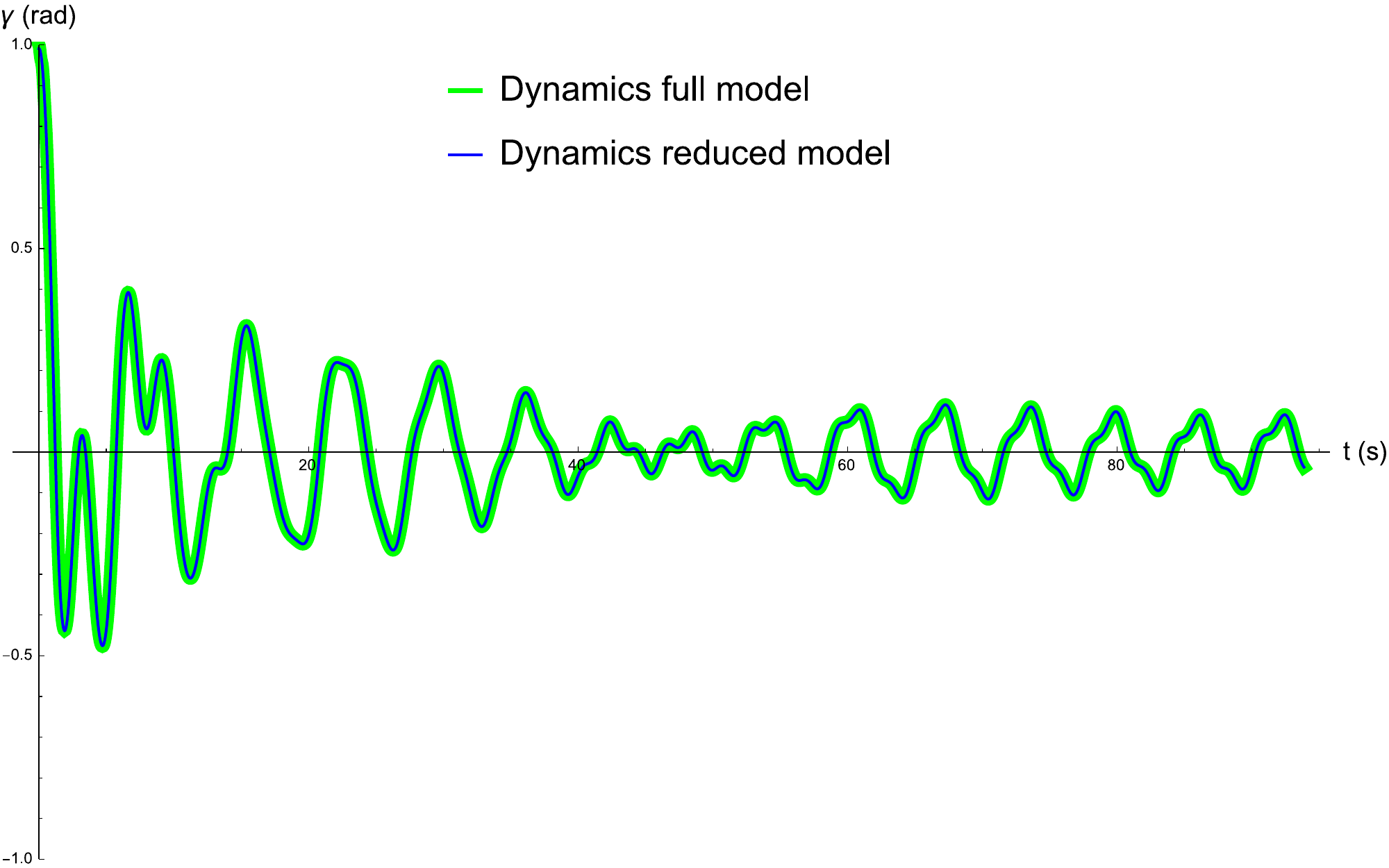}
\par\end{centering}
}\hfill{}\subfloat[\label{fig:gammad_t}]{\centering{}\includegraphics[scale=0.35]{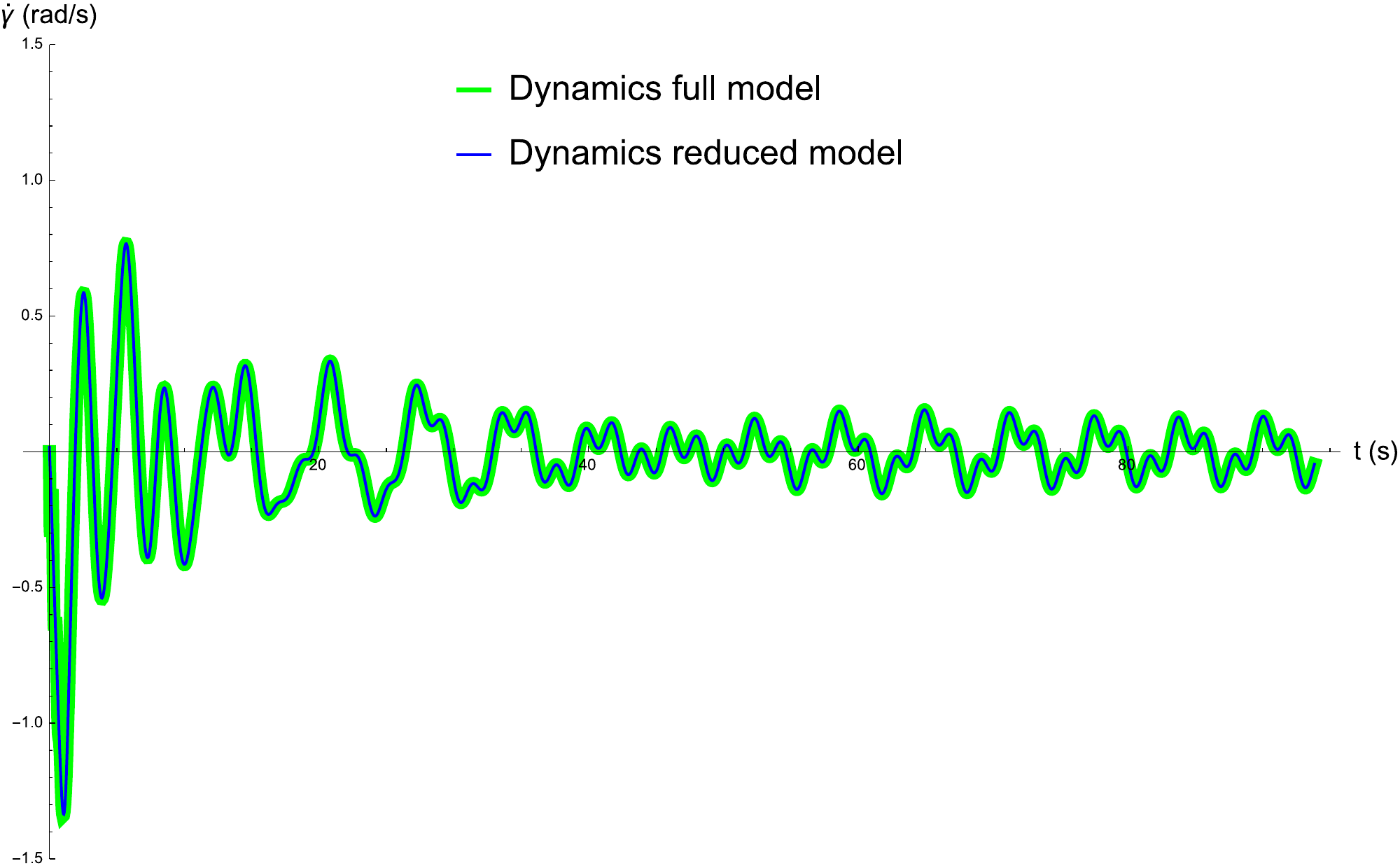}}
\par\end{centering}
\caption{Exponentially fast synchronization of the soft $(\gamma,\dot{\gamma})$
coordinates of the full trajectory and of a reduced model trajectory.
\label{fig:minipage_gamma}}
\end{figure}

\begin{figure}[H]
\begin{centering}
\subfloat[\label{fig:d_t}]{\begin{centering}
\includegraphics[scale=0.35]{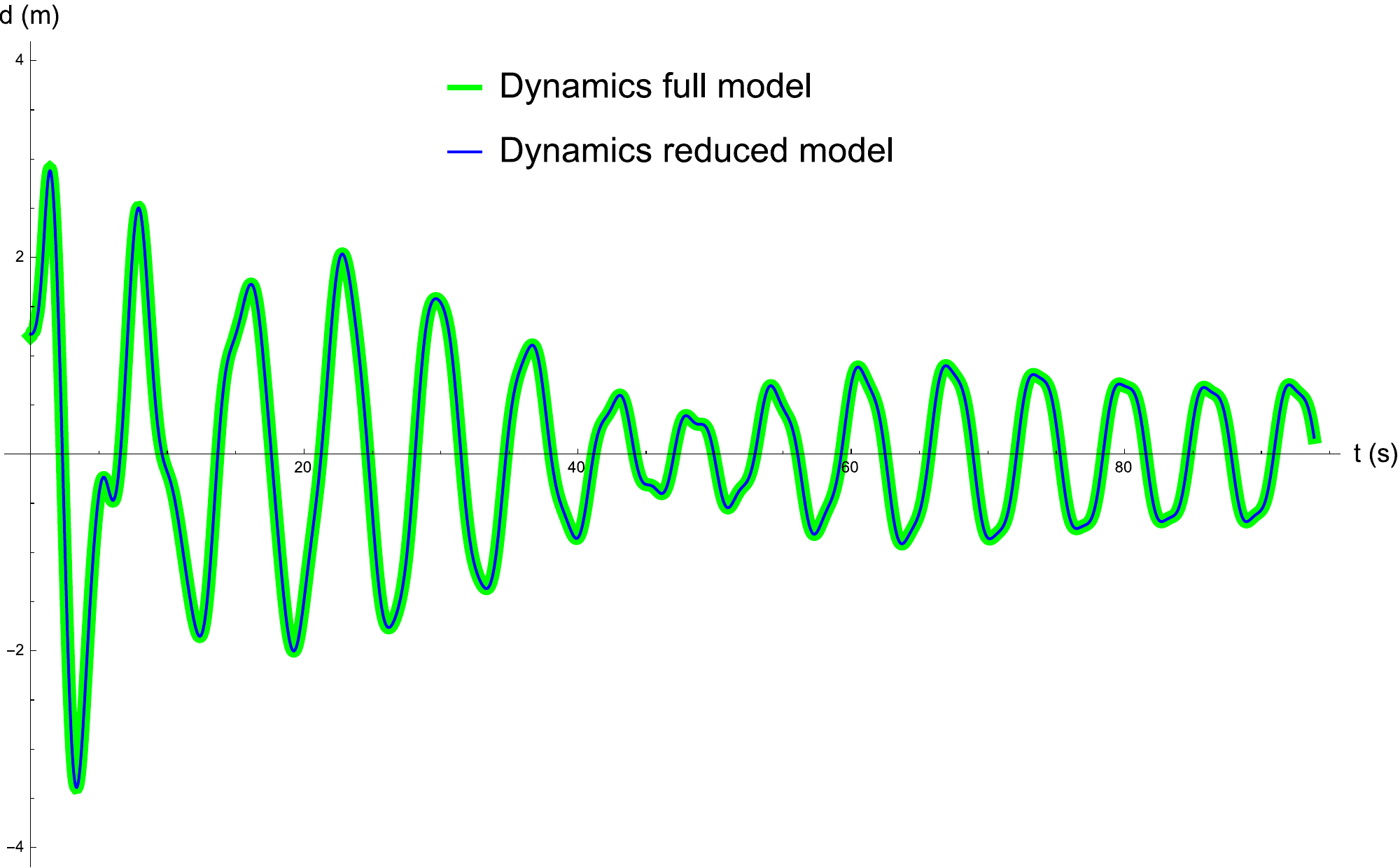}
\par\end{centering}
}\hfill{}\subfloat[\label{fig:dd_t}]{\centering{}\includegraphics[scale=0.35]{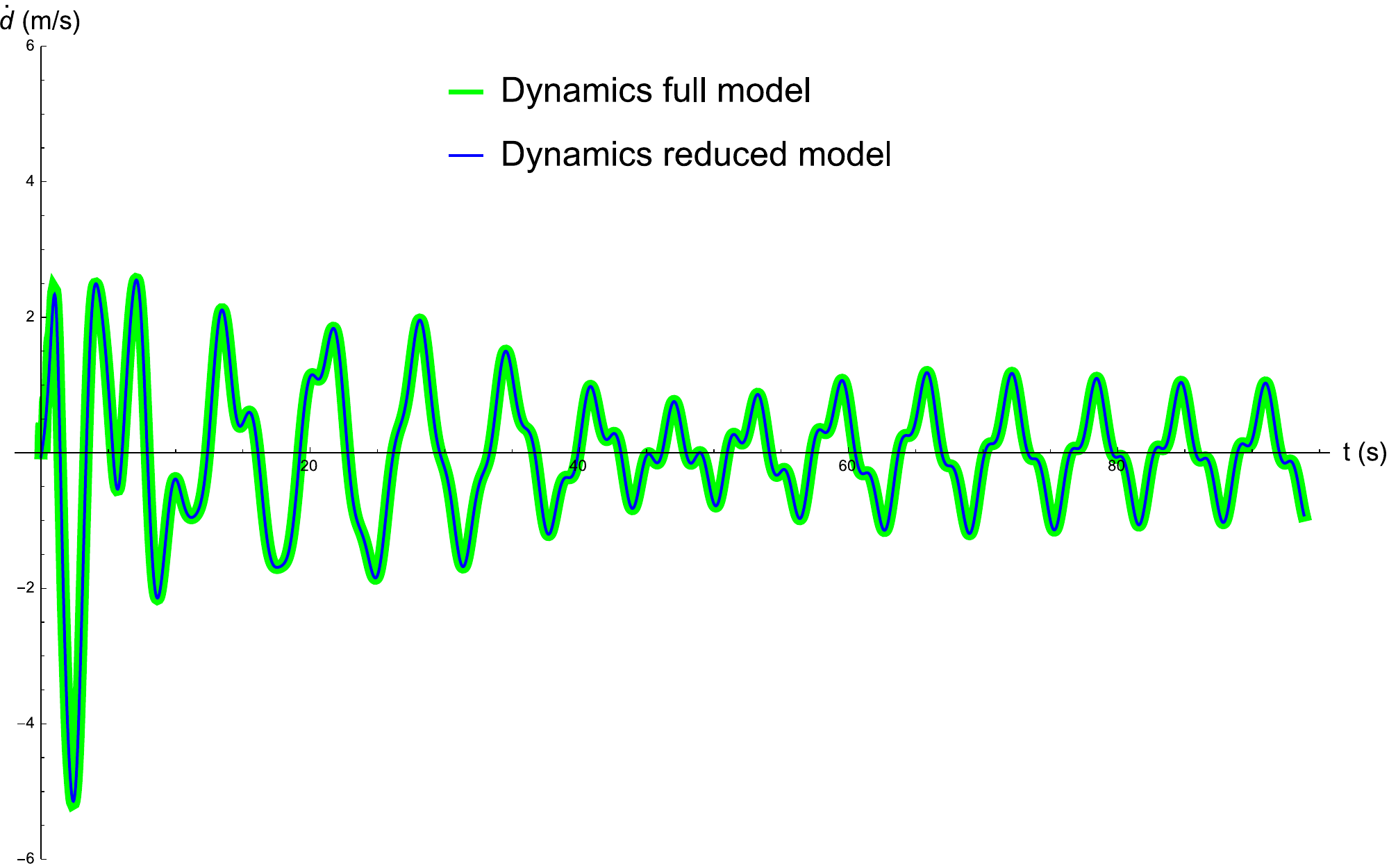}}
\par\end{centering}
\caption{Exponentially fast synchronization of the soft $(d,\dot{d})$ coordinates
of the full trajectory and of a reduced model trajectory.\label{fig:minipage_d}}
\end{figure}

\begin{figure}[H]
\begin{centering}
\includegraphics[width=0.8\textwidth]{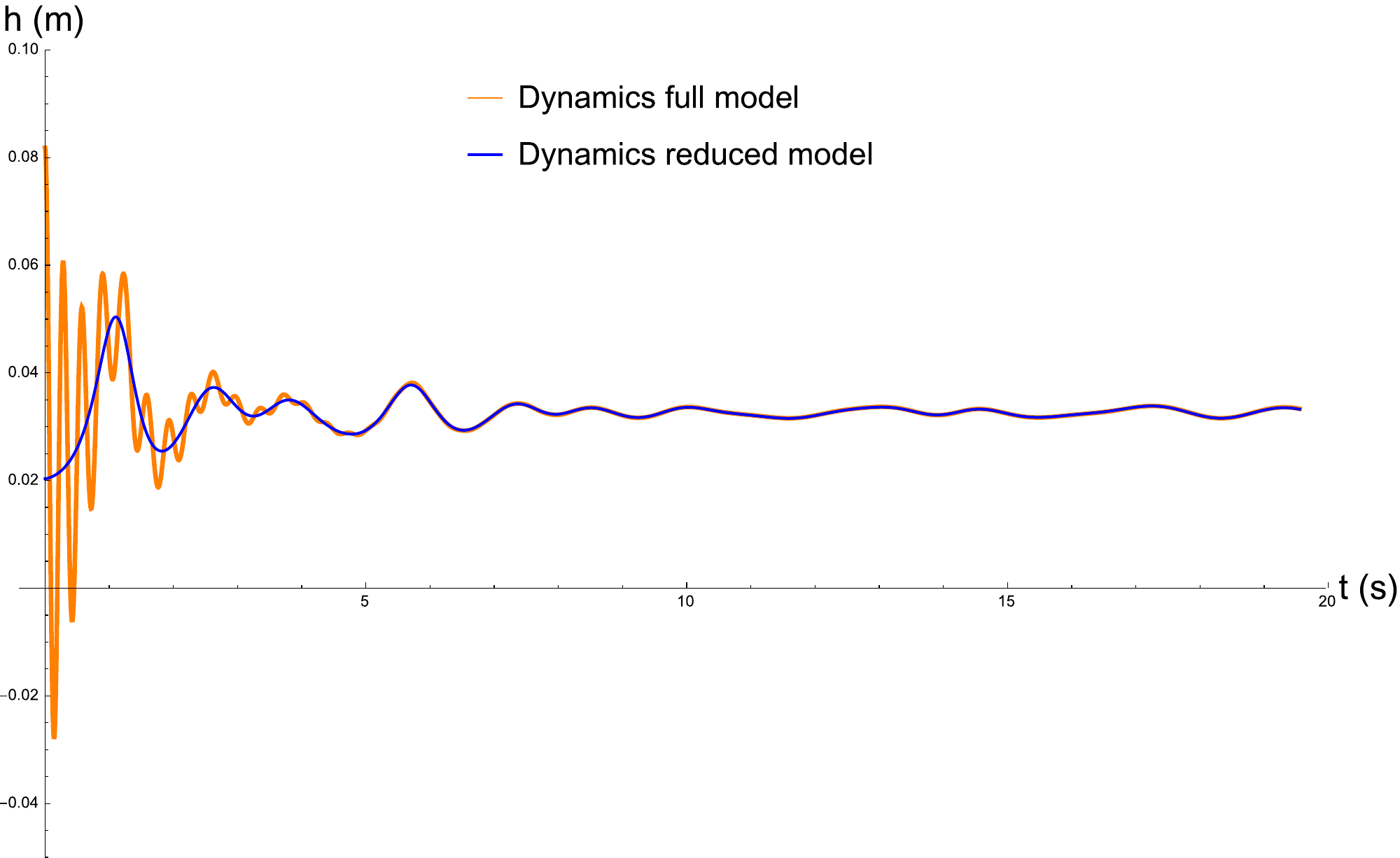}
\par\end{centering}
\caption{Exponentially fast convergence of the fast coordinate $h$ along the
full trajectory to the same coordinate along a trajectory of the reduced
system. \label{fig:minipage_h}}
\end{figure}

For a numerical illustration of the accuracy of the reduced model,
we choose the following values for the system parameters:

\begin{gather*}
l=D=6\text{ m},\quad L=1\text{ m},\quad M=m=1\text{ kg},\\
K_{h}=600\text{ N/m},\quad\Gamma_{h}=0.5\text{ N/m}^{3},\quad K_{d}=2\text{ N/m},\\
C_{d}=0.33\cdot\omega_{p}\cdot M\text{ kg/s},\quad C_{h}=3\cdot\omega_{p}\cdot M\text{ kg/s},\\
c_{p}=0.33\cdot\omega_{p}\cdot m\cdot L^{2}\text{ }(\text{kg}\cdot\text{m}^{2})\text{/s},\quad g=9.81\text{ m/s}^{2},\\
f_{p}(t)=0.5\cdot\sin(t)\text{ N},\quad f_{h}(t)=f_{d}(t)=0.5\cdot\sin(3t)\text{ N},\\
\epsilon=1\cdot10^{-8}.
\end{gather*}
We give the full system the initial condition

\[
q_{0}=(\gamma_{0},d_{0},\dot{\gamma}_{0},\dot{d}_{0},h_{0},\dot{h}_{0})=(1.000,1.200,0.000,0.000,0.08182,0.005301),
\]
which lies off the slow manifold $\mathcal{M}_{\epsilon}$, then integrate
the trajectory starting from this initial condition in forward time.
We track the Euclidean distance between the fast variables ($h(t)$,
$\dot{h}(t)$) and the explicitly computable slow manifold $\mathcal{M}_{\epsilon}$
for the given slow variables $(\gamma(t),d(t),\dot{\gamma}(t),\dot{d}(t))$.
When the fast variables are $\mathcal{O}(10^{-5})$ close to $\mathcal{M}_{\epsilon}$
after the time value $t\geq t_{\epsilon}=15.6\,s$, we take the point
$x_{\epsilon}(t_{\epsilon})$ belonging to the full trajectory and
use the slow coordinates $(\gamma(t_{\epsilon}),d(t_{\epsilon}),\dot{\gamma}(t_{\epsilon}),\dot{d}(t_{\epsilon}))$
of this point as an initial position for the reduced model (\ref{eq:red_mod_1})
and (\ref{eq:red_mod_2}). Consecutively, we simulate the reduced
model in backward and forward time and compare the results with the
results obtained from the full model (see Figs. \ref{fig:minipage_gamma},
\ref{fig:minipage_d}, and \ref{fig:minipage_h}). 

\subsection{Two stiff degrees of freedom\label{subsec:One-slow-degree}}

We reconsider here the same mechanical system as in section \ref{subsec:Two-slow-degrees},
but assume now that both the $d$ and $h$ variables represent stiff
degrees of freedom, while $\gamma$ still describes a soft degree
of freedom. In this setting, the anticipated slow variable $x$ and
fast variable $y=(y_{d},y_{h})$ are defined as

\[
x_{\gamma}=\gamma,\quad y_{d}=d,\quad y_{h}=h,
\]
We express the stiffness of the $y$ degree of freedom by letting
\[
\Delta=\frac{l}{L}=\frac{\delta}{\epsilon},\qquad q_{h}=\frac{\omega_{h}^{2}}{\omega_{p}^{2}}=\frac{\Omega_{h}^{2}}{\epsilon^{2}},\quad q_{d}=\frac{\omega_{d}^{2}}{\omega_{p}^{2}}=\frac{\Omega_{d}^{2}}{\epsilon^{2}},\quad a_{h}=\frac{\alpha_{h}}{\epsilon^{4}},
\]
\begin{equation}
\pi_{h}=\frac{C_{h}}{\omega_{p}M}=\frac{\mu_{h}}{\epsilon},\quad\pi_{d}=\frac{C_{d}}{\omega_{p}M}=\frac{\mu_{d}}{\epsilon},\quad\pi_{p}=\frac{c_{p}}{\omega_{p}mL^{2}}=\frac{c_{p}}{\omega_{p}m\left(\frac{\epsilon}{\delta}l\right)^{2}}=\frac{c_{p}\delta^{2}}{\omega_{p}m\epsilon^{2}l^{2}}=\frac{\mu_{p}}{\epsilon^{2}}.\label{eq:rescaling for fast-fast-slow}
\end{equation}
As we show in Appendix \ref{sec:Appendix:-Details-for-fast-fast-slow system},
assumptions (A1)-(A3) are satisfied in the parameter regime represented
by the above scaling for $0<\epsilon\ll1$. In the scaled variables,
we have
\[
M_{1}=\frac{\delta^{2}}{\epsilon^{2}\left(1+\beta\right)},
\]
thus the mass matrix $M_{1}$ associated with the slow degree of freedom
is not differentiable at $\epsilon=0.$ Therefore, only the more general
form \eqref{eq:reduced model} of the reduced model is applicable,
giving

\begin{eqnarray*}
\ddot{x} & = & P_{1}\left(x,\dot{x},G_{0}(x,\dot{x},t),0,t;0\right)+\mathcal{O}(\epsilon)\\
 & = & -\frac{\mu_{p}}{\delta^{2}}\dot{x}-\sin x+G_{p}(t)+\mathcal{O}(\epsilon).
\end{eqnarray*}
Scaling back to the original time, we conclude that at leading order,
the exact reduced-order model on the two-dimensional, attracting slow
manifold $\mathcal{M}_{\epsilon}$ is given by
\[
\ddot{x}+\frac{\mu_{p}}{\delta^{2}}\dot{x}+\sin x=G_{p}(t)+\mathcal{O}(\epsilon),
\]
 or, equivalently, 
\begin{equation}
ml^{2}\ddot{x}+c_{p}\dot{x}+mgl\sin x=f_{p}(t)l+\mathcal{O}(\epsilon).\label{eq:fast-fast}
\end{equation}

\begin{figure}[H]
\begin{centering}
\subfloat[\label{fig:h_gamma_gammad}]{\begin{centering}
\includegraphics[scale=0.35]{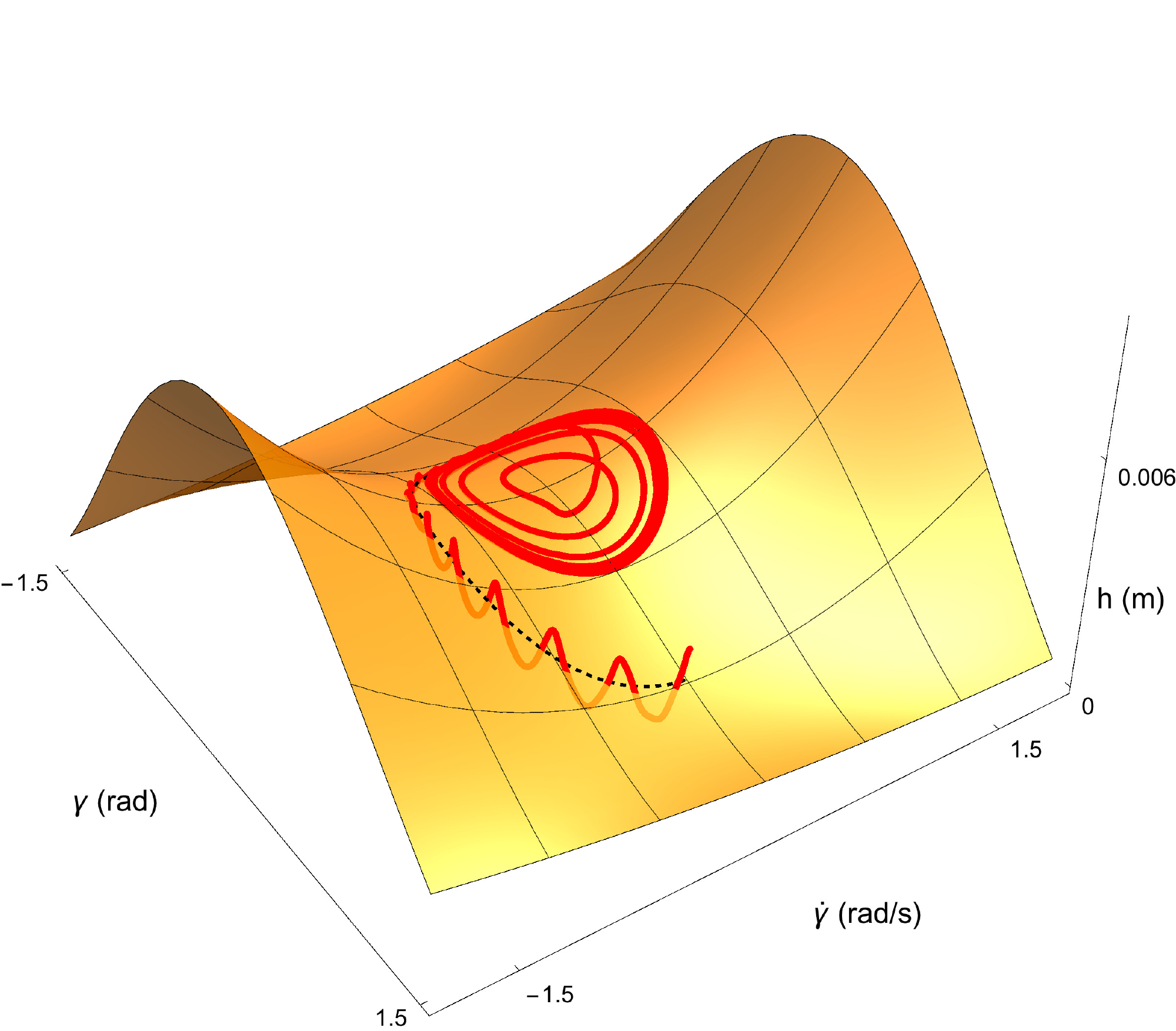}
\par\end{centering}
}\hfill{}\subfloat[\label{fig:d_gamma_gammad}]{\centering{}\includegraphics[scale=0.35]{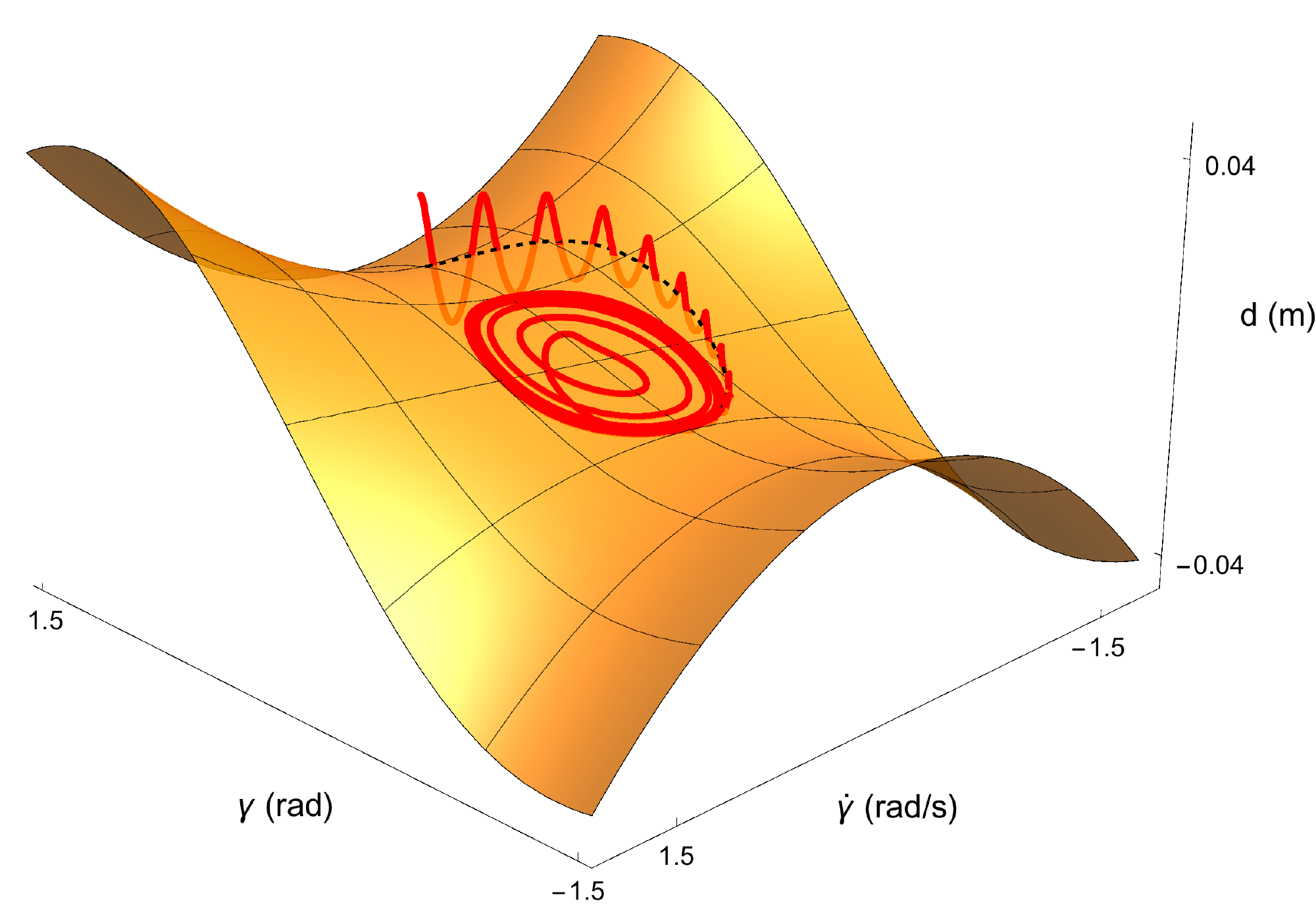}}
\par\end{centering}
\caption{Instantaneous projection of the slow manifold $\mathcal{M_{\epsilon}}$
for the periodically forced, stiff-stiff-soft mechanical system at
$t=20$ s. (\ref{fig:h_gamma_gammad}) A trajectory of the full system
is launched at the initial condition $x_{0}=(1.000,0.000,0.002842,0.02296,0.0005551,-0.002546)$
and integrated in forward time. Displayed in red is the $h$ component
of the corresponding full system trajectory that converges to the
slow manifold and synchronizes with the reduced order model (dashed),
shown up to time $t_{\epsilon}=15.6$ s. (\ref{fig:d_gamma_gammad})
Convergence of the horizontal coordinate $d$ (red) to the slow manifold,
synchronizing with the dynamics of the reduced order model (dashed).
\label{fig:minipage_2DM}}
\end{figure}

As for the example treated in section \ref{subsec:Two-slow-degrees},
we illustrate numerically that trajectories of the full system synchronize
exponentially fast with those of the reduced-order model. For the
parameter values 

\begin{gather*}
l=6\text{ m},\quad L=3\text{ m},\quad M=0.25\text{ kg},\quad m=0.5\text{ kg},\\
K_{h}=2000\text{ N/m},\quad\Gamma_{h}=0.5\text{ N/m}^{3},\quad K_{d}=280\text{ N/m},\\
C_{d}=3\cdot\omega_{p}\cdot M\text{ kg/s},\quad C_{h}=3\cdot\omega_{p}\cdot M\text{ kg/s},\\
c_{p}=\omega_{p}\cdot m\cdot L^{2}\text{ }(\text{kg}\cdot\text{m}^{2})\text{/s},\quad g=9.81\text{ m/s}^{2},\\
f_{p}(t)=0.6\cdot\sin(\omega_{p}t),\quad f_{h}(t)=f_{d}(t)=0,\\
\epsilon=1\cdot10^{-8},
\end{gather*}
and the initial condition, 

\[
x_{0}=(\gamma_{0},\dot{\gamma}_{0},h_{0},d_{0},\dot{h}_{0},\dot{d}_{0})=(1.000,0.000,0.002842,0.02296,0.0005551,-0.002546),
\]
we illustrate the convergence of the trajectory to a trajectory of
the reduced model on the slow manifold in Fig. \ref{fig:minipage_2DM}).
The target model trajectory was identified as earlier in the soft-soft-stiff
version of the same example.

\section{Conclusions}

We have developed a methodology for exact model reduction in multi-degree-of-freedom
mechanical systems with soft and stiff degrees of freedom. This Slow-Fast
Decomposition (SFD) approach allows for a systematic identification
of parameter regimes in which an attracting slow manifold exists.
On this invariant manifolds, the stiff variables are enslaved to the
remaining soft variables. 

We have derived explicit expressions for the slow manifold and for
the first two orders of the reduced flow on this manifold. The latter
formulas provide a mathematically exact reduced order-model with which
trajectories of the full system synchronize at an a priori predictable
exponential rate. We have also identified a domain boundary over which
the slow manifold generically loses its stability and hence the dynamics
on it no longer serves as a reduced-order model for the mechanical
system.

Slow-fast reduction has previously been carried out with varying levels
of mathematical rigor in several specific mechanical model problems
(see the Introduction for a review). Our contributions here are: (i)
explicit conditions under which an attracting slow manifold in guaranteed
to exist in a general, multi-degree-of-freedom mechanical system;
(ii) readily applicable general formulas for reduced-order models
on such manifolds. All these results follow from the application of
classic results from geometric singular perturbation theory (see,
e.g., Fenichel \cite{fenichel79}, Jones \cite{jones06}). 

We have found that the SFD conditions yield reduced-order models that
satisfy the basic requirements (R1)-(R2) we have formulated for a
mathematically exact model reduction procedure. As we has shown explicitly
in Section \ref{sec:Approximate-SFD-near-equilbria}, the formal methods
of static condensation and modal derivatives in structural dynamics
can only be justified if the conditions of SFD are satisfied. When
these conditions do not hold, the reduced-order models produced by
these methods are inaccurate or even qualitatively incorrect.

Importantly, the SFD approach does not require the explicit identification
of eigenfrequencies and normal modes for a linearized system, which
is a numerically costly undertaking for high-degree-of-freedom systems.
Instead, the SFD can be carried out based on a general identification
of stiff and soft vibratory modes, without an explicit decoupling
of these modes. This flexibility for the method enables its application
in structural vibrations problems such as those including forced and
damped beams (cf. Jain et al. \cite{jain16} for a detailed example
involving the von Kármán beam model). 

An extension of the SFD methodology to stiff-soft continuum vibrations
described by partial differential equations should also be possible
through an appropriate extension of the necessary geometric singular
perturbation results to infinite dimensions (see, e.g., Menon and
Haller \cite{menon01}).\\
\\

\textbf{Acknowledgments}

We would like to thank Paolo Tiso, Daniel Rixen and Shobhit Jain for
useful conversations and for their insights on the subject of this
paper.

\section{Appendix : Proof of the main result\label{app: proof of main}}

\subsection{First-order autonomous form}

By the nondegeneracy of $M$, the matrices $M_{11}$ and $M_{22}$
are necessarily invertible, which enables us to split \eqref{eq:systtem00}
in the form
\begin{eqnarray*}
\left[M_{11}-M_{12}M_{22}^{-1}M_{21}\right]\ddot{x} & = & -F_{1}-M_{12}M_{22}^{-1}F_{2},\\
\left[M_{22}-M_{21}M_{11}^{-1}M_{12}\right]\ddot{y} & = & -F_{2}-M_{21}M_{11}^{-1}F_{1}.
\end{eqnarray*}
The nondegeneracy of $M$ also implies that the two matrices on the
left-hand side of this system must be invertible, leading to the explicit
second-order dynamical system
\begin{eqnarray}
\ddot{x} & = & M_{1}^{-1}\left(x,\frac{y}{\epsilon},t;\epsilon\right)Q_{1}\left(x,\dot{x},\frac{y}{\epsilon},\dot{y},t;\epsilon\right),\nonumber \\
\ddot{y} & = & M_{2}^{-1}\left(x,\frac{y}{\epsilon},t;\epsilon\right)Q_{2}\left(x,\dot{x},\frac{y}{\epsilon},\dot{y},t;\epsilon\right),\label{eq:explicit second order system}
\end{eqnarray}
with $M_{i}$ and $Q_{i}$ defined in \eqref{eq:M_i and Q_i def}.

In order to convert this system into a first-order autonomous system,
we first introduce a phase variable $\varphi\in\mathcal{C}$ such
that
\[
\mathcal{C}=\begin{cases}
S^{1}, & M_{i},Q_{i}\,\,\text{\,\,are periodic in }t,\\
\mathbb{T}^{k}, & M_{i},Q_{i}\,\,\text{\,\,are quasi-periodic with\,\,\emph{k}\,\,independent frequencies in }t,\\{}
[a,b], & M_{i},Q_{i}\,\,\text{\,\,are aperiodic in }t.
\end{cases}
\]

We then let 
\[
v=\dot{x},\quad w=\dot{y},
\]
and rewrite equation \eqref{eq:explicit second order system} as a
first-order autonomous system on the extended phase space $\mathcal{P}=\mathbb{R}^{s}\times\mathbb{R}^{s}\times\mathbb{R}^{f}\times\mathbb{R}^{f}\times\mathcal{C}$
in the form
\begin{eqnarray*}
\dot{x} & = & v,\\
\dot{v} & = & M_{1}^{-1}\left(x,\frac{y}{\epsilon},\varphi;\epsilon\right)Q_{1}\left(x,v,\frac{y}{\epsilon},w,\varphi;\epsilon\right),\\
\dot{y} & = & w,\\
\dot{w} & = & M_{2}^{-1}\left(x,\frac{y}{\epsilon},\varphi;\epsilon\right)Q_{2}\left(x,v,\frac{y}{\epsilon},w,\varphi;\epsilon\right),\\
\dot{\varphi} & = & \omega,
\end{eqnarray*}
where
\[
\mathcal{\omega}=\begin{cases}
\omega_{1}, & \mathcal{C}=S^{1},\\
\left(\omega_{1},\ldots,\omega_{k}\right), & \mathcal{C}=\mathbb{T}^{k},\\
1, & \mathcal{C}=[a,b].
\end{cases}
\]

\subsection{Time-scale separation }

Up to this point, the splitting $q=(x,y)$ has been arbitrary. We
now seek conditions under which the $x$-degrees of freedom serve
as coordinates for a reduced-order model. For such a reduced-order
model to capture effectively the long-term system dynamics, we require
the $y$ variables to become enslaved to the $x$ variables and to
the phase variable $\varphi$ over a time scale that is an order of
magnitude faster than the characteristic time scale of the reduced-order
model (cf. the requirement (R2) in the Introduction). To this end,
we introduce a characteristic fast time scale $\tau$ by letting $t=\epsilon\tau,$
with small, non-dimensional parameter $0<\epsilon\ll1.$ Denoting
differentiation with respect to $\tau$ by prime, we obtain the rescaled
equations
\begin{eqnarray}
x^{\prime} & = & \epsilon v,\nonumber \\
v^{\prime} & = & \epsilon M_{1}^{-1}\left(x,\frac{y}{\epsilon},\varphi;\epsilon\right)Q_{1}\left(x,v,\frac{y}{\epsilon},w,\varphi;\epsilon\right),\nonumber \\
\varphi^{\prime} & = & \epsilon\omega,\nonumber \\
y^{\prime} & = & \epsilon w,\nonumber \\
W^{\prime} & = & \epsilon M_{2}^{-1}\left(x,\frac{y}{\epsilon},\varphi;\epsilon\right)Q_{2}\left(x,v,\frac{y}{\epsilon},w,\varphi;\epsilon\right).\label{eq:fully slow eqs}\\
\nonumber 
\end{eqnarray}

In this new scale, the evolution in the $(y,w)$ variables should
be taking place at an $\mathcal{O}(1)$ speed with respect to $\epsilon,$
whereas the $(x,v)$ variables should experience an $\mathcal{O}(\epsilon)$
rate of change. By the structure of system \eqref{eq:fully slow eqs},
this time-scale separation will only arise if we localize $y$ by
letting $y=\epsilon\eta$. With this scaling, we obtain the equations
\begin{eqnarray}
x^{\prime} & = & \epsilon v,\nonumber \\
v^{\prime} & = & \epsilon M_{1}^{-1}\left(x,\eta,\varphi;\epsilon\right)Q_{1}\left(x,v,\eta,w,\varphi;\epsilon\right),\nonumber \\
\varphi^{\prime} & = & \epsilon\omega,\nonumber \\
\eta^{\prime} & = & w,\nonumber \\
w^{\prime} & = & \epsilon M_{2}^{-1}\left(x,\eta,\varphi;\epsilon\right)Q_{2}\left(x,v,\eta,w,\varphi;\epsilon\right).\label{eq:slow-fast-1-1}
\end{eqnarray}
 To ensure that $w$ also varies at $\mathcal{O}(1)$ speeds for small
enough $\epsilon$, the function $\epsilon M_{2}^{-1}\left(x,\eta,\varphi;\epsilon\right)Q_{2}\left(x,v,\eta,w,\varphi;\epsilon\right)$
must have a smooth, $\mathcal{O}(1)$ limit as $\epsilon\to0.$ We,
therefore, must require the function
\[
P_{2}\left(x,v,\eta,w,\varphi;\epsilon\right)=\epsilon M_{2}^{-1}\left(x,\eta,\varphi;\epsilon\right)Q_{2}\left(x,v,\eta,w,\varphi;\epsilon\right)
\]
 to have a smooth limit at $\epsilon=0$, defined by a smooth function
\begin{equation}
P_{2}(x,v,\eta,w,\varphi;0):=\lim_{\epsilon\to0}P_{2}\left(x,v,\eta,w,\varphi;\epsilon\right)\label{eq:limitcond}
\end{equation}
 on an open and bounded subset of the extended phase space $\mathcal{P}$.
In order to be able to carry out a perturbation argument from this
limit, we also require that 
\[
P_{1}\left(x,v,\eta,w,\varphi;\epsilon\right)=M_{1}^{-1}\left(x,\eta,\varphi;\epsilon\right)Q_{1}\left(x,v,\eta,w,\varphi;\epsilon\right)
\]
 has a similar smooth limit at $\epsilon=0$, defined as
\[
P_{1}\left(x,v,\eta,w,\varphi;0\right):=\lim_{\epsilon\to0}P_{1}\left(x,v,\eta,w,\varphi;\epsilon\right).
\]
 With these quantities and assumptions, \eqref{eq:slow-fast-1-1}
becomes
\begin{eqnarray}
x^{\prime} & = & \epsilon v,\nonumber \\
v^{\prime} & = & \epsilon P_{1}\left(x,v,\eta,w,\varphi;\epsilon\right),\nonumber \\
\varphi^{\prime} & = & \epsilon\omega,\nonumber \\
\eta^{\prime} & = & w,\nonumber \\
w^{\prime} & = & P_{2}\left(x,v,\eta,w,\varphi;\epsilon\right).\label{eq:slow-fast-1-1-2}
\end{eqnarray}

\subsection{Existence of a critical manifold}

We want to ensure the existence of a reduced-order model in which
the $(\eta(t),w(t))$ dynamics can be uniquely expressed, at least
for large enough times, as a function of the $(x(t),v(t))$ dynamics
and the time $t$. In geometric terms, this amounts to the existence
of an invariant manifold $\mathcal{M}_{\epsilon}$ that is a graph
over the $(x,v,t)$ variables and attracts all nearby solutions of
the full system. 

We require our reduced model to be smooth in $\epsilon$, which is
equivalent to requiring a smooth limit $\mathcal{M}_{0}=\lim_{\epsilon\to0}\mathcal{M}_{\epsilon}$
for the invariant manifold in the $\epsilon=0$ limit of system \eqref{eq:slow-fast-1-1-2}.
This limiting system can be written as
\begin{eqnarray}
x^{\prime} & = & 0,\nonumber \\
v^{\prime} & = & 0,\nonumber \\
\varphi^{\prime} & = & 0,\nonumber \\
\eta^{\prime} & = & w,\nonumber \\
w^{\prime} & = & P_{2}\left(x,v,\eta,w,\varphi;0\right).\label{eq:slow-fast-1-1-1-1}
\end{eqnarray}

In this limit, therefore, $(x,v,\varphi)\equiv(x_{0},v_{0},\varphi_{0})$
plays the role of a constant parameter vector. Any trajectory of the
fast dynamics
\begin{eqnarray}
\eta^{\prime} & = & w,\nonumber \\
w^{\prime} & = & P_{2}\left(x_{0},v_{0},\eta,w,\varphi_{0};0\right),\label{eq:fast dynamics}
\end{eqnarray}
therefore, gives rise to a $(2s+1)$ -dimensional invariant manifold
for the full system. Along nontrivial trajectories of \eqref{eq:fast dynamics},
however, the $(\eta,v)$ variables change and hence are \emph{not}
uniquely enslaved to $(x_{0},v_{0},\varphi_{0})$, as required for
the smooth limit of a reduced-order model. Consequently, only invariant
manifolds arising from fixed points of \eqref{eq:fast dynamics} can
be considered as limits of reduced-order models. 

Such fixed points of \eqref{eq:fast dynamics} form a set 
\[
\mathcal{M}_{0}=\left\{ (x,v,\eta,w,\varphi)\in\mathcal{P}\,:\,w=0,\quad P_{2}\left(x,v,\eta,w,\varphi;0\right)=0\right\} .
\]
To be a limit of a slow manifold carrying a reduced-order model, $\mathcal{M}_{0}$
must be a smooth graph over an open domain $\mathcal{D}_{0}\subset\mathbb{R}^{m}\times\mathbb{R}^{m}\times\mathcal{C}$
of the space $(x,v,t)$ variables. By the implicit function theorem,
this is equivalent to the requirement that
\begin{equation}
\det\left[\partial_{\eta}P_{2}\left(x,v,\eta,0,\varphi;0\right)\right]\neq0,\label{eq:detnonzero}
\end{equation}
should hold at all points $(x,v,\eta,w,\varphi)\in\mathcal{M}_{0}$.
This condition ensures that if $\mathcal{M}_{0}$ is nonempty, then
it is a $2s+1$ dimensional differentiable manifold that can locally
be expressed as a smooth graph 
\begin{equation}
\left(\begin{array}{c}
\eta\\
w
\end{array}\right)=\left(\begin{array}{c}
G_{0}(x,v,\varphi)\\
0
\end{array}\right),\quad(x,v,\varphi)\in\mathcal{D}_{0}\label{eq:G_0_def-1}
\end{equation}
with the function $G_{0}:\mathcal{D}_{0}\to\mathbb{R}^{f}$ satisfying
the identity 
\begin{equation}
P_{2}\left(x,v,G_{0}(x,v,\varphi),0,\varphi;0\right)=0.\label{eq:identity-1}
\end{equation}
We refer to the part of $\mathcal{M}_{0}$ satisfying \eqref{eq:detnonzero}
as the critical manifold associated with the limiting system \eqref{eq:slow-fast-1-1-1-1}.
In our discussion of assumption (A2), we use the term \emph{critical
manifold }for the $t=const.$ times slice $\mathcal{M}_{0}(t)$ of
$\mathcal{M}_{0}$. 

\subsection{Stability of $\mathcal{M}_{0}$ }

The critical manifold must be normally attracting to persist as an
attracting invariant slow manifold in the full system \eqref{eq:slow-fast-1-1-2}.
The stability type of $\mathcal{M}_{0}$ can be identified by analyzing
the linearization of the fast flow \eqref{eq:fast dynamics} at the
fixed points forming $\mathcal{M}_{0}$. 

The stability of the manifold $\mathcal{M}_{0}$ at the fixed point
family $(\eta_{0},w_{0})=(G(x_{0},v_{0},\varphi_{0}),0)$ of the decoupled
equations is governed by the eigenvalues of the Jacobian 
\begin{equation}
J=\left[\begin{array}{cc}
0 & I\\
\partial_{\eta}P_{2} & \partial_{w}P_{2}
\end{array}\right]_{\left(x,v,\eta,w,t;\epsilon\right)=\left(x_{0},v_{0},G_{0}(y_{0},w_{0},\varphi_{0}),0,\varphi_{0};0\right)}.\label{eq:Jacobian}
\end{equation}

The matrix $J$ has eigenvalues with strictly negative real parts
precisely when the fixed point of the linear vibratory system
\begin{equation}
u^{\prime\prime}-\partial_{w}P_{2}\left(x_{0},v_{0},G_{0}(x_{0},v_{0},t_{0}),0,\varphi_{0};0\right)u^{\prime}-\partial_{\eta}P_{2}\left(x_{0},v_{0},G_{0}(x_{0},v_{0},\varphi_{0}),0,\varphi_{0};0\right)u=0\label{eq:associated linear problem}
\end{equation}
is asymptotically stable for the parameter values $\left(x_{0},v_{0},\varphi_{0}\right)\in\mathcal{D}_{0}$,
which is guaranteed by assumption (A3). In that case, a compact subset
of the critical manifold $\mathcal{M}_{0}$ is a compact normally
hyperbolic invariant manifold with boundary when $(x_{0},v_{0},\varphi_{0})$
is restricted to a domain with a smooth boundary. (In case of $\mathcal{C}=[a,b]$,
one has to select $a$ and\textbf{ $b$ }as\textbf{ }smooth functions
of $(y_{0},t_{0})$ to eliminate non-smooth corners in $\partial\mathcal{M}_{0}$.
This can always be done without loss of generality.)

\subsection{Existence of a slow manifold}

Under the above conditions, the results of Fenichel \cite{fenichel79}
guarantee for the full system \eqref{eq:slow-fast-1-1-2} the existence
of an attracting slow manifold $\mathcal{M}_{\epsilon}$ that is $\mathcal{O}(\epsilon)$
$C^{r}$ -close to $\mathcal{M}_{0}$, and hence continues to be a
graph of the form
\[
\left(\begin{array}{c}
\eta\\
w
\end{array}\right)=\left(\begin{array}{c}
G_{\epsilon}(x,v,\varphi)\\
\epsilon H_{\epsilon}(x,v,\varphi)
\end{array}\right)=\left(\begin{array}{c}
G_{0}(x,v,\varphi)+\epsilon G_{1}(x,w,\varphi)+\mathcal{O}(\epsilon^{2})\\
\epsilon H_{0}(x,v,\varphi)+\epsilon^{2}H_{1}(x,v,\varphi)+\mathcal{O}(\epsilon^{3})
\end{array}\right),\quad(x,v,\varphi)\in\mathcal{D}_{0},
\]
with appropriate smooth functions $G_{\epsilon}$ and $H_{\epsilon}$.
The relation $\eta^{\prime}=w$ in \eqref{eq:slow-fast-1-1-2} imposes
the relationships 
\begin{eqnarray*}
\frac{d}{d\tau}\left[G_{0}(x,v,\varphi)+\epsilon G_{1}(x,v,\varphi)+\mathcal{O}(\epsilon^{2})\right] & = & \epsilon H_{0}(x,v,\varphi)+\epsilon^{2}H_{1}(x,v,\varphi)+\mathcal{O}(\epsilon^{3}),\\
\frac{d}{d\tau}\left[\epsilon H_{0}(x,v,\varphi)+\epsilon^{2}H_{1}(x,v,\varphi)+\mathcal{O}(\epsilon^{3})\right] & = & P_{2}\left(x,v,G_{\epsilon}(x,v,\varphi),\epsilon H_{\epsilon}(x,v,\varphi),\varphi;\epsilon\right).
\end{eqnarray*}
Carrying out the differentiation in these two equations gives
\[
\left(\epsilon\partial_{x}G_{0}+\epsilon^{2}\partial_{x}G_{1}\right)v+\left(\epsilon\partial_{v}G_{0}+\epsilon^{2}\partial_{v}G_{1}\right)P_{1}\left(x,v,G_{\epsilon},\epsilon H_{\epsilon},\varphi;\epsilon\right)+\left(\epsilon\omega\partial_{\varphi}G_{0}+\epsilon^{2}\omega\partial_{\varphi}G_{1}\right)+\mathcal{O}(\epsilon^{3})
\]
 
\[
=\epsilon H_{0}+\epsilon^{2}H_{1}+\mathcal{O}(\epsilon^{3}),
\]
\[
\left(\epsilon^{2}\partial_{x}H_{0}+\epsilon^{3}\partial_{x}H_{1}\right)v+\left(\epsilon^{2}\partial_{v}H_{0}+\epsilon^{3}\partial_{v}H_{1}\right)P_{1}\left(x,v,G_{\epsilon},\epsilon H_{\epsilon},\varphi;\epsilon\right)+\left(\epsilon^{2}\omega\partial_{\varphi}H_{0}+\epsilon^{3}\omega\partial_{\varphi}H_{1}\right)+\mathcal{O}(\epsilon^{4})
\]
 
\[
=P_{2}\left(x,v,G_{\epsilon}(x,v,\varphi),\epsilon H_{\epsilon}(x,v,\varphi),\varphi;\epsilon\right).
\]
We Taylor-expand these two equations, then equate the $\mathcal{O}(\epsilon)$
and $\mathcal{O}(\epsilon^{2})$ terms in the first equation, as well
as $\mathcal{O}(\epsilon)$ terms in the second equation, to obtain
\begin{eqnarray*}
H_{0}(x,v,\varphi) & = & \partial_{x}G_{0}(x,v,\varphi)v+\partial_{v}G_{0}(x,v,\varphi)P_{1}\left(x,v,G_{0}(x,v,\varphi),0,\varphi;0\right)+\omega\partial_{\varphi}G_{0}(x,v,\varphi),\\
H_{1}(x,w,\varphi) & = & \partial_{x}G_{1}(x,v,\varphi)v+\partial_{v}G_{1}(x,v,\varphi)P_{1}\left(x,v,G_{0}(x,v,\varphi),0,\varphi;0\right)+\omega\partial_{\varphi}G_{1}(x,v,\varphi),\\
G_{1}(x,v,\varphi) & = & -\left[D_{\eta}P_{2}\left(x,v,G_{0}(x,v,\varphi),0,\varphi;0\right)\right]^{-1}D_{w}P_{2}\left(x,v,G_{0}(x,v,\varphi),0,\varphi;0\right)H_{0}(x,v,\varphi)\\
 &  & -\left[D_{\eta}P_{2}\left(x,v,G_{0}(x,v,\varphi),0,\varphi;0\right)\right]^{-1}D_{\epsilon}P_{2}\left(x,v,G_{0}(x,v,\varphi),0,\varphi;0\right).
\end{eqnarray*}

In terms of the original variables, therefore, the slow manifold satisfies
\begin{eqnarray*}
y & = & \epsilon G_{0}(x,\dot{x},t)+\epsilon^{2}G_{1}(x,\dot{x},t)+\mathcal{O}(\epsilon^{3}),\\
\dot{y} & = & \epsilon H_{0}(x,\dot{x},t)+\epsilon^{2}H_{1}(x,\dot{x},t)+\mathcal{O}(\epsilon^{3}),
\end{eqnarray*}
where the functions $H_{0},G_{1}$ and $H_{1}$ are those listed in
\eqref{eq:H0G1H1}.

\subsection{The reduced flow on the slow manifold}

The slow manifold $\mathcal{M}_{\epsilon}$ attracts all nearby solutions,
thus the reduced flow on $\mathcal{M}_{\epsilon}$ will serve as the
type of reduced-order model we have been seeking to construct (cf.
requirement (R1) in the Introduction). The reduced equations on $\mathcal{M}_{\epsilon}$
can be written by restricting the $(x,v,\phi)$ components of our
system \eqref{eq:slow-fast-1-1-2} to $\mathcal{M}_{\epsilon}$, which
yields 
\begin{eqnarray*}
x^{\prime} & = & \epsilon v,\\
v^{\prime} & = & \epsilon P_{1}\left(x,v,G_{0}(x,v,\varphi),0,\varphi;0\right)\\
 &  & +\epsilon^{2}\left[D_{\eta}P_{1}\left(x,v,G_{0}(x,v,\varphi),0,\varphi;0\right)G_{1}(x,v,\varphi)+D_{w}P_{1}\left(x,v,G_{0}(x,v,\varphi),0,\varphi;0\right)H_{0}(x,v,\varphi)\right.\\
 &  & \left.\,\,\,\,\,\,\,\,\,+D_{\epsilon}P_{2}\left(x,v,G_{0}(x,v,\varphi),0,\varphi;0\right)\right]+\mathcal{O}(\epsilon^{2}),\\
\varphi^{\prime} & = & \epsilon\omega.
\end{eqnarray*}
In the original set of coordinates, this reduced flow can be written
as in eq. \eqref{eq:reduced model}. 

Using the definition of $P_{1}$, we find that if $M_{1}(x,G_{0}(x,\dot{x},t),t)$
has a smooth limit at $\epsilon=0$, then the reduced equation can
be multiplied by $M_{1}(x,G_{0}(x,\dot{x},t),t)$ to yield the leading-order
equivalent form of \eqref{eq:reduced model} as given in eq. \eqref{eq:reduced model 1}.
When necessary, the $\mathcal{O}(\epsilon)$ terms in \eqref{eq:reduced model 1}
can also be computed from the formulas we have given above.

\subsection{Convergence to the reduced trajectories}

By the invariant foliation results of Fenichel \cite{fenichel79},
for small enough $\epsilon$ and for motions close enough to the critical
manifold, the $y(t)$ component of all solutions of equation \eqref{eq:systtem00}
synchronize exponentially fast with solutions of the reduced-order
model \eqref{eq:reduced model}. 

Specifically, the local stable manifold $W_{loc}^{s}(\mathcal{M}_{\epsilon})$
is foliated by an invariant family of class $C^{r}$ \emph{stable
fibers} $f^{s}(p)$. This $\left(2s+\dim\mathcal{C}\right)$-parameter
fiber-family is parametrized by the \emph{base points} $p\in\mathcal{M}_{\epsilon}$
of the fibers. Each fiber is a class $C^{r-1}$ manifold whose dimension
is $2f$. The invariance of the fiber family means that for the flow
map $F^{\tau}\colon\mathcal{P}\to\mathcal{P}$ of system \eqref{eq:slow-fast-1-1-2},
we have 
\[
F^{\tau}\left(f^{s}(p)\right)\subset f^{s}(F^{\tau}(p))
\]
for all $\tau>0.$ Furthermore, the trajectory of the reduced flow
through a fiber base point $p$ attracts exponentially all trajectories
that cross the fiber $f^{s}(p)$. Specifically, if $p=\left(x_{R}(\tau_{0}),v_{R}(\tau_{0}),\varphi_{R}(\tau_{0})\right)$
and $\left(x(\tau_{0}),v(\tau_{0}),\varphi(\tau_{0}),\eta(\tau_{0}),w(\tau_{0})\right)\in f^{s}(p)$,
then for all $\tau$ values satisfying 
\[
\left(x(\tau),v(\tau),\varphi(\tau),\eta(\tau),w(\tau)\right)\in W_{loc}^{s}(\mathcal{M}_{\epsilon}),
\]
 we have the estimate 
\begin{equation}
\left|\left(\begin{array}{c}
x(\tau)-x_{R}(\tau)\\
v(\tau)-v_{R}(\tau)\\
\varphi(\tau)-\varphi_{R}(\tau)\\
\eta(\tau)-\eta_{R}(\tau)\\
w(\tau)-w_{R}(\tau)
\end{array}\right)\right|\leq C\left|\left(\begin{array}{c}
x(\tau_{0})-x_{R}(\tau_{0})\\
v(\tau_{0})-v_{R}(\tau_{0})\\
\varphi(\tau_{0})-\varphi_{R}(\tau_{0})\\
\eta(\tau_{0})-\eta_{R}(\tau_{0})\\
w(\tau_{0})-w_{R}(\tau_{0})
\end{array}\right)\right|e^{-\Lambda(\tau-\tau_{0})},\quad\tau>\tau_{0}.\label{eq:synchronization inequaltiy-2}
\end{equation}
Here $\Lambda>0$ can be selected as any constant satisfying 
\[
\max_{j\in[1,2f],\,(x,v,\varphi)\in\mathcal{D}_{0}}\mathrm{Re}\,\lambda_{j}(x,v,\varphi)<-\Lambda<0,
\]
with $\lambda_{j}(x,v,\varphi),$ $j=1,\ldots,2f$, denoting the eigenvalues
of the Jacobian $J$, or equivalently, of the associated linear system
\eqref{eq:associated linear problem}. The constant $C>0$ depends
on $\Lambda$ but is independent of the choice of the fiber base point
$p$ and the times $\tau$ and $\tau_{0}$. 

By the form of system system \eqref{eq:slow-fast-1-1-2}, we have
$\left|\varphi(\tau)-\varphi_{R}(\tau)\right|=\left|\varphi(\tau_{0})-\varphi_{R}(\tau_{0})\right|$.
This is only consistent with \eqref{eq:synchronization inequaltiy-2},
if $\varphi(\tau_{0})\equiv\varphi_{R}(\tau_{0})$, which implies
that the fibers $f^{s}(p)$ are necessarily flat (i.e, constant) in
the coordinate $\varphi.$ Using this fact in \eqref{eq:synchronization inequaltiy-2}
and passing back to the original coordinates gives
\[
\left|\left(\begin{array}{c}
x(t)-x_{R}(t)\\
\dot{x}(t)-\dot{x}_{R}(t)\\
\frac{1}{\epsilon}y(t)-\frac{1}{\epsilon}y_{R}(t)\\
\dot{y}(t)-\dot{y}(t)
\end{array}\right)\right|\leq C\left|\left(\begin{array}{c}
x(t_{0})-x_{R}(t_{0})\\
\dot{x}(t_{0})-\dot{x}_{R}(t_{0})\\
\frac{1}{\epsilon}y(t_{0})-\frac{1}{\epsilon}y_{R}(t_{0})\\
\dot{y}(t_{0})-\dot{y}_{R}(t_{0})
\end{array}\right)\right|e^{-\frac{\Lambda}{\epsilon}(t-t_{0})},\quad\tau>\tau_{0}.
\]

Along the reduced flow on the slow manifold $\mathcal{M}_{\epsilon}$,
the $(y,\dot{y})$ variables are enslaved to the $(x,v,t)$ variables,
thus we can further rewrite this last inequality as
\[
\left|\left(\begin{array}{c}
x(t)-x_{R}(t)\\
\dot{x}(t)-\dot{x}_{R}(t)\\
\frac{1}{\epsilon}y(t)-G_{\epsilon}\left(x_{R}(t),\dot{x}_{R}(t),t\right)\\
\dot{y}(t)-\epsilon H_{\epsilon}\left(x_{R}(t),\dot{x}_{R}(t),t\right)
\end{array}\right)\right|\leq C\left|\left(\begin{array}{c}
x(t_{0})-x_{R}(t_{0})\\
\dot{x}(t_{0})-\dot{x}_{R}(t_{0})\\
\frac{1}{\epsilon}y(t_{0})-G_{\epsilon}\left(x_{R}(t_{0}),\dot{x}_{R}(t_{0}),t\right)\\
\dot{y}(t_{0})-\epsilon H_{\epsilon}\left(x_{R}(t_{0}),\dot{x}_{R}(t_{0}),t\right)
\end{array}\right)\right|e^{-\frac{\Lambda}{\epsilon}(t-t_{0})},\quad\tau>\tau_{0}.
\]
Applying the triangle inequality to the left-hand-side and using the
definition of $G_{\epsilon}$ and $H_{\epsilon}$ on the right-hand
side of this inequality proves formula \eqref{eq:synchronization inequaltiy}.

\section{Appendix: Proof of Proposition \ref{prop:Guyan and modal derivatives}
\label{app:Guyan and modal derivatives}}

We start by noting that, as a consequence of assumption \eqref{eq:P_2-mod-der},
the graph
\[
\eta=G_{0}(x,t)
\]
of the critical manifold $\mathcal{M}_{0}$ depends only on the slow
positions $x$ and the time $t$. Near the unperturbed equilibrium,
$\mathcal{M}_{0}(t)$ can therefore be approximated by its Taylor
expansion with respect to $x$. Specifically, we have
\begin{eqnarray}
\eta=G_{0}(x,t) & = & G_{0}(0,t)+\partial_{x}G_{0}(0,t)x+\frac{1}{2}\left(\partial_{xx}^{2}G_{0}(0,t)x\right)x+\mathcal{O}\left(\left|x\right|^{3}\right).\label{eq:Taylor1}
\end{eqnarray}

Differentiation of the implicit equation $P_{2}(x,G_{0}(x,t),0,t;0)=0$
with respect to $x$ gives

\begin{equation}
\partial_{x}P_{2}+\partial_{\eta}P_{2}\partial_{x}G_{0}=0.\label{eq:Taylor2}
\end{equation}
Substitution of \eqref{eq:Taylor1} into\eqref{eq:Taylor2} and setting
$x=0$ yields 
\[
\partial_{x}G_{0}(0,t)=-\left.\left[\partial_{\eta}P_{2}\right]^{-1}\partial_{x}P_{2}\right|_{x=0,\eta=G(0,t),\dot{y}=0,\epsilon=0},
\]
where the inverse of $\partial_{\eta}P_{2}(x,G_{0}(x,t),0,t;0)$ is
guaranteed to exist by assumption (A3). Differentiating \eqref{eq:Taylor2}
once more in $x$ gives
\[
\partial_{xx}^{2}P_{2}+\left(2\partial_{x\eta}^{2}P_{2}+\partial_{\eta\eta}^{2}P_{2}\partial_{x}G_{0}\right)\partial_{x}G_{0}+\partial_{\eta}P_{2}\partial_{x}^{2}G_{0}=0,
\]
enabling us to express the three-tensor $\partial_{xx}^{2}G_{0}(0,t)$
as
\begin{equation}
\partial_{xx}^{2}G_{0}(0,t)=-\left.\left[\partial_{\eta}P_{2}\right]^{-1}\left[\partial_{xx}^{2}P_{2}+\left(2\partial_{x\eta}^{2}P_{2}+\partial_{\eta\eta}^{2}P_{2}\partial_{x}G_{0}\right)\partial_{x}G_{0}\right]\right|_{x=0,\eta=G(0,t),\dot{y}=0,\epsilon=0}.\label{eq:Taylor3}
\end{equation}

Therefore, with the help of the formulas \eqref{eq:P_2-mod-der},
the critical manifold $\mathcal{M}_{0}$ can be written near the origin
as a smooth, codimension-$2f$ graph of the form
\begin{equation}
\mathcal{M}_{0}(t)=\left\{ (x,\dot{x},\eta,\dot{y},t)\in\mathcal{P}\,:\,\eta=G_{0}(x,t)=\Gamma(t)+\Phi(t)x+\left(\Theta(t)x\right)x+\mathcal{O}\left(\left|x\right|^{3}\right),\quad\dot{y}=0\right\} ,\label{eq:crit_man_local}
\end{equation}
where
\begin{eqnarray}
P_{2}(0,\Gamma(t),0,t;0) & = & 0,\nonumber \\
\Phi(t) & = & -\left.\left[\partial_{\eta}P_{2}\right]^{-1}\partial_{x}P_{2}\right|_{x=0,\eta=\Gamma(t),\dot{y}=0,\epsilon=0},\nonumber \\
\Theta(t) & = & -\left.\frac{1}{2}\left[\partial_{\eta}P_{2}\right]^{-1}\left[\partial_{xx}^{2}P_{2}+\left(2\partial_{x\eta}^{2}P_{2}+\partial_{\eta\eta}^{2}P_{2}\Phi(t)\right)\Phi(t)\right]\right|_{x=0,\eta=\Gamma(t),\dot{y}=0,\epsilon=0},\label{eq:MD-general_form}
\end{eqnarray}
as claimed in statement (i) of the Proposition. These expressions
in \eqref{eq:MD-general_form} can then be used in the reduced\textendash order
models \eqref{eq:reduced model}-\eqref{eq:reduced model 1} to obtain
more specific local approximations to the reduced dynamics, in case
a global expression for the critical manifold is not explicitly available.
Specifically, \eqref{eq:reduced model} can be localized near $x=0$
as 
\begin{eqnarray}
\ddot{x}-P_{1}\left(x,\dot{x},\left[\Gamma(t)+\Phi(t)x+\left(\Theta(t)x\right)x\right],0,t;0\right)+\mathcal{O}(\epsilon,\left|x\right|^{3}) & = & 0.\label{eq:reduced model-1}
\end{eqnarray}

Under the further assumptions in statement (ii) of the Proposition,
we have the following simplifications in formulas \eqref{eq:MD-general_form}:
\begin{equation}
\Gamma(t)\equiv0,\qquad\Phi(t)\equiv0,\qquad\Theta(t)\equiv-\frac{1}{2}\left[\partial_{\eta}P_{2}(0,0,0;0)\right]^{-1}\partial_{xx}^{2}P_{2}(0,0,0;0).\label{eq:simple_MD_formulas}
\end{equation}
Substituting these quantities into \eqref{eq:reduced model-1} and
truncating the expression for $\mathcal{M}_{0}(t)$ at linear and
then at quadratic order proves the leading-order forms of the reduced
equations in statements (iii) and (iv) of the Proposition, respectively.
To obtain the order of the error terms in these equations, note that
if $x$ and $y$ are modal coordinates of the linearized system, then
we have 
\begin{eqnarray}
P_{1}(x,\dot{x},\eta,\dot{y},t;\epsilon) & = & P_{1}(x,\dot{x},\eta,0,t;0)+\mathcal{\mathcal{O}}(\epsilon)\nonumber \\
 & = & P_{1}(x,\dot{x},0,0,t;0)+\mathcal{\mathcal{O}}(\left|x\right|\left|\eta\right|)+\mathcal{\mathcal{O}}(\epsilon).\label{eq:P1expansion}
\end{eqnarray}
Substitution of $\eta=0+\mathcal{O}(\left|x\right|^{2})$ and $\eta=\left(\Theta(t)x\right)x+\mathcal{O}(\left|x\right|^{3})$,
respectively, into the $\mathcal{\mathcal{O}}(\left|x\right|\left|\eta\right|)$
term in \eqref{eq:P1expansion} then proves the order of the higher-order
terms, as listed in statements (iii) and (iv) of the Proposition.

\section{Appendix: Details for Example \ref{ex:reduction failure}\label{sec:Appendix:-Details-for-failure}}

For the system 
\begin{eqnarray}
\ddot{x}+\left(c_{1}+\mu_{1}x^{2}\right)\dot{x}+k_{1}x+axy+bx^{3} & = & 0,\qquad x\in\mathbb{R},\nonumber \\
\ddot{y}+c_{2}\dot{y}+k_{2}y+cx^{2} & = & 0,\qquad y\in\mathbb{R},\label{eq:2DOF example-1}
\end{eqnarray}
we consider reduction by static condensation via the linear change
of variables
\begin{equation}
\left(\begin{array}{c}
x\\
y
\end{array}\right)=U\hat{x},\qquad U=\left(\begin{array}{c}
1\\
0
\end{array}\right),\quad\hat{x}\in\mathbb{R}.\label{eq:Guyan example-1}
\end{equation}
Dropping the tilde from $\hat{x}$ and substituting $y=0$ from \eqref{eq:Guyan example-1}
into the first equation of \eqref{eq:2DOF example} gives the statically
condensed model \eqref{eq:ex guyan reduced}. 

Next, applying the idea of modal derivatives, we seek a quadratic
invariant manifold of the form \eqref{eq:G0 quadratic}, with the
coefficients computed in the unscaled variables as\emph{ 
\begin{eqnarray}
\Phi & = & 0,\nonumber \\
\Theta & = & -\left.\frac{1}{2}\left[\partial_{y}\left(c_{2}\dot{y}+k_{2}y+cx^{2}\right)\right]^{-1}\left[\partial_{xx}^{2}\left(c_{2}\dot{y}+k_{2}y+cx^{2}\right)\right]\right|_{x=0,y=0,\dot{y}=0}=-\frac{c}{k_{2}}.\label{eq:MD-general_form-1-1-1}
\end{eqnarray}
}Substitution of $y=\Theta x^{2}$ into into the first equation of
\eqref{eq:2DOF example} gives the modal-derivate-based reduced-order
model \eqref{eq:ex md reduced}, representing only a slight correction
to \eqref{eq:ex guyan reduced} at cubic order. All this appears reasonable
at this point, with the statically condensed system \eqref{eq:ex guyan reduced}
offering a leading-order model that is subsequently refined at cubic
order by the modal derivatives approach in \eqref{eq:ex md reduced}. 

At the same time, there exists a slow spectral submanifold (SSM),
the unique smoothest, nonlinear continuation of the $y=0$ modal subspace
of the equilibrium. This unique, two-dimensional analytic invariant
manifold is tangent to the modal subspace of the $x$-degree of freedom
at the origin (cf. Haller and Ponsioen \cite{haller16}). The slow
SSM, therefore, can locally be written as a two-dimensional invariant
graph $(y,\dot{y})=\left(g_{1}(x,\dot{x}),g_{2}(x,\dot{x})\right)=\mathcal{O}\left(x^{2},x\dot{x},\dot{x}^{2}\right)$
over $(x,\dot{x})$, as originally envisioned by Shaw and Pierre \cite{shaw93}.
Differentiating the general form 
\begin{equation}
y=g_{1}(x,\dot{x})=\alpha x^{2}+\beta x\dot{x}+\gamma\dot{x}^{2}+\mathcal{O}\left(3\right)\label{eq:g1-1}
\end{equation}
of such an invariant graph twice in time, with $\ddot{x}$ substituted
from the first equation of system \eqref{eq:2DOF example}, we obtain
\begin{eqnarray*}
\ddot{y} & = & -k_{1}(2\alpha-2\gamma k_{1}-\beta c_{1})x^{2}\\
 &  & -\left[2\beta k_{1}+c_{1}\left(2\alpha-2\gamma k_{1}-\beta c_{1}\right)+2\left(\beta-2\gamma c_{1}\right)k_{1}\right]x\dot{x}\\
 &  & +\left[\left(2\alpha-2\gamma k_{1}-\beta c_{1}\right)-2c_{1}\left(\beta-2\gamma c_{1}\right)\right]\dot{x}^{2}\\
 &  & +\mathcal{O}\left(3\right).
\end{eqnarray*}
A comparison of this differential equation with the second equation
of system \eqref{eq:2DOF example}, with $y$ and $\dot{y}$ substituted
from \eqref{eq:g1-1}, leads to the linear system of algebraic equations
\[
\left(\begin{array}{ccc}
k_{2}-2k_{1} & k_{1}\left(c_{1}-c_{2}\right) & 2k_{1}^{2}\\
2\left(c_{2}-c_{1}\right) & k_{2}-4k_{1}+c_{1}^{2}-c_{1}c_{2} & 2k_{1}\left(3c_{1}-c_{2}\right)\\
2 & c_{2}-3c_{1} & k_{2}-2k_{1}+4c_{1}^{2}-2c_{1}c_{2}
\end{array}\right)\left(\begin{array}{c}
\alpha\\
\beta\\
\gamma
\end{array}\right)=-\left(\begin{array}{c}
c\\
0\\
0
\end{array}\right)
\]
for the unknown coefficients $\alpha$, $\beta$ and $\gamma$ in
the expression \eqref{eq:g1-1} of the slow SSM. The solution of this
system of equations is given by
\begin{eqnarray}
\alpha & = & -\frac{c}{D}\left(4c_{1}^{4}-6c{}_{1}^{3}c_{2}+2c{}_{1}^{2}c{}_{2}^{2}+5c{}_{1}^{2}k_{2}-c_{1}c_{2}\left(2k_{1}+3k_{2}\right)+2c{}_{2}^{2}k_{1}+8k_{1}^{2}-6k_{1}k_{2}+k{}_{2}^{2}\right),\nonumber \\
\beta & = & -\frac{2c}{D}\left(4c_{1}k_{1}+k_{2}\left(c_{1}-c_{2}\right)+2c_{1}c_{2}^{2}-6c_{1}^{2}c_{2}+4c_{1}^{3}\right),\nonumber \\
\gamma & = & -\frac{2c}{D}\left(2c_{1}^{2}-3c_{1}c_{2}+c_{2}^{2}+4k_{1}-k_{2}\right),\label{eq:alpha/beta/gamma-1}
\end{eqnarray}
with
\begin{equation}
D=\left(c_{1}^{2}-c_{1}c_{2}+k_{2}\right)\left(4c_{1}^{2}k_{2}-8c_{1}c_{2}k_{1}-2c_{1}c_{2}k_{2}+4c_{2}^{2}k_{1}+16k_{1}^{2}-8k_{1}k_{2}+k_{2}^{2}\right).\label{eq:Ddef-1}
\end{equation}
With these coefficients, substitution of \eqref{eq:g1-1} into the
first equation of system \eqref{eq:2DOF example} gives the exact
reduced system on the slow SSM, up to cubic order, in the form 
\[
\ddot{x}+\left[c_{1}+\left(\mu_{1}+a\beta\right)x^{2}\right]\dot{x}+\left(k_{1}+a\gamma\dot{x}^{2}\right)x+\left(b+a\alpha\right)x^{3}+\mathcal{O}\left(4\right)=0.
\]
 Substitution of the formulas \eqref{eq:alpha/beta/gamma-1} into
this last equation gives the final form \eqref{eq:exact reduced example}
of the exact reduced model on the SSM.

\section{Appendix: Details for Section \ref{subsec:Two-slow-degrees} \label{sec:Appendix:-Details-for-fast-slow-slow system}}

For the parameter range described by the scalings \eqref{eq:rescalings two slow},
we take the $\epsilon\to0$ limit in the expressions for $P_{1}$
and $P_{2}$ in \eqref{eq:P1_slow_slow}-\eqref{eq:P2_slow_slow}.
We then obtain

\begin{align*}
P_{1}\left(x,v,\eta,w,t;0\right) & =\left[\begin{array}{c}
p_{1}^{0}\\
p_{2}^{0}
\end{array}\right],
\end{align*}

\begin{align}
p_{1}^{0} & =\frac{1+\beta}{\delta^{2}}\biggl(-\mu_{p}v_{\gamma}-\delta^{2}\sin x_{\gamma}+\delta^{2}G_{p}(t)\\
 & +\frac{\delta\sin x_{\gamma}}{1+\beta}\bigl[\beta\delta\cos x_{\gamma}v_{\gamma}^{2}-\mu_{h}w-\Omega_{h}^{2}\eta-\alpha_{h}\eta^{3}+(1+\beta)\delta+F_{h}(t)\delta-F_{p}(t)\delta\sin x_{\gamma}\bigr]\biggr)\nonumber \\
 & -\frac{\phi}{\delta}\cos x_{\gamma}\left(\beta\frac{\delta}{\phi}\sin x_{\gamma}v_{\gamma}^{2}-\mu_{d}v_{d}-\Omega_{d}^{2}\left(1+x_{d}\right)Q^{0}(x_{d})+F_{d}(t)\frac{\delta}{\phi}+F_{p}(t)\frac{\delta}{\phi}\cos x_{\gamma}\right),\nonumber \\
p_{2}^{0} & =-\frac{\beta}{\delta\phi}\cos x_{\gamma}\biggl(-\mu_{p}v_{\gamma}-\delta^{2}\sin x_{\gamma}+\delta^{2}G_{p}(t)\\
 & +\frac{\delta\sin x_{\gamma}}{1+\beta}\bigl[\beta\delta\cos x_{\gamma}v_{\gamma}^{2}-\mu_{h}w-\Omega_{h}^{2}\eta-\alpha_{h}\eta^{3}+(1+\beta)\delta+F_{h}(t)\delta-F_{p}(t)\delta\sin x_{\gamma}\bigr]\biggr)\nonumber \\
 & +\frac{1}{1+\beta}\left(1+\beta\cos^{2}x_{\gamma}\right)\left(\beta\frac{\delta}{\phi}\sin x_{\gamma}v_{\gamma}^{2}-\mu_{d}v_{d}-\Omega_{d}^{2}\left(1+x_{d}\right)Q^{0}(x_{d})+F_{d}(t)\frac{\delta}{\phi}+F_{p}(t)\frac{\delta}{\phi}\cos x_{\gamma}\right),\nonumber 
\end{align}

\begin{eqnarray*}
P_{2}\left(x,v,\eta,w,t;0\right) & = & \left(\frac{1+\beta\sin^{2}x_{\gamma}}{1+\beta}\right)\biggl(\beta\delta\cos x_{\gamma}v_{\gamma}^{2}-\mu_{h}w-\Omega_{h}^{2}\eta-\alpha_{h}\eta^{3}+(1+\beta)\delta+F_{h}(t)\delta-F_{p}(t)\delta\sin x_{\gamma}\\
 &  & +\frac{\left(1+\beta\right)\beta\sin x_{\gamma}}{\delta\left(1+\beta\sin^{2}x_{\gamma}\right)}\left[-\mu_{p}v_{\gamma}-\delta^{2}\sin x_{\gamma}+\delta^{2}G_{p}(t)\right]\\
 &  & -\frac{\beta\phi\sin x_{\gamma}\cos x_{\gamma}}{1+\beta\sin^{2}x_{\gamma}}\biggl[\beta\frac{\delta}{\phi}\sin x_{\gamma}v_{\gamma}^{2}-\mu_{d}v_{d}-\Omega_{d}^{2}\left(1+x_{d}\right)Q^{0}(x_{d})+F_{d}(t)\frac{\delta}{\phi}+F_{p}(t)\frac{\delta}{\phi}\cos x_{\gamma}\biggr]\biggr),\\
\end{eqnarray*}
where $Q^{0}(x_{d})$ is defined as

\[
Q^{0}(x_{d})=\left(1-\frac{1}{1+x_{d}}\right),\quad x_{d}>-1.
\]
We observe that both $P_{1}$ and $P_{2}$ continue to be smooth in
$\epsilon$ at the $\epsilon=0$ limit, thereby satisfying assumption
(A1). 

For the critical manifold defined through the relationship $\eta=G_{0}(x,v,t)$
in assumption (A2), we have the equation

\[
P_{2}\left(x,v,\eta,0,t;0\right)=0\quad\Longleftrightarrow\quad\Omega_{h}^{2}\eta+\alpha_{h}\eta^{3}=T(x,v,t),
\]
where

\begin{eqnarray*}
T(x,v,t) & = & \beta\delta\cos x_{\gamma}v_{\gamma}^{2}+(1+\beta)\delta+F_{h}(t)\delta-F_{p}(t)\delta\sin x_{\gamma}\\
 &  & +\frac{\left(1+\beta\right)\beta\sin x_{\gamma}}{\delta\left(1+\beta\sin^{2}x_{\gamma}\right)}\left[-\mu_{p}v_{\gamma}-\delta^{2}\sin x_{\gamma}+\delta^{2}G_{p}(t)\right]\\
 &  & -\frac{\beta\phi\sin x_{\gamma}\cos x_{\gamma}}{1+\beta\sin^{2}x_{\gamma}}\biggl[\beta\frac{\delta}{\phi}\sin x_{\gamma}v_{\gamma}^{2}-\mu_{d}v_{d}-\Omega_{d}^{2}\left(1+x_{d}\right)Q^{0}(x_{d})+F_{d}(t)\frac{\delta}{\phi}+F_{p}(t)\frac{\delta}{\phi}\cos x_{\gamma}\biggr].
\end{eqnarray*}
Using the cubic formula, the real root of this equation can be expressed
explicitly as
\[
\eta=G_{0}(x,v,t)=\sqrt[3]{\frac{T(x,v,t)}{2\alpha_{h}}+\sqrt{\frac{T^{2}(x,v,t)}{4\alpha_{h}^{2}}+\frac{\Omega_{h}^{6}}{27\alpha_{h}^{3}}}}-\sqrt[3]{-\frac{T(x,v,t)}{2\alpha}+\sqrt{\frac{T^{2}(x,v,t)}{4\alpha_{h}^{2}}+\frac{\Omega_{h}^{6}}{27\alpha_{h}^{3}}}},
\]
assuming that $\Omega_{h}^{2}$ and $\alpha_{h}$ are greater than
zero. 

The oscillatory system \eqref{eq:associated linear system} determining
the stability of the critical manifold takes the specific form
\begin{eqnarray}
A(x,v,t) & = & -\partial_{w}P_{2}\left(x,v,G_{0}(x,v,t),0,t;0\right)=\left(\frac{1+\beta\sin^{2}x_{\gamma}}{1+\beta}\right)\mu_{h},\label{eq:ABdef-1-2-1}\\
B(x,v,t) & = & -\partial_{\eta}P_{2}\left(x,v,G_{0}(x,v,t),0,t;0\right)=\left(\frac{1+\beta\sin^{2}x_{\gamma}}{1+\beta}\right)\left(\Omega_{h}^{2}+3\alpha_{h}G_{0}^{2}(x,v,t)\right).
\end{eqnarray}
The equilibrium solution of the unforced linear oscillatory system
\eqref{eq:associated linear system} is, therefore, always asymptotically
stable, given that
\[
\mu_{h}>0,\quad\beta>0,\quad\Omega_{h}^{2}>0,\quad\alpha_{h}>0.
\]

We conclude that assumptions (A1)-(A3) hold, and hence a global reduced-order
model \eqref{eq:reduced model} exists over the slow variables in
the specific form

\begin{align*}
\ddot{x}= & \left[\begin{array}{l}
\frac{1+\beta}{\delta^{2}\left(1+\beta\sin^{2}x_{\gamma}\right)}\mathcal{A}-\frac{\phi\cos x_{\gamma}}{\delta\left(1+\beta\sin^{2}x_{\gamma}\right)}\mathcal{B}\\
\frac{1}{1+\beta\sin^{2}x_{\gamma}}\mathcal{B}-\frac{\beta\cos x_{\gamma}}{\phi\delta\left(1+\beta\sin^{2}x_{\gamma}\right)}\mathcal{A}
\end{array}\right]+\mathcal{O}(\epsilon),
\end{align*}
where

\begin{eqnarray*}
\mathcal{A}(x_{\gamma},\dot{x}_{\gamma}) & = & -\mu_{p}\dot{x}_{\gamma}-\delta^{2}\sin x_{\gamma}+\delta^{2}G_{p}(t),\\
\mathcal{B}(x_{\gamma},x_{d},\dot{x}_{\gamma}) & = & \beta\frac{\delta}{\phi}\sin x_{\gamma}\dot{x}_{\gamma}^{2}-\mu_{d}\dot{x}_{\gamma}-\Omega_{d}^{2}x_{d}+F_{d}(t)\frac{\delta}{\phi}+F_{p}(t)\frac{\delta}{\phi}\cos x_{\gamma}.
\end{eqnarray*}
Scaling back to the original time and substituting the physical parameters
back into the non-dimensionalized equations, we obtain that the exact
reduced-order model on the slow manifold of the form \eqref{eq:red_mod_1}-\eqref{eq:red_mod_2}

\section{Appendix: Details for Section \ref{subsec:One-slow-degree} .\label{sec:Appendix:-Details-for-fast-fast-slow system}}

Here we verify assumptions (A1)-(A3) in detail for the fast-fast-slow
setting treated in Section \ref{subsec:One-slow-degree}. To make
the horizontal spring stiff, we choose its length as $D=L$, so that
the original equations of motion \eqref{eq:spring-pendulum} now become

\begin{eqnarray*}
ml^{2}\ddot{\gamma}-ml\sin\gamma\ddot{h}+ml\cos\gamma\ddot{d}+c_{p}\dot{\gamma}+mgl\sin\gamma & = & f_{p}(t)l,\\
(M+m)\ddot{h}-ml\sin\gamma\ddot{\gamma}-ml\cos\gamma\dot{\gamma}^{2}+C_{h}\dot{h}+K_{h}h+K_{d}Q(d,h)h+\Gamma_{h}h^{3} & = & (M+m)g+f_{h}(t)-f_{p}(t)\sin\gamma,\\
(M+m)\ddot{d}+ml\cos\gamma\ddot{\gamma}-ml\sin\gamma\dot{\gamma}^{2}+C_{d}\dot{d}+K_{d}\left(L+d\right)Q(d,h) & = & f_{d}(t)+f_{p}(t)\cos\gamma,
\end{eqnarray*}
with 

\begin{equation}
Q(d,h)=\left(1-\frac{L}{\sqrt{\left(L+d\right)^{2}+h^{2}}}\right).
\end{equation}

The linearized oscillation frequencies of the uncoupled springs and
pendulum remain the same as in (\ref{eq:lineigen-1-1}). We adopt
the same scaling as in section \ref{subsec:Two-slow-degrees}, except
that we now scale the $d$ coordinate with the unstretched length
$L$ of the vertical spring. Denoting differentiation with respect
to the new time $\tilde{t}$ still by a dot, then dropping all the
tildes, we obtain the non-dimensionalized equations of motions

\begin{eqnarray*}
\Delta^{2}\ddot{\gamma}-\Delta\sin\gamma\ddot{h}+\Delta\cos\gamma\ddot{d}+\pi_{p}\dot{\gamma}+\Delta^{2}\sin\gamma & = & \Delta^{2}G_{p}(t),\\
(1+\beta)\ddot{h}-\beta\Delta\sin\gamma\ddot{\gamma}-\beta\Delta\cos\gamma\dot{\gamma}^{2}+\pi_{h}\dot{h}+q_{h}h+q_{d}hQ(d,h)+a_{h}h^{3} & = & (1+\beta)\Delta+F_{h}(t)\Delta-F_{p}(t)\Delta\sin\gamma,\\
(1+\beta)\ddot{d}+\beta\Delta\cos\gamma\ddot{\gamma}-\beta\Delta\sin\gamma\dot{\gamma}^{2}+\pi_{d}\dot{d}+q_{d}\left(1+d\right)Q(d,h) & = & F_{d}(t)\Delta+F_{p}(t)\Delta\cos\gamma.
\end{eqnarray*}

In the notation used for system \eqref{eq:systtem00}, we now have

\begin{eqnarray*}
M(q,t;\epsilon) & = & \left(\begin{array}{ccc}
\Delta^{2} & -\Delta\sin x & \Delta\cos x\\
-\beta\Delta\sin x & 1+\beta & 0\\
\beta\Delta\cos x & 0 & 1+\beta
\end{array}\right),\\
F(q,\dot{q},t;\epsilon) & = & \left(\begin{array}{l}
-\pi_{p}\dot{x}-\Delta^{2}\sin x+\Delta^{2}G_{p}(t)\vspace{3pt}\\
\beta\Delta\cos x\dot{x}^{2}-\pi_{h}\dot{y}_{h}-q_{h}\epsilon\frac{y_{h}}{\epsilon}-q_{d}\epsilon\frac{y_{h}}{\epsilon}Q(\frac{y_{d}}{\epsilon},\frac{y_{h}}{\epsilon})-a_{h}\epsilon^{3}\left(\frac{y_{h}}{\epsilon}\right)^{3}\\
+(1+\beta)\Delta+F_{h}(t)\Delta-F_{p}(t)\Delta\sin x\vspace{5pt}\\
\beta\Delta\sin x\dot{x}^{2}-\pi_{d}\dot{y}_{d}-q_{d}\left(1+\epsilon\frac{y_{d}}{\epsilon}\right)Q(\frac{y_{d}}{\epsilon},\frac{y_{h}}{\epsilon})+F_{d}(t)\Delta+F_{p}(t)\Delta\cos x
\end{array}\right),
\end{eqnarray*}
with the parameter $\epsilon>0$ yet to be determined based on the
assumptions of the SFD approach. Note that the mass-matrix above is
not symmetric due to the scalings we have employed, but it is nevertheless
nonsingular, as we generally assume in this paper.

With the above quantities at hand, we obtain the modified mass matrices
$M_{i}$ and the forcing terms $Q_{i}$ defined in \eqref{eq:M_i and Q_i def}
in the specific form
\begin{eqnarray*}
M_{1} & = & M_{11}-M_{12}M_{22}^{-1}M_{21}=\frac{\Delta^{2}}{1+\beta},\\
M_{2} & = & M_{22}-M_{21}M_{11}^{-1}M_{12}=\left(\begin{array}{cc}
1+\beta\cos^{2}x & \beta\sin x\cos x\\
\beta\sin x\cos x & 1+\beta\sin^{2}x
\end{array}\right),\\
Q_{1} & = & F_{1}-M_{12}M_{22}^{-1}F_{2}=-\pi_{p}\dot{x}-\Delta^{2}\sin x+\Delta^{2}G_{p}(t)+\frac{\Delta}{1+\beta}\sin x\biggl[\beta\Delta\cos x\dot{x}^{2}-\pi_{h}\dot{y}_{h}\\
 &  & -q_{h}\epsilon\frac{y_{h}}{\epsilon}-q_{d}\epsilon\frac{y_{h}}{\epsilon}Q(\frac{y_{d}}{\epsilon},\frac{y_{h}}{\epsilon})-a_{h}\epsilon^{3}\left(\frac{y_{h}}{\epsilon}\right)^{3}+(1+\beta)\Delta+F_{h}(t)\Delta-F_{p}(t)\Delta\sin x\biggr]\\
 &  & -\frac{\Delta}{1+\beta}\cos x\biggl[\beta\Delta\sin x\dot{x}^{2}-\pi_{d}\dot{y}_{d}-q_{d}\left(1+\epsilon\frac{y_{d}}{\epsilon}\right)Q(\frac{y_{d}}{\epsilon},\frac{y_{h}}{\epsilon})+F_{d}(t)\Delta+F_{p}(t)\Delta\cos x\biggr],\\
Q_{2} & = & F_{2}-M_{21}M_{11}^{-1}F_{1}\\
 & = & \left[\begin{array}{l}
\beta\Delta\cos x\dot{x}^{2}-\pi_{h}\dot{y}_{h}-q_{h}\epsilon\frac{y_{h}}{\epsilon}-q_{d}\epsilon\frac{y_{h}}{\epsilon}Q(\frac{y_{d}}{\epsilon},\frac{y_{h}}{\epsilon})-a_{h}\epsilon^{3}\left(\frac{y_{h}}{\epsilon}\right)^{3}+(1+\beta)\Delta\\
+F_{h}(t)\Delta-F_{p}(t)\Delta\sin x-\frac{\beta}{\Delta}\pi_{p}\sin x\dot{x}-\beta\Delta\sin^{2}x+\beta\Delta\sin xG_{p}(t)\vspace{5pt}\\
\beta\Delta\sin x\dot{x}^{2}-\pi_{d}\dot{y}_{d}-q_{d}\left(1+\epsilon\frac{y_{d}}{\epsilon}\right)Q(\frac{y_{d}}{\epsilon},\frac{y_{h}}{\epsilon})+F_{d}(t)\Delta+F_{p}(t)\Delta\cos x\\
+\frac{\beta}{\Delta}\pi_{p}\cos x\dot{x}+\beta\Delta\sin x\cos x-\beta\Delta\cos xG_{p}(t)
\end{array}\right].
\end{eqnarray*}
We therefore obtain
\begin{eqnarray*}
P_{1}\left(x,v,\eta,w,t;\epsilon\right) & = & \frac{1+\beta}{\Delta^{2}}\Biggl[-\pi_{p}v-\Delta^{2}\sin x+\Delta^{2}G_{p}(t)+\frac{\Delta}{1+\beta}\sin x\biggl[\beta\Delta\cos xv^{2}-\pi_{h}w_{h}\\
 &  & -q_{h}\epsilon\eta_{h}-q_{d}\epsilon\eta_{h}Q(\eta_{d},\eta_{h})-a_{h}\epsilon^{3}\eta_{h}^{3}+(1+\beta)\Delta+F_{h}(t)\Delta-F_{p}(t)\Delta\sin x\biggr]\\
 &  & -\frac{\Delta}{1+\beta}\cos x\biggl[\beta\Delta\sin xv^{2}-\pi_{d}w_{d}-q_{d}\left(1+\epsilon\eta_{d}\right)Q(\eta_{d},\eta_{h})+F_{d}(t)\Delta+F_{p}(t)\Delta\cos x\biggr]\Biggr],\\
P_{2}\left(x,v,\eta,w,t;\epsilon\right) & = & \epsilon M_{2}^{-1}\left[\begin{array}{l}
\beta\Delta\cos xv^{2}-\pi_{h}w_{h}-q_{h}\epsilon\eta_{h}-q_{d}\epsilon\eta_{h}Q(\eta_{d},\eta_{h})-a_{h}\epsilon^{3}\eta_{h}^{3}+(1+\beta)\Delta\\
+F_{h}(t)\Delta-F_{p}(t)\Delta\sin x-\frac{\beta}{\Delta}\pi_{p}\sin xv-\beta\Delta\sin^{2}x+\beta\Delta\sin xG_{p}(t)\\
\\
\beta\Delta\sin xv^{2}-\pi_{d}w_{d}-q_{d}\left(1+\epsilon\eta_{d}\right)Q(\eta_{d},\eta_{h})+F_{d}(t)\Delta+F_{p}(t)\Delta\cos x\\
+\frac{\beta}{\Delta}\pi_{p}\cos xv+\beta\Delta\sin x\cos x-\beta\Delta\cos xG_{p}(t)
\end{array}\right],
\end{eqnarray*}
where $M_{2}^{-1}$ is equal to 

\[
M_{2}^{-1}=\frac{1}{1+\beta}\left[\begin{array}{cc}
1+\beta\sin^{2}x & -\beta\sin x\cos x\\
-\beta\sin x\cos x & 1+\beta\cos^{2}x
\end{array}\right].
\]

Recall that $\epsilon>0$ has been a completely arbitrary small parameter
so far. We now need to define $\epsilon$ in a way that assumptions
(A1)-(A3) are satisfied. Since at present we have $\lim_{\epsilon\to0}P_{2}\left(x,v,\eta,w,t;\epsilon\right)\equiv0,$
these assumptions will not hold. We can only satisfy (A1)-(A3) by
making the system parameters appropriate functions of $\epsilon.$ 

With the parameter choices listed in \eqref{eq:rescaling for fast-fast-slow},
we have

\begin{eqnarray*}
P_{1}\left(x,v,\eta,w,t;\epsilon\right) & = & \frac{1+\beta}{\delta^{2}}\Biggl[-\mu_{p}v-\delta^{2}\sin x+\delta^{2}G_{p}(t)+\frac{\delta}{1+\beta}\sin x\biggl[\beta\delta\cos xv^{2}-\mu_{h}w_{h}\\
 &  & -\Omega_{h}^{2}\eta_{h}-\Omega_{d}^{2}\eta_{h}Q(\eta_{d},\eta_{h})-\alpha_{h}\eta_{h}^{3}+(1+\beta)\delta+F_{h}(t)\delta-F_{p}(t)\delta\sin x\biggr]\\
 &  & -\frac{\delta}{1+\beta}\cos x\biggl[\beta\delta\sin xv^{2}-\mu_{d}w_{d}-\frac{\Omega_{d}^{2}}{\epsilon}\left(1+\epsilon\eta_{d}\right)Q(\eta_{d},\eta_{h})+F_{d}(t)\delta+F_{p}(t)\delta\cos x\biggr]\Biggr],\\
P_{2}\left(x,v,\eta,w,t;\epsilon\right) & = & M_{2}^{-1}\left[\begin{array}{l}
\beta\delta\cos xv^{2}-\mu_{h}w_{h}-\Omega_{h}^{2}\eta_{h}-\Omega_{d}^{2}\eta_{h}Q(\eta_{d},\eta_{h})-\alpha_{h}\eta_{h}^{3}+(1+\beta)\delta\\
+F_{h}(t)\delta-F_{p}(t)\delta\sin x-\frac{\beta}{\delta}\mu_{p}\sin xv-\beta\delta\sin^{2}x+\beta\delta\sin xG_{p}(t)\\
\\
\beta\delta\sin xv^{2}-\mu_{d}w_{d}-\frac{\Omega_{d}^{2}}{\epsilon}\left(1+\epsilon\eta_{d}\right)Q(\eta_{d},\eta_{h})+F_{d}(t)\delta+F_{p}(t)\delta\cos x\\
+\frac{\beta}{\delta}\mu_{p}\cos xv+\beta\delta\sin x\cos x-\beta\delta\cos xG_{p}(t)
\end{array}\right].
\end{eqnarray*}
where $M_{2}^{-1}$ remains unchanged.

Noting that

\begin{align*}
 & \lim_{\epsilon\rightarrow0}\frac{\Omega_{d}^{2}}{\epsilon}\left(1+\epsilon\eta_{d}\right)Q(\eta_{d},\eta_{h})\\
 & =\lim_{\epsilon\rightarrow0}\frac{\Omega_{d}^{2}}{\epsilon}\left(1+\epsilon\eta_{d}\right)\left(1-\frac{1}{\sqrt{\left(1+\epsilon\eta_{d}\right)^{2}+\left(\epsilon\eta_{h}\right)^{2}}}\right)\\
 & =\lim_{\epsilon\rightarrow0}\frac{\Omega_{d}^{2}\left(1+\epsilon\eta_{d}\right)\left(\sqrt{\left(1+\epsilon\eta_{d}\right)^{2}+\left(\epsilon\eta_{h}\right)^{2}}-1\right)}{\epsilon\sqrt{\left(1+\epsilon\eta_{d}\right)^{2}+\left(\epsilon\eta_{h}\right)^{2}}}\\
 & =\lim_{\epsilon\rightarrow0}\frac{f(\epsilon)}{g(\epsilon)}=\lim_{\epsilon\rightarrow0}\frac{\partial_{\epsilon}f(\epsilon)}{\partial_{\epsilon}g(\epsilon)}\\
 & =\lim_{\epsilon\rightarrow0}\frac{\Omega_{d}^{2}\eta_{d}\left(\sqrt{\left(1+\epsilon\eta_{d}\right){}^{2}+\left(\epsilon\eta_{h}\right)^{2}}-1\right)+\Omega_{d}^{2}\left(1+\epsilon\eta_{d}\right)\left(\left((1+\epsilon\eta_{d})\eta_{d}+\epsilon\eta_{h}^{2}\right)\left(\left(1+\epsilon\eta_{d}\right)^{2}+\left(\epsilon\eta_{h}\right)^{2}\right)^{-\frac{1}{2}}\right)}{\sqrt{\left(1+\epsilon\eta_{d}\right)^{2}+\left(\epsilon\eta_{h}\right)^{2}}+\epsilon\left((1+\epsilon\eta_{d})\eta_{d}+\epsilon\eta_{h}^{2}\right)\left(\left(1+\epsilon\eta_{d}\right)^{2}+\left(\epsilon\eta_{h}\right)^{2}\right)^{-\frac{1}{2}}}\\
 & =\Omega_{d}^{2}\eta_{d},
\end{align*}
we conclude that both $P_{1}$ and $P_{2}$ continue to be smooth
in $\epsilon$ at the $\epsilon=0$ limit, thereby satisfying assumption
(A1). 

For the critical manifold defined through the relationship $\eta=G_{0}(x,v,t)$
in assumption (A2), we have the equations 

\begin{eqnarray*}
P_{2}\left(x,v,\eta,0,t;0\right) & = & M_{2}^{-1}\left[\begin{array}{l}
\beta\delta\cos xv^{2}-\Omega_{h}^{2}\eta_{h}-\alpha_{h}\eta_{h}^{3}+(1+\beta)\delta+F_{h}(t)\delta\\
-F_{p}(t)\delta\sin x-\frac{\beta}{\delta}\mu_{p}\sin xv-\beta\delta\sin^{2}x+\beta\delta\sin xG_{p}(t)\\
\\
\beta\delta\sin xv^{2}-\Omega_{d}^{2}\eta_{d}+F_{d}(t)\delta+F_{p}(t)\delta\cos x\\
+\frac{\beta}{\delta}\mu_{p}\cos xv+\beta\delta\sin x\cos x-\beta\delta\cos xG_{p}(t)
\end{array}\right]=\left[\begin{array}{c}
0\\
0
\end{array}\right].
\end{eqnarray*}
Since $M_{2}^{-1}$ is invertible, the critical manifold can be found
by solving the following equations for $\eta_{h}$ and $\eta_{d}:$

\begin{gather*}
\Omega_{h}^{2}\eta_{h}+\alpha_{h}\eta_{h}^{3}=T_{h}(x,v,t)=\beta\delta\cos xv^{2}+(1+\beta)\delta+F_{h}(t)\delta-F_{p}(t)\delta\sin x-\frac{\beta}{\delta}\mu_{p}\sin xv-\beta\delta\sin^{2}x+\beta\delta\sin xG_{p}(t),\\
\Omega_{d}^{2}\eta_{d}=T_{d}(x,v,t)=\beta\delta\sin xv^{2}+F_{d}(t)\delta+F_{p}(t)\delta\cos x+\frac{\beta}{\delta}\mu_{p}\cos xv+\beta\delta\sin x\cos x-\beta\delta\cos xG_{p}(t).
\end{gather*}
The real roots of these two equations can be expressed explicitly
as
\[
\eta_{h}=\sqrt[3]{\frac{T_{h}(x,v,t)}{2\alpha_{h}}+\sqrt{\frac{T_{h}^{2}(x,v,t)}{4\alpha_{h}^{2}}+\frac{\Omega_{h}^{6}}{27\alpha_{h}^{3}}}}-\sqrt[3]{-\frac{T_{h}(x,v,t)}{2\alpha}+\sqrt{\frac{T_{h}^{2}(x,v,t)}{4\alpha_{h}^{2}}+\frac{\Omega_{h}^{6}}{27\alpha_{h}^{3}}}},
\]

\[
\eta_{d}=\frac{T_{d}(x,v,t)}{\Omega_{d}^{2}},
\]
assuming that $\Omega_{h}^{2}$, $\Omega_{d}^{2}$ and $\alpha_{h}$
are greater than zero. The stability of this critical manifold is
determined by the associated oscillatory system \eqref{eq:associated linear system},
whose coefficient matrices now take the specific form
\begin{eqnarray*}
A(x,v,t) & = & -\partial_{w}P_{2}\left(x,v,G_{0}(x,v,t),0,t;0\right)\\
 & = & \frac{1}{1+\beta}\left[\begin{array}{cc}
1+\beta\sin^{2}x & -\beta\sin x\cos x\\
-\beta\sin x\cos x & 1+\beta\cos^{2}x
\end{array}\right]\left[\begin{array}{cc}
\mu_{h} & 0\\
0 & \mu_{d}
\end{array}\right],\\
B(x,v,t) & = & -\partial_{\eta}P_{2}\left(x,v,G_{0}(x,v,t),0,t;0\right)\\
 & = & \frac{1}{1+\beta}\left[\begin{array}{cc}
1+\beta\sin^{2}x & -\beta\sin x\cos x\\
-\beta\sin x\cos x & 1+\beta\cos^{2}x
\end{array}\right]\left[\begin{array}{cc}
\Omega_{h}^{2}+3\alpha_{h}\eta_{h}^{2} & 0\\
0 & \Omega_{d}^{2}
\end{array}\right].
\end{eqnarray*}
Consequently, the equilibrium solution of the unforced linear oscillatory
system \eqref{eq:associated linear system} is always asymptotically
stable, given that
\[
\mu_{h}>0,\quad\mu_{d}>0\quad\beta>0,\quad\Omega_{h}^{2}>0,\quad\Omega_{d}^{2}>0,\quad\alpha_{h}>0.
\]
We conclude that assumptions (A1)-(A3) hold, and hence a global reduced-order
model exists over the flexible variables $(x,v,t)\in$$\mathcal{D}_{0}=\mathbb{R}\times\mathbb{R}\times S^{1}$.

\end{document}